\documentclass{amsart}
\usepackage{amssymb, graphicx}

\theoremstyle{plain}
\newtheorem{prop}{Proposition}[section]
\newtheorem{coro}[prop]{Corollary}

\newtheorem{lemm}[prop]{Lemma}
\newtheorem{ques}[prop]{Question}

\theoremstyle{definition}
\newtheorem{defi}[prop]{Definition}

\newtheorem{exam}[prop]{Example}
\newtheorem{rema}[prop]{Remark}
\newtheorem*{remn}{Remark on notation}

\def\Reff#1; #2; #3; #4; #5; #6; #7\par{%
\bibitem{#1} #2, {\it #3}, #4 {\bf #5} (#6) #7}

\def\Ref#1; #2; #3; #4\par{%
\bibitem{#1} #2, {\it #3}, #4}

\catcode`\¥=\active\def¥{\relax}

\def\aa#1{a_{#1}}
\def\aaa#1{a_{#1}\inv}
\def\aaaa#1{a_{#1}\pmo}
\def\act{\mathbin{\scriptscriptstyle\bullet}}
\newcommand{\Aut}{{\rm Aut}}
\def\auto(#1){\widehat{#1}}

\newcommand{\MyIndex}{\bullet}
\newcommand{\BBs}{\widetilde{B}_\MyIndex}
\newcommand{\BBsp}{\BBs^\smp}
\def\Bi{B_\infty}
\def\Bisp{B_\infty^{sp}}
\def\Bip{B_\infty^\smp}
\newcommand{\Bs}{B_\MyIndex}
\newcommand{\Bsp}{B_\MyIndex^\smp}
\newcommand{\bx}{\beta}

\newcommand{\bxx}{\bx'}
\newcommand{\bxxx}{\bx''}
\newcommand{\by}{\gamma}

\newcommand{\bz}{\alpha}

\newcommand{\Can}{\mathbf{K}}
\newcommand{\cc}[2]{c_{{#1}, {#2}}}
\newcommand{\CCC}{\mathbf{C}}
\newcommand{\cl}[1]{\overline#1}
\newcommand{\cL}{/}
\newcommand{\Col}{\mathrm{Col}}
\newcommand{\Cone}{C}
\newcommand{\cR}{\backslash}
\newcommand{\ct}{\boldsymbol{t}}
\newcommand{\cti}{\ct_1}
\newcommand{\ctii}{\ct_2}
\newcommand{\ctt}{\ct'}
\newcommand{\ctti}{\ct'_1}
\newcommand{\cttii}{\ct'_2}

\def\dbl{\mathrm{db}}		
\newcommand{\dd}[1]{{#1}^{\scriptscriptstyle\mathtt{\#}}}
\newcommand{\ddd}[1]{(#1)^{\scriptscriptstyle\mathtt{\#}}}
\newcommand{\dddsmall}[1]{\scriptstyle(#1)^{\!\scriptscriptstyle\mathtt{\#}}}
\newcommand{\ddsmall}[1]{\scriptstyle(#1)}
\newcommand{\DD}{\mathcal{D}}
\newcommand{\DDs}{\mathcal{D}_\MyIndex}
\newcommand{\Dec}{\mathrm{dec}}
\newcommand{\Dya}{\mathrm{Dyad}}

\def\e{\varepsilon}
\def\eg{{\it e.g.}}
\newcommand{\eps}{\epsilon}
\newcommand{\et}{\mathord{\bullet}}
\newcommand{\etwhite}{\mathord{\circ}}
\def\etc{{\it etc.}}
\newcommand{\ev}{\mathrm{ev}}
\newcommand{\Evv}{\ev^*}
\newcommand{\ex}{x}
\newcommand{\ey}{y}
\newcommand{\ez}{z}

\newcommand{\Fp}{F^\smp}

\newcommand{\FGs}{F_\MyIndex}
\newcommand{\fw}{w}
\newcommand{\fx}{f}
\newcommand{\fxi}{\fx_1}
\newcommand{\fxii}{\fx_2}
\newcommand{\fxx}{\fx'}
\newcommand{\fy}{g}
\newcommand{\fyy}{\fy'}
\newcommand{\fz}{h}

\let\g=\gamma
\let\ge=\geqslant
\def\Gr(#1;#2){\langle#1\,; #2\rangle}
\newcommand{\gx}{x}
\newcommand{\gy}{y}
\newcommand{\gz}{z}

\newcommand{\ha}{\widehat{a}}
\let\Hat=\widehat
\renewcommand{\hom}[1]{H(#1)}
\newcommand{\hs}{\widehat{\s}}

\def\ie{{\it i.e.}}
\let\ince=\subseteq
\newcommand{\inv}{^{-1}}
\newcommand{\invv}{^{\!-\!1}}
\newcommand{\Inv}{{}^{-1}}

\newcommand{\lB}{<_B}
\def\LD#1#2{{#1}\lceil{#2}\rceil}
\def\LDbar#1#2{{#1}\lfloor{#2}\rfloor}
\def\LDD#1#2{{#1}\lceil{#2}\rceil}
\def\le{\leqslant}
\newcommand{\lF}{<_F}
\newcommand{\lp}{<^{\scriptscriptstyle+}}

\newcommand{\MCG}{M\!C\!G}
\def\Mon(#1;#2){\langle#1\,; #2\rangle^\smp}

\newcommand{\nn}{n}
\newcommand{\nni}{\nn_1}
\newcommand{\nnii}{\nn_2}

\newcommand{\NNN}{\mathbf{N}}
\newcommand{\NNNs}{\NNN_\MyIndex}

\def\op{\relax}
\newcommand{\OP}{\circ}
\let\opp=\cdot

\newcommand{\PBi}{P\!B_\infty}
\newcommand{\PBs}{P\!B_\MyIndex}
\newcommand{\pmo}{^{\pm1}}
\let\pp=\ldots

\newcommand{\Pos}[1]{\mathrm{Pos}(#1)}
\newcommand{\pw}{w}
\newcommand{\pwi}{\pw_1}
\newcommand{\pwii}{\pw_2}
\newcommand{\pwiii}{\pw_3}
\newcommand{\pww}{\pw'}
\newcommand{\px}{x}

\newcommand{\pxx}{\px'}
\newcommand{\py}{y}

\newcommand{\pyy}{\py'}
\newcommand{\pz}{z}
\newcommand{\pzi}{\pz_1}
\newcommand{\pzii}{\pz_2}

\newcommand{\RR}{R}
\newcommand{\RRo}{R_0}
\newcommand{\RRs}{R_\MyIndex}
\newcommand{\red}{\mathrm{red}}

\def\revr#1{\curvearrowright_{#1}}
\newcommand{\revL}{\curvearrowleft}
\newcommand{\revR}{\curvearrowright}
\newcommand{\RRR}{\mathbf{R}}
\def\resp{{\it resp.\ }}

\let\s=\sigma
\newcommand{\SetOfCol}{S}
\renewcommand{\ss}[1]{\sigma_{#1}}
\def\sss#1{\sigma_#1\inv}
\def\ssss#1{\sigma_#1\pmo}

\newcommand{\SK}{S_\Can}
\def\smp{{\scriptscriptstyle +}}
\def\sp{\!+\!}
\newcommand{\Sym}{S}
\newcommand{\Symi}{\Sym_\infty}
\newcommand{\Syms}{\Sym_\MyIndex}

\def\tt{t}
\def\tti{\tt_1}
\def\ttii{\tt_2}

\def\ttt{\tt'}
\def\ttti{\ttt_1}
\def\tttii{\ttt_2}

\def\ttii{\tt_2}

\def\ttt{\tt'}
\def\ttti{\ttt_1}
\def\tttii{\ttt_2}

\def\uu{u}
\def\uui{u_1}

\newcommand{\uuu}{u'}

\def\vv{v}
\def\vvi{v_1}

\newcommand{\vvv}{v'}

\newcommand{\wa}{a}
\newcommand{\ws}{\sigma}
\newcommand{\wsa}{\sigma\text{,}a}
\newcommand{\ww}{w}
\newcommand{\wwi}{\ww_1}
\newcommand{\wwii}{\ww_2}

\newcommand{\www}{\ww'}

\def\xx{x}
\newcommand{\XX}{X}
\newcommand{\XXo}{\XX_0}
\newcommand{\XXa}{\aa*}
\newcommand{\XXopm}{\XXo^{\pm}}
\newcommand{\XXpm}{\XX^{\pm}}
\newcommand{\XXs}{\ss*}
\newcommand{\XXsa}{\XXa,\XXs}

\def\yy{y}

\def\zz{z}

\begin{document}

\author{Patrick DEHORNOY}
\address{Laboratoire de Math\'ematiques Nicolas
Oresme UMR 6139\\ Universit\'e de Caen,
14032~Caen, France}
\email{dehornoy@math.unicaen.fr}
\urladdr{//www.math.unicaen.fr/\!\!\!\hbox{$\sim$}dehornoy}

\title{The group of parenthesized braids}



\begin{abstract}
We investigate a group~$\Bs$ that includes Artin's
braid group~$\Bi$ and Thompson's group~$F$. The
elements of~$\Bs$ are represented by braids diagrams
in which the distances between the strands are not
uniform and, besides the usual crossing generators,
new rescaling operators shrink or strech the
distances between the strands. We prove that
$\Bs$ is a group of fractions, that it is orderable,
admits a non-trivial self-distributive structure, \ie,
one involving the law $x(yz)=(xy)(xz)$, it embeds in the
mapping class group of a sphere with a Cantor set of
punctures, and that Artin's representation of~$\Bi$
into the automorphisms of a free group extends to~$\Bs$.
\end{abstract}

\maketitle

The aim of this paper is to study a certain group,
denoted~$\Bs$, which includes both Artin's braid
group~$\Bi$ \cite{Bir, BuZ, Dgd} and Thompson's
group~$F$~\cite{Tho, McT, CFP}. The group~$\Bs$ is
generated by (the copies of)~$\Bi$ and~$F$,
and its seemingly rich and deep properties appear to be
a mixture of those of~$\Bi$ and~$F$. Here, starting
from a geometric approach in terms of parenthesized
braid diagrams, we give an explicit presentation
of~$\Bs$ that extends the standard presentations
of~$\Bi$ and~$F$, we prove that
$\Bs$ is a group of fractions, is an orderable group,
and embeds into the mapping class group of a sphere
with a Cantor set of punctures and into the
automorphisms of a free group. Besides its group
multiplication, $\Bs$ is also equipped with a second
binary operation satisfying the self-distributivity law
$x(yz)=(xy)(xz)$. We prove that every element of~$\Bs$
generates a free subsystem with respect to that second
operation---which shows that the self-distributive
structure of~$\Bs$ is highly non-trivial---and we
deduce canonical decompositions for the elements
of~$\Bs$. The self-distributive structure is
instrumental in proving most of the above results
about the group structure of~$\Bs$.

Here the elements of~$\Bs$ are seen as {\it
parenthesized braids}, which are braids in which the
distances between the strands are not uniform. An
ordinary braid diagram connects an initial sequence of
equidistant positions to a similar final sequence, as
for instance in
\begin{center}
\begin{picture}(20,10)(0, 0)
\put(0,0){\includegraphics{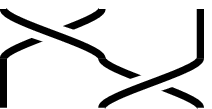}}
\end{picture}
\end{center}
where the initial and final set of positions can
be denoted~$\et\et\et$. A parenthesized braid diagram
connects a parenthesized sequence of positions to
another possibly different parenthezied sequence of
positions, the intuition being that grouped positions
are (infinitely) closer than ungrouped ones. An
example is
\begin{center}
\begin{picture}(20,20)(0, 0)
\put(0,0){\includegraphics{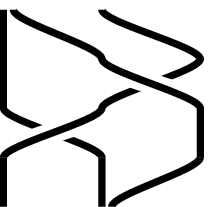}}
\end{picture}
\end{center}
where the initial positions are $(\et\et)\et$ and
the final positions are $\et(\et\et)$. Arranging such
objects into a group leads to introducing, besides the
usual braid generators~$\ss i$ that create crossings,
new rescaling generators~$\aa i$ that shrink the
distances between the strands in the vicinity of
position~$i$: as one can expect, the~$\ss i$'s
generate the copy of~$\Bi$, while the~$\aa i$'s generate
the copy of Thompson's group~$F$.

Parenthesized braids have been considered by D.\,Bar
Natan in~\cite{Bar1, Bar2} in connection with
Vassiliev's invariants of knots and the computation
of a Drinfeld associator. In these papers,
parenthesized braids, and more generally
parenthesized tangles, are studied as categories, and
the question of finding presentations is not addressed.

The realization of~$\Bs$ as a group of
parenthesized braids is not the only possible one,
and this group recently appeared in various
frameworks. In~\cite{Bri1, Bri2}, M.\,Brin
investigates a certain group~${BV}$ introduced as a 
torsion-free version of Thompson's group~$V$, and which
admits a subgroup~$\widehat{BV}$ that is isomorphic to~$\Bs$.
In~\cite{Dhb}, an independent approach leads to
introducing~$\Bs$ as the so-called geometry group for the
associativity law together with a twisted version of the
semi-commutativity law. All these approaches are more or
less equivalent, but we think that parenthesized braids
provide an especially intuitive and natural description.
Larger groups extending both braid groups and Thompson's
groups appear in~\cite{GrS, FuK, KaS}.

The current paper is self-contained in that it
requires no knowledge of the above mentioned papers
(by contrast, \cite{Dhb} resorts to results from
the current paper). As for results, the only
overlap with other papers is the result that $\Bs$ is
a group of fractions, which is established using
Zappa-Sz\'ep products of monoids in~\cite{Bri1}, while
we deduce it from general results involving the word
reversing technique. 

\begin{remn}
We follow the usual braid conventions: our
generators~$\ss i$ are numbered from~$1$, and the
product corresponds to an action on the right ($xy$
means $x$ followed by~$y$). For coherence, we adopt
a similar notation for Thompson's group~$F$, thus
shifting indices and reversing products: what we
denote~$\aa i$ is~$x_{i-1}\inv$ (or $X_{i-1}\inv$) in
the standard presentation of~$F$ \cite{CFP}.
An index of terms and notation is given at the end of the
paper.
\end{remn}

The author thanks Matthew Brin for helpful comments
and suggestions.

\section{Parenthesized braids}
\label{S:PBraids}

Throughout the paper, $\NNN$ denotes the set
of all positive integers ($0$ excluded). 

We construct a new group~$\Bs$ using the approach
that is standard for braids, namely starting with
isotopy classes of braid diagrams. The  difference
is that we consider diagrams in which the distances
between the endpoints of the strands need not be
uniform. Such sets of positions can be specified
using parenthesized expressions, like
$\et((\et\et)\et)$, where grouped positions are to
be seen as infinitely closed than the adjacent ones.
This principle is implemented by considering
positions that are indexed by finite sequence of integers.

The current construction of~$\Bs$ is exactly as simple
as that of~$\Bi$. Although making it precise requires
some notation, needed in particular in subsequent
proofs, the ideas should be clear, and many details can
be skipped.

\subsection{An intuitive description}

A braid diagram consists of curves that connect an
initial sequence of positions to a similar final
sequence of positions. In an ordinary braid diagram,
the positions are indexed by positive integers
\begin{center}
\begin{picture}(83,4)(0,0)
\put(10,-1){\includegraphics{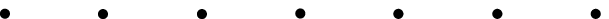}}
\put(9.5,2){$1$}
\put(19.5,2){$2$}
\put(77,0){\pp}
\end{picture}
\end{center}
and a generic diagram is obtained by stacking
elementary diagrams of the type
\begin{center}
\begin{picture}(83,11)(0,0)
\put(10,0){\includegraphics{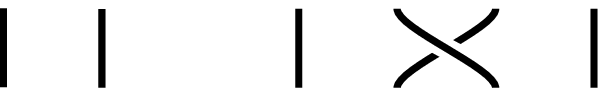}}
\put(0,3){$\ss i$~:}
\put(9.5,9){$1$}
\put(19.5,9){$2$}
\put(50,9){$i$}
\put(58,9){$i{+}1$}
\put(27,3){\pp}
\put(77,3){\pp}
\end{picture}
\end{center}
or their reflections in a horizontal mirror. 
 
\begin{figure}[htb]
\begin{picture}(120,22)(0, 0)
\put(10,0){\includegraphics{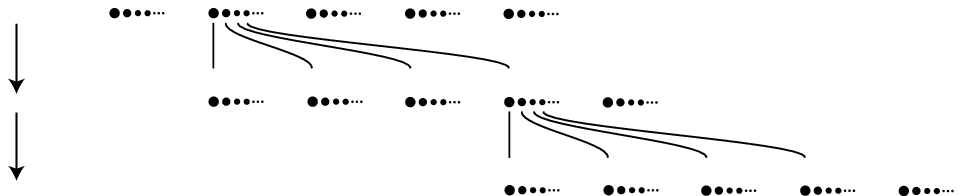}}
\put(21,20){$1$}
\put(31,20){$2$}
\put(41,20){$3$}
\put(31,11){$\scriptstyle 2$}
\put(40,11){$\scriptstyle 2\sp\eps$}
\put(50,11){$\scriptstyle 2\sp2\eps$}
\put(60,11){$\scriptstyle 2\sp3\eps$}
\put(60,2){$\scriptscriptstyle 2\sp3\eps$}
\put(70,2){$\scriptscriptstyle 2\sp3\eps\sp\eps^2$}
\put(80,2){$\scriptscriptstyle 2\sp3\eps\sp2\eps^2$}
\put(90,2){$\scriptscriptstyle 2\sp3\eps\sp3\eps^2$}
\put(0,5){$\times \eps\inv$}
\put(0,13){$\times \eps\inv$}
\put(70,18){\pp}
\put(80,9){\pp}
\put(110,0){\pp}
\end{picture}
\caption{\smaller The set of all positions
realized using infinitesimal distances}
\label{F:Positions}
\end{figure}

Here we consider braid diagrams in which the initial
and final positions need not be equidistant, but
instead the distances may be $1, \eps,
\eps^2, \pp$ with $\eps \ll 1$. This leads to
considering that, between the positions~$1$ and~$2$,
infinitely many new positions $1 + \eps$, $1 +
2\eps$, \pp are possible, and so on iteratively.
Thus $2 + 3\eps + \eps^2$ or $1 + \eps^3$ are
typical positions (Figure~\ref{F:Positions}).

Then, as in the case of ordinary braid diagrams, we
can consider generalized braid diagrams obtained by
stacking (finitely many) elementary crossing diagrams
\begin{center}
\begin{picture}(83,11)(0,0)
\put(10.5,0){\includegraphics{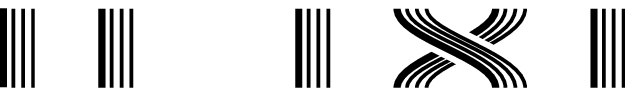}}
\put(0,3){$\ss i$~:}
\put(10,9){$1$}
\put(20,9){$2$}
\put(50,9){$i$}
\put(58,9){$i\sp1$}
\put(27,3){\pp}
\put(77,3){\pp}
\end{picture}
\end{center}
in which all strands near position~$i$ cross over
all strands near position~$i+1$, and rescaling
diagrams
\begin{center}
\begin{picture}(83,10)(0,0)
\put(10.5,0){\includegraphics{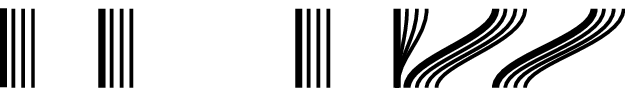}}
\put(0,3){$\aa i$~:}
\put(10,9){$1$}
\put(20,9){$2$}
\put(50,9){$i$}
\put(58,9){$i\sp1$}
\put(27,3){\pp}
\put(77,3){\pp}
\end{picture}
\end{center}
in which the strands near position~$i$ are shrinked
by a factor~$\eps$ and all strands on the right are
translated to fill the gaps. We also allow the mirror
images of the above diagrams. Our claim is that such
diagrams up to isotopy form a group, and this group is
the object we investigate in this paper.

Though intuitive, the previous informal description
is partly misleading in that it involves diagrams
with infinitely many strands. The objects we really
wish to consider are finite subdiagrams obtained by
restricting to a finite set of positions. In this
way, one exactly obtains the diagrams that are
arranged into a small category in~\cite{Bar1, Bar2},
the objects being the possible sets of
positions---which we shall see can be specified by
parenthesized expressions or, equivalently, finite
binary trees---and the morphisms being the isotopy
classes of braid diagrams.

A (minor) problem arises when we wish to make a
group out of the previous objects. In ordinary
braid diagrams, the initial and final positions
coincide, so, for each~$n$, concatenating $n$~strand
diagrams is always possible, which leads to the
braid group~$B_n$. In our extended framework,
concatenating two diagrams~$\DD_1, \DD_2$ is
possible only when the final set of positions
in~$\DD_1$ coincides with the initial set of
positions in~$\DD_2$, and an everywhere defined
product appears only when we consider infinite
completions, a situation similar to that of~$\Bi$:
to make a group out of all ordinary diagrams,
independently on the number of strands, one
embeds~$B_n$ into~$B_{n'}$ for $n < n'$ and the
elements of~$\Bi$ are then represented by infinite
diagrams.

\subsection{Sets of positions, parenthesized
expressions and trees}

For a more formal construction, we first define the
convenient sets of positions. Infinitesimal distances are
intuitive, but there Is no need to use them: the
infinitesimals we consider are polynomials in~$\eps$, and
the simplest solution is to index positions by polynomials,
\ie, by finite  sequences of nonnegative integers. To make
explicit geometric constructions easier, we also embed
positions into the unit interval using a dyadic expansion.

\begin{defi} \label{I:Position}
A finite sequence of nonnegative integers is called a {\it
position} if it does not begin or finish with~$0$. The set
of all positions is denoted by~$\NNNs$. For~$s$ a
position---or, more generally, any finite sequence of
nonnegative integers not beginning with~$0$---say $s = (i_1,
\pp, i_p)$, the {\it dyadic realization} of~$s$ is the
rational number~$\dd s$ with dyadic expansion
$0.1^{i_1-1}01^{i_2}0\pp 1^{i_p}$.
\end{defi}

\begin{figure}[htb]
\begin{picture}(130,24)(0, 0)
\put(0.5,3){\includegraphics{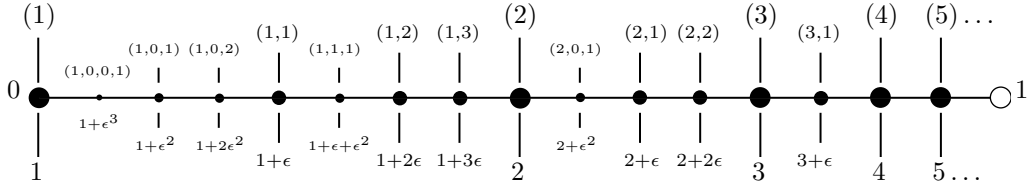}}
\put(1,0){$1$}
\put(65,0){$2$}
\put(97,0){$3$}
\put(113,0){$4$}
\put(121,0){$5 \pp$}
\put(31,2){$\scriptstyle1+\eps$}
\put(47,2){$\scriptstyle1+2\eps$}
\put(55,2){$\scriptstyle1+3\eps$}
\put(80,2){$\scriptstyle2+\eps$}
\put(87,2){$\scriptstyle2+2\eps$}
\put(103,2){$\scriptstyle3+\eps$}
\put(7.5,7){$\scriptscriptstyle1+\eps^3$}
\put(15,4){$\scriptscriptstyle1+\eps^2$}
\put(23,4){$\scriptscriptstyle1+2\eps^2$}
\put(38,4){$\scriptscriptstyle1+\eps+\eps^2$}
\put(71,4){$\scriptscriptstyle2+\eps^2$}
\put(0,21){$(1)$}
\put(64,21){$(2)$}
\put(96,21){$(3)$}
\put(112,21){$(4)$}
\put(120,21){$(5) \pp$}
\put(31,19){$\scriptstyle(1,1)$}
\put(47,19){$\scriptstyle(1,2)$}
\put(55,19){$\scriptstyle(1,3)$}
\put(80,19){$\scriptstyle(2,1)$}
\put(87,19){$\scriptstyle(2,2)$}
\put(103,19){$\scriptstyle(3,1)$}
\put(5.5,14){$\scriptscriptstyle(1,0,0,1)$}
\put(14,17){$\scriptscriptstyle(1,0,1)$}
\put(22,17){$\scriptscriptstyle(1,0,2)$}
\put(38,17){$\scriptscriptstyle(1,1,1)$}
\put(70,17){$\scriptscriptstyle(2,0,1)$}
\put(-2,11){$0$}
\put(132,11){$1$}
\end{picture}
\caption{\smaller Realization of positions
by dyadic numbers in the unit interval~$[0, 1]$, and
the corresponding infinitesimal numbers as in
Figure~\ref{F:Positions}}
\label{F:DyadicPositions}
\end{figure}

Intuitively, $(i_1, \pp, i_p)$ corresponds to
what is denoted $i_1 + i_2\eps + \cdots + i_p\eps^{p-1}$
in Figure~\ref{F:Positions}. Under the dyadic realization,
we find $\ddd1 = 0$, $\ddd2 = \frac12$, $\ddd3 = \frac34$, \pp\
and $\ddd{1,2,1} = 0.01101 = \frac{13}{32}$
(Figure~\ref{F:DyadicPositions}). The
requirement that positions do not finish with~$0$ is needed
to guarantee that both the infinitesimal and the dyadic
realizations be injective on~$\NNNs$---alternatively, we
can allow final~$0$'s at the expense of identifying~$s$
and~$(s, 0)$.

\begin{figure} [htb]
\begin{picture}(100,17)(0, 0)
\put(0,4){\includegraphics{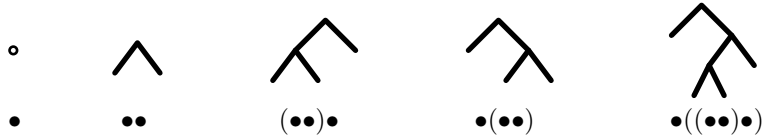}}
\put(0,0){$\et$}
\put(15,0){$\et \op \et$}
\put(36,0){$(\et \op \et) \op \et$}
\put(62,0){$\et \op (\et \op \et)$}
\put(88,0){$\et \op ((\et \op \et) \op \et)$}
\end{picture}
\caption{\smaller Typical trees and the
corresponding parenthesized expressions}
\label{F:ExampleTrees}
\end{figure}

The set of positions involved in an ordinary braid
diagram is an initial interval $\{1, 2, \pp, n\}$
of~$\NNN$. When we turn to~$\NNNs$, the role of such
an interval is played by a finite binary
tree---simply called a tree in the sequel. We
denote by~$\et$ the tree consisting of a single
vertex and by~$\tti \op
\ttii$ the tree with left subtree~$\tti$ and
right subtree~$\ttii$. Every tree has a unique
decomposition in terms of~$\et$, so
we can identify trees and parenthesized expressions
(Figure~\ref{F:ExampleTrees}). The {\it right height}
of a tree is defined to be the length of its
righmost branch.

Then we associate with every tree a finite set of
positions as follows:

\begin{defi} \label{I:Tree}
For~$\tt$ a tree, we define a finite set of dyadic
numbers~$\Dya(\tt)$ by the following rules: $\Dya(\et)$
is $\{0, 1\}$, and $\Dya(\tti\ttii)$ is the union
of~$\Dya(\tti)$ contracted from~$[0, 1]$ to~$[0,
\frac12]$ and of~$\Dya(\ttii)$ contracted
to~$[\frac12, 1]$. Then $\Pos\tt$ is defined to be
the set of all positions~$s$ such that $\dd s$
belongs to~$\Dya(\tt)$ with the largest two elements
removed.
\end{defi}

\begin{exam} \label{I:Vine}
Let $c_n$ denote the size~$n+1$ right vine
$\et(\pp (\et(\et\et))\pp )$, $n+1$~times~$\et$. Then
$\Dya(c_n)$ is $\{0,
\frac12, \frac34, \pp, 1-\frac1{2^n}, 1\}$,  \ie, 
$\{\ddd1, \ddd2, \pp, \ddd{n+1}, 1\}$, and  $\Pos{c_n}
$ is $\{(1), \pp, (n)\}$.  For $\tt =
\et((\et\et)\et)$ (the last tree in
Figure~\ref{F:ExampleTrees}), we find $\Dya(\tt)  =
\{0, \frac12, \frac58, \frac34, 1\}$, hence \break 
$\Dya(\tt) = \{\ddd1, \ddd2,
\ddd{2,1}, \ddd3, 1\}$, and $\Pos\tt = \{(1), (2),
(2,1)\}$. 
\end{exam}

\begin{lemm}
Every tree~$\tt$ is determined by the set of
positions~$\Pos\tt$.
\end{lemm}

\begin{proof}
An obvious induction shows that $\tt$ is
determined by~$\Dya(\tt)$. So the only problem is
that, in~$\Pos\tt$, the last two elements
of~$\Dya(\tt)$ are forgotten. Now the last element
of~$\Dya(\tt)$ is always~$1$, and an induction shows
that the forelast one is $\ddd{n+1}$, where $n$ is
maximal such that $(n)$ belongs to~$\Pos\tt$
(\eg, $\ddd3$, \ie, $\frac34$, in the
example above).
\end{proof}

\begin{rema} \label{R:Nota}
Instead of using~$\Dya(\tt)$, we can attribute to
each node in a binary tree an address that is a
sequence of positive integers as in
Figure~\ref{F:Addresses} below; then $\Pos\tt$
consists of the addresses of the leaves in~$\tt$,
up to removing the last address, diminishing by~$1$
all non-initial factors and removing the final~$0$'s
in each sequence. Our notational convention may seem
strange at first, because the initial and non-initial
entries in a position are not treated similarly in
the dyadic realization: the former is diminished by~$1$, the
latter are not. A more homogeneous definition would force
either to index  positions starting from~$0$---and therefore
numbering the braid generators~$\ss i$ from~$0$, which is
unusual---or to identify~$s$ with~$(s, 1)$ and not
with~$(s, 0)$---which is not intuitive.
\end{rema}

\subsection{Parenthesized braid diagrams} \label{I:Diagrams}

The diagrams we consider are constructed from
two series of elementary diagrams indexed by
letters~$\ssss i$ and~$\aaaa i$, and, therefore, a
diagram will be specified using a word on these
letters. In the sequel, such a word is called a
$\wsa$-word, or, simply, a {\it word}. A word
containing only letters~$\ssss i$ (\resp $\aaaa i$)
will be called a {\it $\ws$-word} (\resp an
{\it $\wa$-word}). Our aim is now to construct a
parenthesized diagram~$\DD_\tt(\ww)$ for~$\ww$ a
word and $\tt$ a large enough tree, exactly as the
ordinary diagram $\DD_n(\ww)$ is defined for $\ww$ a
word in the letters~$\ssss i$ and $n$ a large enough
integer. For~$\tt$ of size~$n+1$, hence defining
$n$~positions, $\DD_\tt(\ww)$ consists of
$n$~strands that connect the positions of~$\Pos\tt$,
considered as embedded in the unit interval, to
$n$~new positions. 

If $[x, y)$ and $[x', y')$ are
subintervals of~$[0,1)$, we say that we connect~$[x,
y)$ to~$[x', y')$ {\it homothetically} to mean that
each point~$(z,0)$ in~$[x, y) \times \{0\}$ is
connected to the point~$(z', 1)$ of~$[x', y') \times
\{1\}$ that satisfies $(z'-x')/(z' - y') = (z -
x)/(z-y)$.

\begin{figure} [htb]
\begin{picture}(85,28)(0, 0)
\put(0,3){\includegraphics{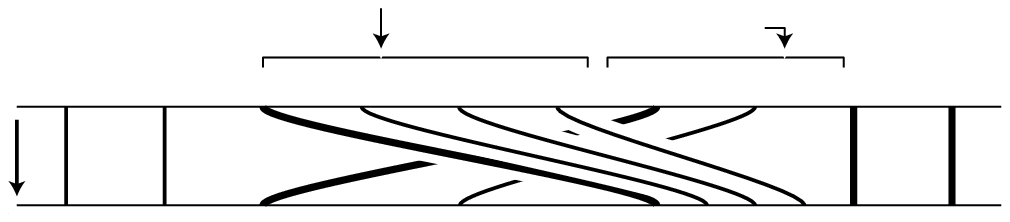}}
\put(24,0){$\scriptstyle \ddd i$}
\put(42,0){$\scriptstyle \ddd{i,s'}$}
\put(62,0){$\scriptstyle \ddd{i+1}$}
\put(71,0){$\scriptstyle \ddd{i\sp 1, s}$}
\put(84,0){$\scriptstyle \ddd{i+2}$}
\put(24,15.5){$\scriptstyle \ddd i$}
\put(42,15.5){$\scriptstyle \ddd{i,s}$}
\put(61,15.5){$\scriptstyle \ddd{i+1}$}
\put(71,15.5){$\scriptstyle \ddd{i\!+\!1 , s'}$}
\put(84,15.5){$\scriptstyle \ddd{i+2}$}
\put(15,25){each strand here crosses}
\put(41,22){over each strand there}
\put(-15,8){$\DD_\tt(\ss i):$}
\end{picture}
\begin{picture}(85,31)(0, 0)
\put(0,3){\includegraphics{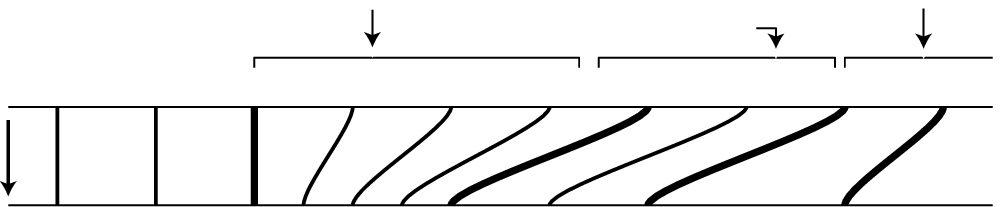}}
\put(24,0){$\scriptstyle \ddd i$}
\put(32,0){$\scriptstyle \ddd{i , 0 , s}$}
\put(42,0){$\scriptstyle \ddd{i, 1}$}
\put(63,0){$\scriptstyle \ddd{i+1}$}
\put(50,0){$\scriptstyle \ddd{i , j\!+\!1 , s'}$}
\put(84,0){$\scriptstyle \ddd{i+2}$}
\put(24,15){$\scriptstyle \ddd i$}
\put(43,15){$\scriptstyle \ddd{i,s}$}
\put(61,15){$\scriptstyle \ddd{i+1}$}
\put(70,15){$\scriptstyle \ddd{i\!+\!1 , j , s'}$}
\put(84,15){$\scriptstyle \ddd{i+2}$}
\put(94,15){$\scriptstyle \ddd{i+3}$}
\put(20,25){contraction by a factor $2$}
\put(67,25){dilatation by a factor $2$}
\put(58,21){translation}
\put(-15,8){$\DD_\tt(\aa i):$}
\end{picture}
\caption{\smaller The diagrams~$\DD_\tt(\ss i)$
and~$\DD_\tt(\aa i)$: in~$\DD_\tt(\ss i)$, the
positions coming from $[\ddd i, \ddd{i+1})$ and from
$[\ddd{i+1},
\ddd{i+2})$ are exchanged, with
a contraction/dilatation factor~$2$ due to the dyadic
realization; in~$\DD_\tt(\aa i)$, the positions
in~$[\ddd i, \ddd{i+1})$ are contracted by~$2$, those
in $[\ddd{i+1}, \ddd{i+2})$ are translated to the
left, and those in~$[\ddd k, \ddd{k+1})$ are
translated to the left and dilated by~$2$. In terms
of positions, $\DD_\tt(\ss i)$ exchanges $(i, s)$ and
$(i+1, s)$ for every~$s$, while
$\DD_\tt(\aa i)$ connects $(i, s)$ to~$(i, 0, s)$,
then $(i+1, j, s)$ to~$(i, j+1, s)$, and $(k, s)$
to~$(k-1, s)$ for $k \ge i+2$.}
\label{F:DefDiagrams}
\end{figure}

\begin{defi} (Figure~\ref{F:DefDiagrams}) \label{I:Diagram}
For $\tt$ a tree of right height at least~$i+1$, the
diagram~$\DD_\tt(\ss i)$ homothetically connects
$[\ddd i, \ddd{i+1})$ with $[\ddd{i+1}, \ddd{i+2})$,
then $[\ddd{i+1}, \ddd{i+2})$ with $[\ddd i, \ddd{i+1})$
with strands crossing under those of the previous
family, and, finally, \break $[\ddd k, \ddd{k+1})$ with
itself for $k \not= i, i+1$.

The diagram~$\DD_\tt(\aa i)$ homothetically connects
$[\ddd k, \ddd{k+1})$ with itself for
$k < i$, then $[\ddd i, \ddd{i+1})$ with $[\ddd{i},
\ddd{i, 1})$, next $[\ddd{i+1}, \ddd{i+2})$ with
$[\ddd{i, 1}, \ddd{i+1})$, and, finally, $[\ddd k,
\ddd{k+1})$ with $[\ddd{k-1}, \ddd k)$ for $k > i+1$.
\end{defi}

In contrast to the case of~$\Bi$, the
diagrams~$\DD_\tt(\ss i)$ or~$\DD_\tt(\aa i)$ so
defined cannot be carelessly stacked since the final
positions of the strands need not coincide with the
initial ones. Now, the changes correspond to an
easily described (partial) action on trees. 

\begin{figure} [htb]
\begin{picture}(130,28)(0, 0)
\put(0,3){\includegraphics{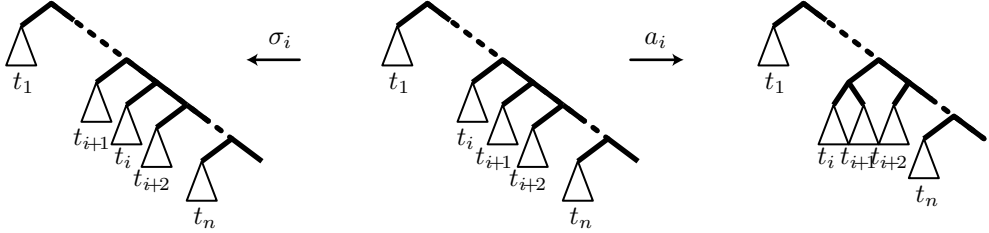}}
\put(1,18){$\tt_1$}
\put(9,10.5){$\tt_{i\!+\!1}$}
\put(14.5,8){$\tt_i$}
\put(17,5){$\tt_{i\!+\!2}$}
\put(25,0){$\tt_n$}
\put(35,24){$\ss i$}
\put(51,18){$\tt_1$}
\put(60,10.5){$\tt_i$}
\put(62.5,8){$\tt_{i\!+\!1}$}
\put(67,5){$\tt_{i\!+\!2}$}
\put(75,0){$\tt_n$}
\put(85,24){$\aa i$}
\put(101,18){$\tt_1$}
\put(108,8){$\tt_i$}
\put(111,8){$\tt_{i\!+\!1}$}
\put(115,8){$\tt_{i\!+\!2}$}
\put(121,3){$\tt_n$}
\end{picture}
\caption{\smaller Action of~$\ss i$ and~$\aa i$ on
a tree: $\ss i$ switches the $i$th and the
$(i+1)$st factors in the right decomposition, while
$\aa i$ glues them.}
\label{F:ActionOnTrees}
\end{figure}

\begin{defi} (Figure~\ref{F:ActionOnTrees})
\label{I:RD}
For $\tt$ a tree, the unique sequence of trees
$(\tt_1, \pp, \tt_n)$ such that $\tt$ factorizes as
$\tt_1 \op (\tt_2 \op ( \pp \op (\tt_n \op \et) \pp
))$ is called the {\it (right) decomposition}
of~$\tt$, and denoted by~$\Dec(\tt)$. For~$\tt$
a tree with $\Dec(\tt) = (\tt_1, \pp, \tt_n)$ with $n
> i$, we define the trees~$\tt \act \ss i$ and $\tt
\act \aa i$ by:
\begin{gather}
\label{E:ActionSOnTree}
\Dec(\tt \act \ss i) =
(\tt_1, \pp, \tt_{i-1}, \tt_{i+1}, \tt_i, \tt_{i+2},
\pp, \tt_n), \\
\label{E:ActionAOnTree}
\Dec(\tt \act \aa i) =
(\tt_1, \pp, \tt_{i-1}, \tt_i \op \tt_{i+1},
\tt_{i+2}, \pp, \tt_n).
\end{gather}
Then, one inductively defines $\tt \act \ww$ for
$\ww$ a word so that $\tt \act \ww\inv =
\ttt$ is equivalent to
$\ttt \act \ww = \tt$ and $\tt \act (\wwi
\wwii)$ is equal to $(\tt \act \wwi) \act \wwii$.
\end{defi}

The definition implies that the final
positions of the strands in~$\DD_\tt(\ss
i)$ and $\DD_\tt(\aa i)$ are $\Pos{\tt \act \ss i}$
and $\Pos{\tt \act \aa i}$, respectively.
Completing the construction of the
diagram~$\DD_\tt(\ww)$ is now obvious. 

\begin{defi}
The diagrams~$\DD_\tt(\sss i)$
and~$\DD_\tt(\aa i\inv)$ are defined to be the
mirror images of~$D_{\tt \act \ss i}(\ss i)$
and $D_{\tt \act \aa i}(\aa i)$, respectively.
Then, for~$\ww$ a word and~$\tt$ a binary
tree such that $\tt \act \ww$ is defined, the
{\it parenthesized braid diagram}~$\DD_\tt(\ww)$ is
inductively defined by the rule that, if~$\ww$
is~$\xx\ww'$ where $\xx$ is one of~$\ss i\pmo$, $\aa
i\pmo$, then
$\DD_\tt(\ww)$ is obtained by
stacking~$\DD_\tt(\xx)$ over $D_{\tt
\act
\xx}(\ww')$.
\end{defi}

\begin{figure} [htb]
\begin{picture}(138,53)(0, 0)
\put(5,3){\includegraphics{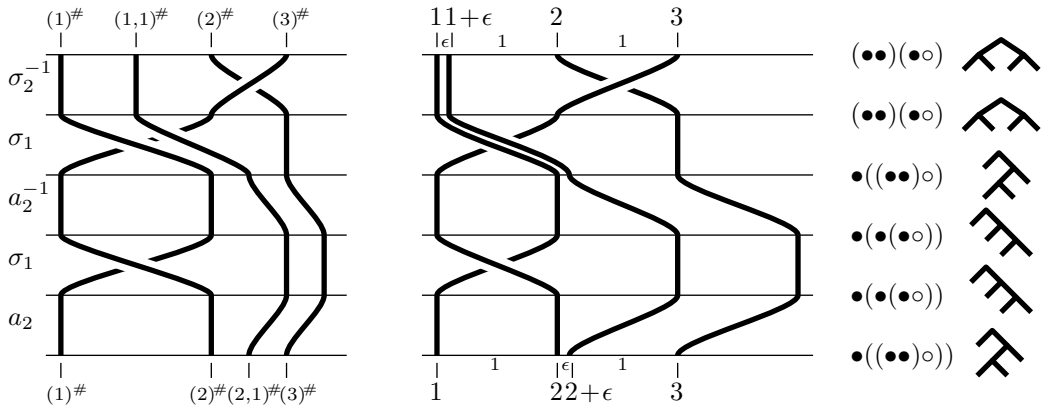}}
\put(0,42){$\sss2$}
\put(0,34){$\ss1$}
\put(0,26){$\aaa2$}
\put(0,18){$\ss1$}
\put(0,10){$\aa2$}
\put(5,0){$\scriptstyle \ddd1$}
\put(24,0){$\dddsmall 2$}
\put(29,0){$\dddsmall {2, 1}$}
\put(36,0){$\dddsmall 3$}
\put(5,50){$\dddsmall 1$}
\put(25,50){$\dddsmall 2$}
\put(14,50){$\dddsmall {1, 1}$}
\put(35,50){$\dddsmall 3$}
\put(56,0){$1$}
\put(72,0){$2$}
\put(74,0){$2\!+\!\eps$}
\put(88,0){$3$}
\put(56,50){$1$}
\put(72,50){$2$}
\put(58,50){$1\!+\!\eps$}
\put(88,50){$3$}
\put(73.5,4){$\scriptstyle\eps$}
\put(64,4){$\scriptstyle1$}
\put(81,4){$\scriptstyle1$}
\put(57.5,47){$\scriptstyle\eps$}
\put(65,47){$\scriptstyle1$}
\put(81,47){$\scriptstyle1$}
\put(112,45){$(\et\et)(\et\etwhite)$}
\put(112,37){$(\et\et)(\et\etwhite)$}
\put(112,29){$\et((\et\et)\etwhite)$}
\put(112,21){$\et(\et(\et\etwhite))$}
\put(112,13){$\et(\et(\et\etwhite))$}
\put(112,5){$\et((\et\et)\etwhite))$}
\end{picture}
\caption{\smaller The dyadic realization
of the diagram
$\DD_{(\et\et)(\et\et)}(\sss2 \ss1 \aaa2
\ss1 \aa2)$ and its infinitesimal version,
which (of course) is topologically equivalent; at
each step, the corresponding set of positions is
displayed, both as a parenthesized expression
(the last node is marked~$\circ$ because it
contributes no position) and as a binary tree.}
\label{F:Example}
\end{figure}

An example is displayed in
Figure~\ref{F:Example}. Ordinary braid diagrams are
special cases of parenthesized braid
diagrams: an $n$~strand braid diagram is a
diagram of the form~$\DD_\tt(\ww)$ where $\tt$
is the right vine of size~$n+1$ and $\ww$ is a
$\ws$-word.

An easy induction gives:

\begin{lemm}
For every tree~$\tt$ and every word~$\ww$, the
diagram~$\DD_\tt(\ww)$ is defined if and only if the
tree~$\tt \act \ww$ is, and, in this case, the final
positions in~$\DD_\tt(\ww)$ are~$\Pos{\tt \act
\ww}$.
\end{lemm}

\subsection{The group of parenthesized braids}

\begin{figure}[htb]
\begin{picture}(110,22)(0,2)
\put(0,0){\includegraphics[scale=1]{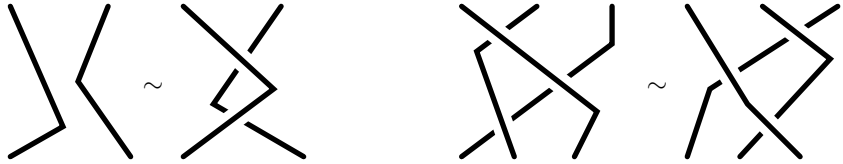}}
\put(18,10){$\sim$}
\put(85,10){$\sim$}
\end{picture}
\caption{\smaller Reidemeister moves of type II
(left) and III (right); the only requirement is that the
endpoints remain fixed}
\label{F:Reidemeister}
\end{figure}

According to Artin's original construction, braids
can be introduced as equivalence classes of braid
diagrams. Viewing a diagram as the projection of a
3D-figure, one considers the equivalence relation
corresponding to ambient isotopy of 3D-figures. As is
well-known, this amounts to declaring equivalent
those diagrams that can be connected by a finite
sequence of Reidemeister moves of types~II and~III
(Figure~\ref{F:Reidemeister}).

>From a topological point of view, parenthesized braid
diagrams are just ordinary diagrams, so they are
eligible for the same notion of equivalence: 

\begin{defi} \label{I:Equivalence}
Two parenthesized braid diagrams are declared {\it
equivalent} if and only if they can be transformed
one into the other by using Reidemeister moves of
types~II and~III (and keeping the endpoints
fixed). 
\end{defi} 

Our aim is to make a group out of parenthesized
braids---not only a groupoid, \ie, a category, as
in~\cite{Bar1, Bar2}. As mentioned above, the
problem is that we cannot compose arbitrary
diagrams. It can be solved easily by
introducing a completion procedure and defining the
group operation on the completed objects. In the
case of ordinary braids, the only parameter is the
number of strands, and, in order to compose two
diagrams~$\DD_{\nni}(\wwi)$, $\DD_{\nnii}(\wwii)$
with, say, $\nnii > \nni$, one first
completes~$\DD_{\nni}(\wwi)$ into the
$\nnii$-diagram~$\DD_{\nnii}(\wwi)$ obtained
from~$\DD_{\nni}(\wwi)$ by adding $\nnii - \nni$
unbraided strands on the right. 
The previous construction amounts to working with
infinite diagrams. For each braid word~$\ww$,
the diagrams~$\DD_\nn(\ww)$ make an inductive system
when $\nn$ varies, and, defining~$\DD_\infty(\ww)$
to be the limit of this system, we obtain a
well-defined product on infinite diagrams.
Moreover, as the completion preserves equivalence,
the product so defined induces a group structure,
namely that of~$\Bi$.

The procedure is similar for parenthesized braid
diagrams, the appropriate ordering being the
inclusion of trees viewed as sets of nodes. 

\begin{defi} \label{I:Completion}
For~$\tt, \ttt$ trees with $\tt \ince \ttt$, we
denote by~$\cc\tt\ttt$ the {\it completion} that
maps~$\DD_\tt(\ww)$ to~$\DD_{\ttt}(\ww)$ whenever
$\DD_\tt(\ww)$ exists.
\end{defi}

\begin{figure} [htb]
\begin{picture}(138,44)(0, -3)
\put(5,0){\includegraphics{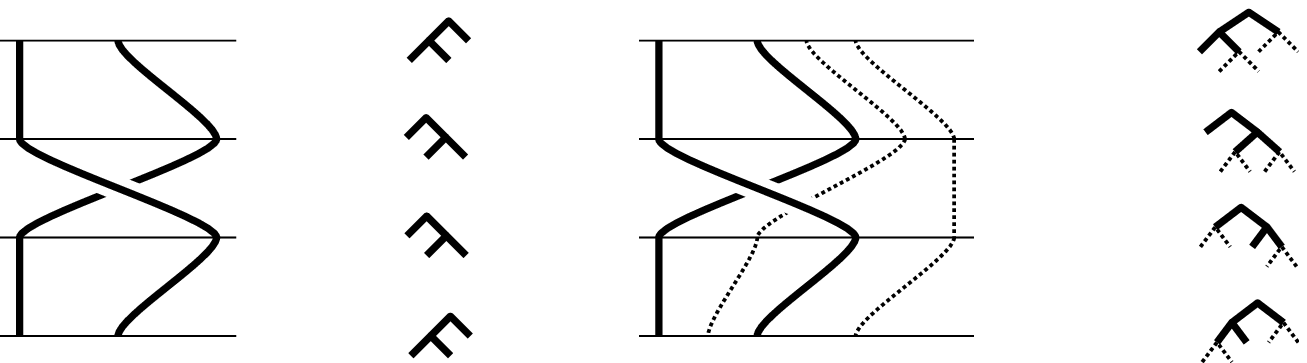}}
\put(0,28){$\aaa1$}
\put(0,18){$\ss1$}
\put(0,8){$\aa1$}
\put(65,28){$\aaa1$}
\put(65,18){$\ss1$}
\put(65,8){$\aa1$}
\put(5,36){$\ddsmall 1$}
\put(14,36){$\ddsmall {1, 1}$}
\put(5,0){$\ddsmall 1$}
\put(14,0){$\ddsmall {1, 1}$}
\put(35,33){$(\et\et)\etwhite$}
\put(35,23){$\et(\et\etwhite)$}
\put(35,13){$\et(\et\etwhite)$}
\put(35,3){$(\et\et)\etwhite$}
\put(70,35.5){$\ddsmall 1$}
\put(70,0.5){$\ddsmall 1$}
\put(78,35.5){$\ddsmall {1, 1}$}
\put(79,0.5){$\ddsmall {1, 1}$}
\put(83.5,39){$\ddsmall {1, 2}$}
\put(73,-3){$\ddsmall {1, 0, 1}$}
\put(90.5,39){$\ddsmall 2$}
\put(90,-3){$\ddsmall 2$}
\put(107,33){$(\et(\et\et))(\et\etwhite)$}
\put(107,23){$\et((\et\et)(\et\etwhite))$}
\put(107,13){$(\et\et)(\et(\et\etwhite))$}
\put(107,3){$((\et\et)\et)(\et\etwhite)$}
\end{picture}
\caption{\smaller Completion
of $\DD_{(\et\et)\et}(\aaa1 \ss1 \aa1)$ into
$\DD_{(\et(\et\et))(\et\et)}(\aaa1 \ss1
\aa1)$: two more leaves in the tree, two more
strands in the braid}
\label{F:Completion}
\end{figure}

The explicit construction of parenthesized braid diagrams
makes the completion procedure easy: as shown
on Figure~\ref{F:Completion}, the
diagram~$\DD_{\ttt}(\ww)$ for $\ttt \supseteq \tt$
is obtained by keeping the existing strands, and
adding new strands in~$\DD_{\tt}(\ww)$ that always
lie half-way between their left and right neighbours---or
$1$ if there is no right neighbour. The only difference with
ordinary diagrams is that there is in general more than one
basic extension: the only way to extend the interval~$\{1, 2,
\pp, n\}$ into a bigger interval is to add~$n+1$
while, in a tree~$\tt$, each leaf can be split into a
caret with two leaves, so there are $n+1$~basic
extensions when $\tt$ specifies $n$~positions. As
an induction shows, splitting the $k$th leaf amounts
to doubling the $k$th strand. 

The following observations gather what is needed
for mimicking the construction of~$\Bi$:

\begin{lemm}
(i) For each word~$\ww$, the
system~$(\DD_\tt(\ww), \cc\tt\ttt)$ is directed; 

(ii) Diagram concatenation induces a well-defined
product on direct limits;

(iii) The completion maps are compatible with
diagram equivalence.
\end{lemm}

\begin{proof}
For~$(i)$, for any two trees~$\tti, \ttii$, there
exists a tree~$\tt$ that includes both~$\tti$
and~$\ttii$, for instance the tree whose nodes are
the union of the nodes in~$\tti$ and~$\ttii$.
For~$(ii)$, the completion~$\cc\tt\ttt$ is
compatible with the product in that, if
$\DD_{\tt}(\wwi)$ and $\DD_{\tt
\act \wwi}(\wwii)$ exist so that
$\DD_\tt(\wwi \wwii)$ is defined, then, for each
tree~$\ttt$ including~$\tt$, the diagram
$\DD_{\ttt}(\wwi \wwii)$ exists and we have
$$\DD_{\ttt}(\wwi \wwii) = \DD_{\ttt}(\wwi) \opp
\DD_{\ttt \act \wwi}(\wwii).$$
Finally, $(iii)$ follows from the description
of completion in terms of strand addition.
\end{proof}

For each word~$\ww$, let us define $\DDs(\ww)$ to
be the direct limit---actually, by construction, just
the union---of the inductive system of
all~$\DD_\tt(\ww)$'s. We call it an {\it infinite
parenthesized braid diagram}.  Then concatenation
induces an everywhere defined product on infinite
parenthesized braid diagrams, and isotopy induces a well-defined
equivalence relation that is compatible with the
previous product. Then the same argument as for ordinary
braid diagrams gives:

\begin{prop}
Isotopy classes of infinite parenthesized braid diagrams make a
group.
\end{prop}

\begin{defi} \label{I:Braids}
The group of isotopy classes of infinite parenthesized braid
diagrams is called the {\it group of parenthesized
braids}, and denoted~$\Bs$; its elements
are called {\it parenthesized braids}.
\end{defi}

\subsection{Relations in~$\Bs$}

By construction, the group~$\Bs$ is generated by the
elements~$\ss i$ and~$\aa i$. An obvious task is
to look for a presentation in terms of these
elements. For the moment, we just observe that
certain relations are satisfied in~$\Bs$. That these
relations make a presentation of~$\Bs$ will be
established in Section~\ref{S:LD} below. 

\begin{lemm}\label{L:Relations}
For $i \ge 1$ and $j \ge i+2$, the following relations 
induce diagram isotopies, hence equalities in~$\Bs$:
\begin{equation} \label{E:Relations}
\begin{cases}\quad \ss i \ss j  = \ss j \ss i,
\qquad \ss i \aa j  = \aa j \ss i,
\qquad \aa i \aa{j-1} = \aa j \aa i,
\qquad \aa i \ss{j-1} = \ss j \aa i,\\
\quad \ss i \ss{i+1} \ss i 
= \ss{i+1} \ss i \ss{i+1},
\text{\qquad}
\ss{i+1} \ss i \aa{i+1}
= \aa i \ss i,
\text{\qquad}
\ss i \ss{i+1} \aa i 
= \aa{i+1} \ss i.
\end{cases}
\end{equation}
\end{lemm}

\begin{figure}[htb]
\begin{picture}(130,60)(0,0)
\put(0,4){\includegraphics{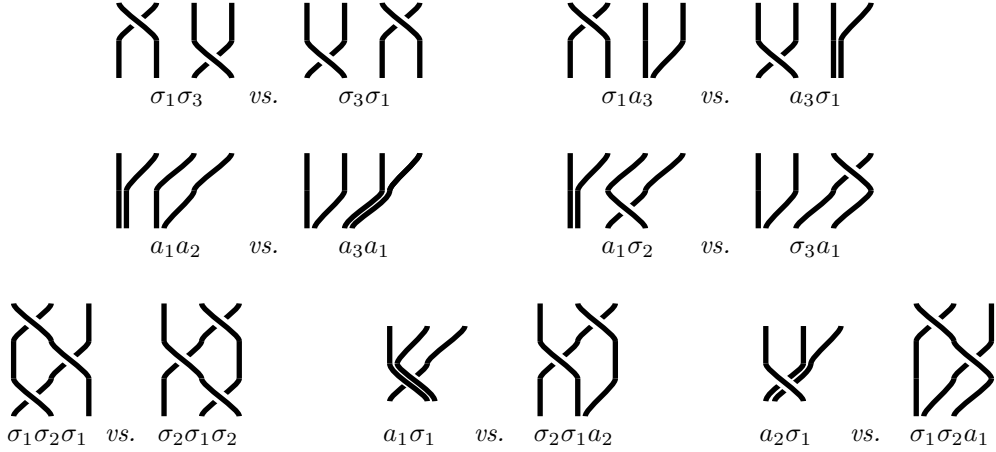}}
\put(19,46){$\ss1 \ss3$}
\put(32,46){\it vs.}
\put(44,46){$\ss3 \ss1$}
\put(79,46){$\ss1 \aa3$}
\put(92,46){\it vs.}
\put(104,46){$\aa3 \ss1$}
\put(19,26){$\aa1 \aa2$}
\put(32,26){\it vs.}
\put(44,26){$\aa3 \aa1$}
\put(79,26){$\aa1 \ss2$}
\put(92,26){\it vs.}
\put(104,26){$\ss3 \aa1$}
\put(0,1){$\ss1 \ss2 \ss1$}
\put(13,1){\it vs.}
\put(20,1){$\ss2 \ss1 \ss2$}
\put(50,1){$\aa1 \ss1$}
\put(62,1){\it vs.}
\put(70,1){$\ss2 \ss1 \aa2$}
\put(100,1){$\aa2 \ss1$}
\put(112,1){\it vs.}
\put(120,1){$\ss1 \ss2 \aa1$}
\end{picture}
\caption{\smaller The relations of~$\RRs$ and
the corresponding diagrams isotopies (here in
infinitesimal realization)}
\label{F:Relations}
\end{figure}

\begin{proof}
The graphical verification is given in
Figure~\ref{F:Relations}.
\end{proof}

Relations~\eqref{E:Relations} include the standard
braid relations, as well as the relations
$\aa i \aa j = \aa{j-1} \aa i$ for $j \ge i+2$,
which correspond to the standard presentation of
Thompson's group~$F$ up to the change of name~$\aa
i = x_{i-1}\inv$.
In order to subsequently prove that \eqref{E:Relations}
gives a presentation of~$\Bs$, it is
convenient to introduce the abstract group
presented by these relations. 

\begin{defi} \label{I:BBraids}
We denote by~$\ss*$ and~$\aa*$ the families of all~$\ss i$'s
and of all~$\aa i$'s, and by~$\RRs$ the
relations~\eqref{E:Relations}. We define~$\BBs$ to be the
group~$\Gr(\XXsa; \RRs)$.
\end{defi}

Lemma~\ref{L:Relations} states that the identity
mapping on~$\XXs$ and~$\XXa$ induces a surjective
morphism of~$\BBs$ onto~$\Bs$. One of our aims will
be to prove that this morphism is an isomorphism.

\section{Algebraic properties of the group~$\BBs$}
\label{S:AlgebraicProp}

A number of algebraic properties of the group~$\BBs$
can be deduced from its explicit presentation, as we
shall easily see using a specific combinatorial
method called word reversing. The main results we
prove are that $\BBs$ is a group of left fractions,
that it is torsion-free, and that it contains copies
of the braid group~$\Bi$ as well as of Thompson's
group~$F$.

\subsection{The word reversing technique}

In order to study the group~$\BBs$, we resort to
general algebraic tools developed in~\cite{Dff, Dgp}
and connected with Garside's seminal work~\cite{Gar}.
This combinatorial method applies to monoid
presentations and it is relevant for establishing
properties like cancellativity or embeddability in a
group of fractions.

For~$\XX$ a nonempty set (of letters), we call 
{\it $\XX$-word} a word made of letters
from~$\XX$, and {\it $\XXpm$-word} a word
made of letters from~$\XX\cup \XX\inv$,
where $\XX\inv$ is a disjoint copy of~$\XX$
containing one letter~$\gx\inv$ for each~$\gx$
in~$\XX$. Then $\XX$-words are called positive,
and we say that a group presentation
$(\XX, \RR)$ is {\it positive} if
$\RR$ exclusively consists of relations $\uu = \vv$
with $\uu, \vv$ nonempty positive words. We
use~$\Gr(\XX; \RR)$ for the group
and~$\Mon(\XX; \RR)$ for the monoid defined
by~$(\XX, \RR)$. Note that the
presentation~$(\XXsa, \RRs)$ is positive.

\begin{defi} \cite{Dff, Dgp} \label{I:Reversing}
Let $(\XX, \RR)$ be a positive group presentation,
and $\ww, \ww'$ be $\XXpm$-words.
We say that $\ww$ is {\it right
$\RR$-reversible} to~$\ww'$, denoted $\ww
\revr\RR \ww'$, if
$\ww'$ can be obtained from~$\ww$ using finitely
many steps consisting either in deleting some
length~$2$ subword~$\gx\inv \gx$, or in
replacing a length~$2$ subword~$\gx\inv \gy$ by a
word~$\vv \uu\inv$ such that $\gx \vv = \gy \uu$
is a relation of~$\RR$.
\end{defi}

Right $\RR$-reversing uses the relations
of~$\RR$ to push the negative letters (those
in~$\XX\inv$) to the right and the positive letters
(those in~$\XX$) to the left by iteratively
reversing $-+$ patterns into $+-$ patterns. Note
that deleting $\gx\inv\gx$ enters the
general scheme if we assume that, for every
letter~$\gx$ in~$\XX$, the trivial relation
$\gx = \gx$ belongs to~$\RR$.

Left $\RR$-reversing is defined symmetrically: the
basic step consists in deleting a subword~$\gx
\gx\inv$, or replacing a subword~$\gx \gy\inv$
with~$\vv\inv \uu$ such that $\vv \gx = \uu \gy$ is a
relation of~$\RR$.

\begin{figure}[htb]
\begin{picture}(48,25)(0, 0)
\put(3,2){\includegraphics{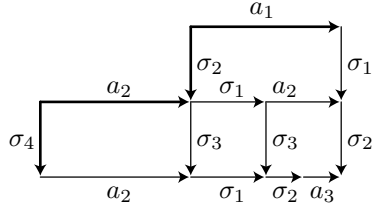}}
\put(13,0){$\aa2$}
\put(28,0){$\ss1$}
\put(35,0){$\ss2$}
\put(40,0){$\aa3$}
\put(0,7){$\ss4$}
\put(25,7){$\ss3$}
\put(35,7){$\ss3$}
\put(45,7){$\ss2$}
\put(13,14){$\aa2$}
\put(28,14){$\ss1$}
\put(35,14){$\aa2$}
\put(25,17){$\ss2$}
\put(45,17){$\ss1$}
\put(32,24.5){$\aa1$}
\end{picture}
\caption{\smaller Right reversing diagram for $\sss4
\aa2 \sss2 \aa1$: one starts with a staircase
labelled $\sss4 \aa2 \sss2 \aa1$ by drawing a
vertical $\gx$-labelled arrow for each
letter~$\gx\inv$, and an horizontal $\gy$-labelled
arrow for each positive letter~$\gy$. Then, when
$\gx\inv \gy$ is replaced with $\vv \uu\inv$, we
complete the open pattern corresponding to  $\gx\inv
\gy$ into a square by adding horizontal
$\vv$-labelled arrows and vertical $\uu$-labelled
arrows.}
\label{F:ExReversing}
\end{figure}

\begin{exam}
Let us consider the presentation~$(\XXsa, \RRs)$, and
let~$\ww$ be the word~$\sss4 \aa2 \sss2 \aa1$.
Then $\ww$ contains two $-+$-subwords, namely
$\sss4 \aa2$ and~$\sss2 \aa1$. So there are 
two ways of starting a right reversing from~$\ww$:
replacing $\sss4 \aa2$ with $\aa2 \sss3$, which is
legal as $\ss4
\aa2 = \aa2 \ss3$ is a relation of~$\RRs$,
or replacing $\sss2 \aa1$ with 
$\ss1 \aa2 \sss1$, owing to the relation $\ss2 (\ss1
\aa2) = \aa1 \ss1$. The reader can check that, in
any case, iterating the process leads in
four steps to $\aa2 \ss1 \ss2 \aa3 \sss2
\sss1$. The latter word is terminal since it contains
no $-+$ subword. It is helpful to visualize the
process using a planar diagram similar to a Van Kampen
diagram as shown in Figure~\ref{F:ExReversing}.
\end{exam}

If $\gx \uu = \gy \vv$ is a relation of~$\RR$,
then $\gx\inv \gy$ and~$\vv \uu\inv$ are
$\RR$-equivalent, hence $\ww \revr\RR \ww'$
implies that $\ww$ and~$\ww'$ represent the same
element of~$\Gr(\XX; \RR)$. A slightly more careful
argument shows that, if $\uu, \vv, \uu', \vv'$ are
positive words, then $\uu\inv \vv \revr\RR \vv'
{\uu'}\inv$ implies that $\uu \vv'$ and~$\vv \uu'$
represent the same element of~$\Mon(\XX; \RR)$.
So, in particular, if $\uu, \vv$ are positive words,
$\uu\inv \vv \revr\RR \e$ (the empty word)
implies that $\uu$ and~$\vv$ represent the same
element of~$\Mon(\XX; \RR)$. The converse need
not be true in general, but the interesting case is
when this happens:

\begin{defi} \cite{Dgp} \label{I:Complete}
A positive presentation~$(\XX, \RR)$ is said to be
{\it complete for right reversing} if right
reversing always detects positive equivalence, in the sense
that, for all $\XX$-words~$\uu, \vv$, one has
$\uu\inv \vv \revr\RR \e$ whenever $\uu$ and~$\vv$
represent the same element of~$\Mon(\XX; \RR)$. 
\end{defi}

Symmetrically, we say that
$(\XX; \RR)$ is complete for left
reversing if $\uu \vv\inv$ is left $\RR$-reversible
to~$\e$ whenever $\uu$ and~$\vv$ represent the same
element of~$\Mon(\XX; \RR)$. The point is that
there exists a tractable criterion for recognizing
whether a given presentation is complete for 
reversing---or for adding new relations  if it is
not.

\begin{defi} \label{I:Homogeneous}
A positive presentation~$(\XX, \RR)$ 
is said to be {\it homogeneous} if there exists a
$\RR$-invariant mapping~$\lambda$
from~$\XX$-words to nonnegative integers such
that $\lambda(\gx) \ge 1$ holds for every~$\gx$
in~$\XX$, and $\lambda(\uu\vv) \ge \lambda(\uu)
+ \lambda(\vv)$ holds for all $\XX$-words~$\uu,
\vv$.
\end{defi}

If all relations in~$\RR$ preserve the length of the
words, then the length satisfies the requirements
for the function~$\lambda$ and the presentation is
homogeneous. 

\begin{prop}\cite{Dgp}\label{P:Cube}
A homogeneous positive presentation $(\XX, \RR)$ is
complete for right reversing if and only if
the following condition holds for each triple~$(\gx,
\gy, \gz)$ of letters:
\begin{equation} \label{E:Cube}
\text{$\gx\inv \gy \gy\inv \gz \revr\RR \vv
\uu\inv$ \quad with~$\uu, \vv$ positive
implies \quad $\vv\inv \gx\inv \gz \uu
\revr\RR \e$.}
\end{equation}
\end{prop}

Condition~\eqref{E:Cube} is called the {\it right
cube condition for~$(\gx, \gy, \gz)$}. Of course, a
symmetric left cube condition guarantees
completeness for left reversing. We shall see now
that the presentation $(\XXsa; \RRs)$ is eligible
for the previous criterion.

\begin{lemm}
The presentation $(\XXsa; \RRs)$ is homogeneous.
\end{lemm}

\begin{proof}
The relations $\ss i \ss{i+1} \aa i = \aa{i+1} 
\ss i$ and $\ss{i+1} \ss i \aa{i+1} = \aa i \ss i$
do not preserve the length, so the latter
cannot be used directly. Instead we construct a
twisted length function~$\lambda$ so that,
in~$\lambda(\ww)$, each letter~$\aa i$
contributes~$1$, but $\ss i$ contributes~$nn'$,
where $n$ and~$n'$ are the numbers of strands
involved in the diagram~$\DD_c(\ww)$ for $c$ a
sufficiently large right vine. Formally, we first
define an action of positive
words on sequences of integers by:
\begin{align*}
(\pp, n_{i-1}, n_i, n_{i+1}, n_{i+2}, \pp) \act
\aa i  
&= (\pp, n_{i-1}, n_i + n_{i+1}, n_{i+2}, \pp),\\
(\pp, n_{i-1}, n_i, n_{i+1}, n_{i+2}, \pp) \act \ss i 
&= (\pp, n_{i-1}, n_{i+1}, n_i, n_{i+2}, \pp).
\end{align*}
Then $n_i$ is the number of strands near
position~$i$, \ie, corresponding to positions~$(i,
s)$, in~$\DD_c(\ww)$, and the action
is compatible with the relations of~$\RRs$. Then,
for~$\ww$ a positive word, we
put $\lambda_\bullet(\aa i, \ww) = 1$ and
$\lambda_\bullet(\ss i, \ww) = n_i n_{i+1}$ for
$(1, 1, \pp) \act \ww = (n_1, \pp, n_p)$. Finally, 
we define $\lambda(\ww) = \sum_k
\lambda_\bullet(\ww(k), \ww_k)$, where $\ww(k)$
denotes the $k$th letter in~$\ww$ and $\ww_k$
denotes the length~$k-1$ prefix of~$\ww$. Then
$\lambda$ witnesses that $(\XXsa, \RRs)$ is
homogeneous.
\end{proof}

\begin{lemm}
The presentation $(\XXsa; \RRs)$ satisfies the right
and the left cube conditions for each triple of
letters.
\end{lemm}

\begin{proof}
As there are infinitely many letters, infinitely
many cases are to be considered. However, it is clear
that only the mutual
distance of the indices matter, and, therefore, only
finitely many types occur. The verification is easy,
and we postpone it to an appendix.
\end{proof}

Applying the criterion of Proposition~\ref{P:Cube},
we deduce:

\begin{prop} \label{P:Completeness}
The presentation $(\XXsa; \RRs)$ is
complete for both right and left
reversing.
\end{prop}

\subsection{The monoid~$\BBsp$}

Once the presentation~$(\XXsa, \RRs)$ is known to be
complete for reversing, a number of
results can be established easily. We begin with
results involving the monoid presented by the
relations~$\RRs$.

\begin{defi} \label{I:Monoid}
We denote by~$\BBsp$ the monoid~$\Mon(\XXsa; \RRs)$.
\end{defi}

The elements of the monoids~$\BBsp$ are represented
by positive words, and, by definition of
completeness, two such words~$\uu, \vv$
represent the  same element in~$\BBsp$ if and only
$\uu\inv \vv$ is right $\RRs$-reversible to the
empty word, if and only if $\uu \vv\inv$ is left
$\RRs$-reversible to the empty word.
Let us begin with cancellativity. The following
criterion tells us that, whenever the presentation
is complete, the monoid is cancellative provided
there is no obvious obstruction.  

\begin{lemm} \cite{Dgp}\label{P:Canc}
Assume that $(\XX, \RR)$ is a positive
presentation that is complete for right
reversing. Then $\Mon(\XX; \RR)$ is left
cancellative whenever $\RR$ contains no relation of
the form $\gx \uu = \gx \vv$.
\end{lemm}

There is no relation of the form $\aa i \uu = \aa i
\vv$, $\ss i \uu = \ss i \vv$, $\uu \aa i = \vv \aa
i$, $\uu \ss i = \vv \ss i$ in~$\RRs$, so, using the
previous criterion and its symmetric counterpart, we
deduce:

\begin{prop}
The monoid~$\BBsp$ admits left and right cancellation.
\end{prop}

Let us now consider common multiples. Say that
$\ez$ is a least common right multiple, or right
lcm, of two elements~$\ex, \ey$ in a monoid~$M$ if
$\ez$ is a right multiple of~$\ex$ and~$\ey$, \ie,
$\ez = \ex \ex' = \ey \ey'$ holds for some~$\ex',
\ey'$, and every common right multiple of~$\ex$
and~$\ey$ is a right multiple of~$\ez$.

\begin{lemm}\cite{Dgp}\label{L:CritLcm}
Assume that $(\XX, \RR)$ is a positive
presentation that is complete for right reversing.
Then a sufficient condition for any two elements
admitting a common right multiple to admit a right
lcm is that, for all~$\gx,
\gy$ in~$\XX$, there is at most one relation of the
form $\gx \uu = \gy \vv$ in~$\RR$. In
that case, for all $\XX$-words~$\uu, \vv$, the
word~$\uu\inv \vv$ is right reversible to a word
of the form~$\vv' {\uu'}\inv$ with $\uu', \vv'$
positive if and only if the elements
represented by~$\uu$ and~$\vv$ in~$\Mon(\XX;
\RR)$ admit a common right multiple, and then
$\uu \vv'$ represents the right lcm of these
elements.
\end{lemm}

The presentation~$(\XXsa, \RRs)$ is
eligible for the previous criterion, and we deduce:

\begin{prop}
Any two elements of the monoid~$\BBsp$ that admit a
common right (\resp left) multiple admit a right
(\resp left) lcm.
\end{prop}

Standard arguments imply:

\begin{coro}
Any two elements of the monoid~$\BBsp$ admit a left
and a right gcd.
\end{coro}

It remains to study whether common multiples do
exist in~$\BBsp$. For right multiples, the answer
is negative: Lemma~\ref{L:CritLcm} tells us that the
elements~$\aa1$ and~$\aa2$ admit a common right
multiple in~$\BBsp$ if and only if the right
reversing of the word~$\aaa1 \aa2$ leads in a finite
number of steps to some positive--negative word. As
there is no relation of the form~$\aa1 \uu = \aa2
\vv$ in~$\RRs$, this cannot happen, and, therefore,
$\aa1$ and~$\aa2$ have no common right multiple
in~$\BBsp$. The situation is different for left
multiples. In order to describe it, we
need some notation.

\begin{defi} \label{I:Double}
For~$\ww$ a $\ws$-word and $k$ a positive integer, we
denote by~$\ww[k]$ the initial position of the strand
that finishes at position~$k$ in the braid
diagram~$\DD(\ww)$, and by~$\dbl_k(\ww)$ the braid
word that encodes the diagram obtained
from~$\DD(\ww)$ by doubling the strand starting at
position~$k$. Similar notations are used for braids,
which is legal as the needed compatibilities are
satisfied.
\end{defi}

Thus we have $\e[k] = k$ and $\dbl_k(\e) = \e$ for
every~$k$, and
\begin{gather}
\ss i[k] = 
\begin{cases}
k & \text{ for $k \not= i, i+1$,}\\
i+1 & \text{ for $k = i$,}\\
i & \text{ for $k = i+1$,}
\end{cases} \qquad
\dbl_k(\ss i) = 
\begin{cases}
\ss{i+1} & \text{ for $k < i$,}\\
\ss{i+1}\ss i & \text{ for $k = i$,}\\
\ss i \ss{i+1} & \text{ for $k = i+1$,}\\
\ss i & \text{ for $k > i+1$,}
\end{cases}\\
 \label{E:Double}
\ww[k] = \wwi[\wwii[k]], \qquad
\dbl_k(\ww) = \dbl_k(\wwi) \opp
\dbl_{\wwi\inv[k]}(\wwii)
\qquad\text{for $\ww = \wwi\wwii$}.
\end{gather}

\begin{lemm} \label{L:Convergence}
Left $\RRs$-reversing always terminates in finitely
many steps.
\end{lemm}

\begin{proof}
The result is not {\it a priori} obvious as the
length of the words appearing during the reversing
may increase. By Garside's theory, any two
elements in the braid monoid~$\Bip$ admit a common
left multiple, and, therefore, the left reversing of
any word~$\uu \vv\inv$ with $\uu, \vv$ positive
$\ws$-words terminates in finitely many steps.
The same is true for $\wa$-words, since, in this
case, the length cannot increase. The only
remaining case is that of mixed words involving both
types of letters. Now, in this case, we can
describe the result of reversing explicitly.
Indeed, we claim that, for every positive
$\ws$-word~$\ww$ and every positive integer~$k$,
\begin{equation} \label{E:LeftRevSA}
\text{$\ww \opp \aa k\inv$ is left
$\RRs$-reversible to $\aa{\ww[k]}\inv \opp
\dbl_{\ww[k]}(\ww)$.}
\end{equation}
We use induction on~$\ww$. For $\ww = \ss
i$, one easily checks \eqref{E:LeftRevSA} in the
various cases. For instance, $\ss1 \aa1\inv$ is left
reversible to $\aa2\inv \ss2 \ss1$, and we have
$\ss1[1] = 2$ and  $\dbl_2(\ss1) = \ss1 \ss2$. Then,
for $\ww = \wwi \wwii$, using the definition of left
reversing and the hypothesis that
\eqref{E:LeftRevSA} holds for~$\wwi$ and~$\wwii$, we
obtain that  
$\wwi \wwii \aaa k$ is left reversible to
$\wwi \aaa{\wwii[k]} 
\dbl_{\wwii[k]}(\wwii)$, and then to
$\aaa{\wwi[\wwii[k]]} 
\dbl_{\wwi[\wwii[k]]}(\wwi) 
\dbl_{\wwii[k]}(\wwii)$, which, by~\eqref{E:Double},
is $\aaa{\ww[k]} \dbl_{\ww[k]}(\ww)$.
\end{proof}

Applying Lemma~\ref{L:CritLcm}, we deduce:

\begin{prop} \label{P:LcmB}
Any two elements in the monoid~$\BBsp$ admit a left
lcm.
\end{prop}

Another merit of word reversing is to make it easy
to recognize what we can call parabolic submonoids
(and, similarly, subgroups).

\begin{lemm} \label{L:Parabolic}
Assume that $(\XX, \RR)$ is a positive presentation
that is complete for left reversing, and $\XXo$
is a subset of~$\XX$. Let~$\RRo$ be the set of
all relations $\vv \gx = \uu \gy$ in~$\RR$ with
$\gx, \gy \in \XXo$. If all words occurring
in~$\RRo$ are $\XXo$-words,
the submonoid of~$\Mon(\XX; \RR)$ generated
by~$\XXo$ admits the presentation~$\Mon(\XXo;
\RRo)$.
\end{lemm}

\begin{proof}
The point is to prove that, if $\uu, \vv$ are
$\RR$-equivalent $\XXo$-words, then
$\uu$ and~$\vv$ also are $\RRo$-equivalent, \ie, no
relation in $\RR \setminus \RRo$ is neeeded to prove
their equivalence. Now, by completeness, $\uu$ and
$\vv$ being $\RR$-equivalent implies that $\vv
\uu\inv$ is left $\RR$-reversible to~$\e$. The
hypothesis on~$\RRo$ implies that
only letters from~$\XXo$ appear during the reversing
process. Therefore, the latter is an
$\RRo$-reversing, and $\uu$ and~$\vv$ are
$\RRo$-equivalent.
\end{proof}

We denote by~$\Fp$ the monoid with presentation
$\Mon(\XXa; \aa i \aa j = \aa{j-1} \aa i
\text{~for~}j \ge i+2)$, and call it {\it
Thompson's monoid}. 

\begin{prop} \label{P:Zappa}
The submonoid of~$\BBsp$ generated by~$\XXs$
is (isomorphic to) the braid monoid~$\Bip$, while the
submonoid generated by~$\XXa$ is
(isomorphic to) Thompson's monoid~$\Fp$. Each
element of~$\BBsp$ admits a unique decomposition
in~$\Bip \times \Fp$. The monoid~$\BBsp$ is the
Zappa-Sz\'ep product of~$\Bip$ and~$\Fp$ associated
with the crossed product defined for~$\bx \in \Bip$
and $k \ge 1$ by
\begin{equation} \label{E:Zappa}
\aa k \opp \bx = \dbl_k(\bx) \opp
\aa{\bx\inv[k]}.
\end{equation}
\end{prop}

\begin{proof}
An inspection of the relations in~$\RRs$ shows that
the families~$\XXs$ and~$\XXa$ are eligible for the
criterion of Lemma~\ref{L:Parabolic}, and the first
part of the proposition follows. We henceforth
identify~$\Bip$ and~$\Fp$ with the subgroups
of~$\BBsp$ generated by~$\XXs$ and~$\XXa$,
respectively.

Formula~\eqref{E:Zappa} is a direct consequence
of~\eqref{E:Double}, and, by a straightforward
induction, it implies $\BBsp = \Bip \opp \Fp$. So
the only point to prove is the uniqueness of the
decomposition in~$\Bip \times \Fp$.
Assume that $\uu \vv$ and $\uuu \vvv$ are
$\RRs$-equivalent, where $\uu, \uuu$ are $\ws$-words
and $\vv, \vvv$ are $\wa$-words. By
completeness, this means that $\uu \vv \uuu\Inv
\vvv\Inv$ is left reversible to the empty word.
Let $\uu_1, \vv_1, \uuu_1, \vvv_1$ be the
intermediate words appearing in the reversing,
as shown in
\vrule width0pt depth9mm
\begin{picture}(31,10)(0,8)
\put(2,2){\includegraphics{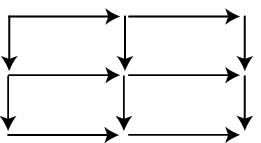}}
\put(7,0){$\uu$}
\put(7,10){$\uui$}
\put(19,10){$\vvi$}
\put(19,0){$\vv$}
\put(7,16){$\e$}
\put(19,16){$\e$}
\put(0,6){$\e$}
\put(0,12){$\e$}
\put(28,6){$\vv'$}
\put(28,12){$\uu'$}
\put(15.5,5.5){$\vv'_1$}
\put(15.5,11.5){$\uu'_1$}
\end{picture}
. As $\vv$ and $\vvv$
are positive $\wa$-words, so are $\vv_1$
and~$\vvv_1$. By~\eqref{E:Double}, the letters~$\aaa
k$ never vanish when they cross $\ss i$'s in a
left reversing. Hence the only possibility for $\uu
\vvv_1\Inv$ to reverse to a positive word~$\uu_1$ is
that $\vvv_1$ is empty. Similarly, $\vv_1$ must be
empty. As $\vv_1$ and~$\vvv_1$ are empty, $\vv$
and~$\vvv$ are $\RRs$-equivalent. On the other
hand, $\vv_1$ and~$\vvv_1$ being empty implies
$\uu_1 = \uu$ and $\uuu_1 = \uuu$, so the
hypothesis that $\uu_1 \uuu_1\Inv$ reverses
to~$\e$ implies that $\uu$ and~$\uuu$ are
$\RRs$-equivalent.
(For general Zappa-Sz\'ep products, see~\cite{Bri3}---or
\cite{Pic} where the name ```crossed product'' is used.)
\end{proof}

\subsection{The group~$\BBs$}

It is now easy to deduce results about the
group~$\BBs$. 

\begin{prop} \label{P:GroupOfFractions}
(i) The monoid~$\BBsp$ embeds in the
group~$\BBs$, and the latter is a group of
left fractions of~$\BBsp$, \ie, every element
of~$\BBs$ can be expressed as $\px\inv \py$
with~$\px, \py$ in~$\BBsp$. Moreover, every element
of~$\BBs$ can be expressed as $\fx\inv \bx\inv \by
\fy$ with~$\bx, \by$ in~$\Bip$ and $\fx, \fy$
in~$\Fp$.

(ii) The group~$\BBsp$ is torsion free.
\end{prop}

\begin{proof}
For~$(i)$, the monoid~$\BBsp$ satisfies Ore's
conditions on the left, \ie, it is cancellative and
any two elements admit a left lcm. The second
decomposition follows from Proposition~\ref{P:Zappa}
and the equality $\BBsp = \Bip \opp \Fp$.
Point~$(ii)$ follows as every torsion element in the
group of fractions of a monoid admitting lcm's is a
conjugate of a torsion element of the
monoid~\cite{Dha}. As $\BBsp$ has no torsion element
but~$1$, the same holds in~$\BBs$.
\end{proof}

Word reversing solves the word problem for the
group~$\BBs$.

\begin{lemm} \label{L:WordProblem}
A word~$\ww$ represents~$1$
in~$\BBs$ if and only if its double left
$\RRs$-reversing ends up with an empty word, where
double left reversing consists in left
reversing~$\ww$ into~$\uu\inv \vv$ with
$\uu, \vv$ positive, and then
left reversing~$\vv \uu\inv$.
\end{lemm}

\begin{proof}
Lemma~\ref{L:Convergence} guarantees that, for every
word~$\ww$, there exist positive words~$\uu,
\vv$ such that $\ww$ is left $\RRs$-reversible
to~$\uu\inv \vv$. Then $\ww$ represents~$1$ in~$\BBs$
if and only if $\uu$ and~$\vv$ represent the same
element of~$\BBs$, hence the same element of~$\BBsp$,
as $\BBsp$ embeds in~$\BBs$. Now, by definition of
completeness, the latter is true if and only if the
left reversing of~$\vv \uu\inv$ ends up with~$\e$.
\end{proof}

Then we have the following group version of
Lemma~\ref{L:Parabolic} for presentation of
subgroups. The point is that word
reversing solves the word problem without
introducing any~$\gx\gx\inv$ or~$\gx\inv\gx$.

\begin{lemm} \label{L:ParabolicBis}
Assume that $(\XX, \RR)$ is a positive presentation
that is complete for left reversing and such that
left reversing always terminates. Let $\XXo$ be a
subset of~$\XX$, and let~$\RRo$ be the set of all
relations $\vv \gx = \uu \gy$ in~$\RR$ with $\gx,
\gy \in \XXo$. If all words occurring in~$\RRo$
are $\XXo$-words, the subgroup of~$\Gr(\XX;
\RR)$ generated by~$\XXo$ admits the
presentation~$\Gr(\XXo; \RRo)$.
\end{lemm}

\begin{proof}
The hypotheses guarantee that an $\XXpm$-word
represents~$1$ in the group~$\Gr(\XX;
\RR)$ if and only if it can be transformed to~$\e$
by double left reversing. Now, as in the proof of
Lemma~\ref{L:Parabolic}, the hypotheses imply that
all words appearing in a (double) reversing from
an $\XXopm$-word are $\XXopm$-words. So, if
such a word is left $\RR$-reversible to~$\e$, it
is also left $\RRo$-reversible to~$\e$, and it
represents~$1$ in~$\Gr(\XXo; \RRo)$. 
\end{proof}

\begin{prop} \label{P:Subgroups}
The subgroup of~$\BBs$ generated by~$\XXs$ is (a
copy of) the braid group~$\Bi$, and the subgroup
generated by~$\XXa$ is (a copy of) Thompson's
group~$F$. These subgroups  generate~$\BBs$, and
their intersection is~$\{1\}$.
\end{prop}

\begin{proof}
The argument is the same as for the submonoids,
replacing Lemma~\ref{L:Parabolic} with
Lemma~\ref{L:ParabolicBis}. Then, by definition,
$\BBs$ is generated by the~$\ss i$'s and the~$\aa
i$'s, hence by the subgroups they generate
(henceforth identified with~$\Bi$ and~$F$). Assume
$\pz \in \Bi \cap F$. Every element of~$F$ is a left
fraction, so we have $\pz = \fx\inv \fxx$ for some
$\fx, \fxx$ in~$\Fp$. By Garside's theory, $\Bi$ is
both a group of left and of right fractions
of~$\Bip$, so we also have $\pz =
\bx \bxx\Inv$ for some $\bx, \bxx$ in~$\Bip$. We
deduce $\bx \fx = \bxx \fxx$ in~$\BBsp$, and the
uniqueness of the decomposition in~$\Fp \times \Bip$
(Proposition~\ref{P:Zappa}) implies $\bx = \bxx$
and $\fx = \fxx$.
\end{proof}

>From now on, we consider~$\Bi$ and~$F$ as
subgroups of~$\BBs$. For future use, we insist that
every element of~$\BBs$ can be represented by a word
in which the~$\aaaa i$ letters are gathered.

\begin{defi} \label{I:Tidy}
A $\wsa$-word is called {\it tidy} if it
consists of letters~$\aaa i$, followed by
letters~$\ssss j$, followed by letters~$\aa k$.
\end{defi}

Propositions~\ref{P:GroupOfFractions} implies:

\begin{coro} \label{C:Tidy}
Every element of~$\BBs$ admits a tidy representative.
\end{coro}

\section{The self-distributive structure on~$\Bs$}
\label{S:LD}

Besides their group structure, parenthesized braids
are equipped with another important algebraic
structure, involving the self-distributivity law. 

A non-trivial property of the braid group~$\Bi$ is
the existence of a binary operation that
obeys the self-distributivity law $x(yz) =
(xy)(xz)$. The importance of this exotic operation
originates from the fact that each element of~$\Bi$
generates a free subsystem with respect to the
self-distributive operation, a property directly
connected with the existence of a canonical ordering
of~$\Bi$ \cite{Dgd, Dgr}. In this section, we show
that the self-distributivity properties of~$\Bi$
extend to~$\Bs$, in an even stronger form as the
structure involves a second related operation that
has no counterpart in the case of ordinary braids.

As an application, we deduce that the groups~$\Bs$
and~$\BBs$ are isomorphic, \ie, we show that the
relations~$\RRs$ of Lemma~\ref{L:Relations} make a
presentation of~$\Bs$.

\subsection{The self-distributive bracket on~$\BBs$}


\begin{defi} \label{I:LDSystem}
An {\it LD-system} is a set equipped with a binary
operation $\xx,\yy \mapsto \LD\xx\yy$ satisfying the
left self-distributivity law
\begin{equation}\label{E:LD}
\LD\xx{\LD\yy\zz} =
\LD{\LD\xx\yy}{\LD\xx\zz}.
\end{equation}
An {\it augmented} LD-system, or ALD-system, is an
LD-system equipped with a second binary
operation~$\OP$ satisfying the mixed laws
\begin{equation} 
\label{E:ELD}
\LD\xx{\LD\yy\zz} = \LD{(\xx \OP \yy)}\zz
\text{\qquad and \qquad}
\LD\xx{\yy \OP \zz} = \LD\xx\yy \OP \LD\xx\zz.
\end{equation}
An LD-system is said to be {\it left cancellative} if
all left translations are injective, \ie, if $\LD x y
= \LD x z$ implies $y = z$; it is called a {\it rack}
\cite{FeR} if all left translations are
bijective, which means that there exists a binary
operation $x, y \mapsto \LDbar x y$ satisfying
$\LDbar x{\LD x y} = \LD x{\LDbar x y} = y$.
\end{defi}

A group equipped with $\LD\xx\yy = \xx \yy \xx\inv$,
$\LDbar\xx\yy = \xx\inv \yy\xx$ and $\xx \OP\yy =
\xx\yy$ is an augmented rack, always satisfying the
additional law
$\LD\xx\xx=\xx$. On the other hand, Artin's
group~$\Bi$ is an LD-system when equipped with the
operation
\begin{equation}\label{E:BraidBracket}
\LD\bx\by = \bx \opp \partial\by \opp \ss1
\opp \partial\bx\inv,
\end{equation}
where $\partial$ is the endomorphism that maps~$\ss i$
to~$\ss{i+1}$ for each~$i$. This operation can
be seen as a sort of twisted conjugacy, and there are
several ways of making the definition natural
\cite{Dgd}. The braid bracket is very different from
a group conjugacy in that
$\LD\bx\bx = \bx$ never holds. Observe that
there is no way to augment the LD-system~$\Bi$,
as, for instance, $\LD1{\LD11} = \LD \bx
1$  would imply $\partial\bx = \sss1 \sss2 \bx
\ss1$, which holds for no~$\bx$ in~$\Bi$.  

We shall see now that the braid bracket extends
to~$\BBs$, and, moreover, it can be augmented. We
begin with a preparatory result.

\begin{defi} \label{I:Shift}
We denote by~$\partial$ the {\it shift} that
maps $\ss i$ to~$\ss{i+1}$ and $\aa i$
to~$\aa{i+1}$ for each~$i$. 
\end{defi}

\begin{lemm} \label{L:dInjective}
The mapping~$\partial$ induces an injective endomorphism
of the group~$\BBs$ into itself.
\end{lemm}

\begin{proof}
As the shift mapping on positive integers is
injective, $\partial$ induces an isomorphism
of the group~$\Gr(\XXsa; \RRs)$ into its
image~$\Gr(\partial(\XXsa); \partial\RRs)$. Now the
explicit form of the relations in~$\RRs$ shows that
$\partial\RRs$ is included in~$\RRs$, and that the
criterion of  Lemma~\ref{L:ParabolicBis} is
satisfied by~$\partial(\XXsa)$ and~$\partial\RRs$. So the
subgroup of~$\BBs$ generated by~$\partial(\XXsa)$
admits the presentation~$\Gr(\partial(\XXsa);
\partial\RRs)$, and, therefore, $\partial$ is an
isomorphism of~$\BBs$ onto the latter subgroup.
\end{proof}

\begin{defi} \label{I:Bracket}
For~$\px, \py$ in~$\BBs$, we set
\begin{equation}
\LD\px\py = \px \opp \partial\py \opp \ss1
\opp \partial\px\inv, \qquad\text{and} \qquad
\px \OP \py = \px \opp \partial\py \opp \aa1.
\end{equation}
\end{defi}

\begin{prop}
The set~$\Bs$ equipped with the
operations~$\LD{}~$ and~$\OP$ is an ALD-system.
Furthermore, the bracket is left-cancellative, \ie,
$\LD\px\py = \LD\px\pz$ implies $\py = \pz$.
\end{prop}

\begin{proof}
A simple verification:
\begin{align*}
\LD{\LD\px\py}{\LD\px\pz}
&= \LD{(\px \opp \partial\py \opp \ss1 \opp
\partial\px\inv)}{\px \opp \partial\pz \opp \ss1
\opp \partial\px\inv}\\
&= \px \opp \partial\py \opp \ss1 \opp
\partial\px\inv \opp \partial\px \opp \partial^2\pz \opp
\ss2 \opp \partial^2\px\inv \opp \ss1 \opp
\partial^2\px \opp \ss2\inv \opp \partial^2\py\inv
\opp \partial\px\inv\\
&=^{(*)} \px \opp \partial\py  \opp \partial^2\pz
\opp \ss1 \ss2 \ss1 \ss2\inv \opp
\partial^2\py\inv \opp \partial\px\inv\\
&= \px \opp \partial\py \opp \partial^2\pz \opp \ss2
\ss1 \opp \partial^2\py\inv \opp \partial\px\inv
=^{(*)} \LD\px{\py \opp \partial\pz \opp \ss1 \opp
\partial\py\inv}
=\LD\px{\LD\py\pz}.
\end{align*}
The reason for~$(*)$ is that $\partial^2x$
commutes with~$\ss1$ for every~$x$. For left
cancellativity, $\LD\px\py =
\LD\px\pz$ implies $\partial\py \opp \ss1 = \partial\pz
\opp \ss1$, hence $\partial\py = \partial\pz$, and,
therefore, $\py = \pz$ by
Lemma~\ref{L:dInjective}.

Then, we find similarly:
\begin{align*}
\LD\px{\LD\py\pz}
&= \px \opp \partial\py \opp \partial^2\pz \opp
\ss2\ss1 \opp \partial^2\py\inv \opp \partial\px\inv
= \px \opp \partial\py \opp \aa1 \aa1\inv \opp
\partial^2\pz \opp \ss2\ss1 \opp \aa2 \aa2\inv
\opp \partial^2\py\inv \opp \partial\px\inv\\
&= (\px \OP \py) \opp \aa1\inv
\opp \partial^2\pz \opp \ss2\ss1\aa2 \opp
\partial(\px \OP \py)\inv
=  (\px \OP \py) \opp 
\partial\pz \opp \aa1\inv\ss2\ss1\aa2 \opp
\partial(\px \OP \py)\inv\\
\intertext{(because $\aa1 \opp \partial\pz = \partial^2\pz
\opp \aa1$ always holds)}
 &=  (\px \OP \py) \opp 
\partial\pz \opp \ss1 \opp
\partial (\px \OP \py)\inv
= \LD{ (\px \OP \py)}\pz,\\
\LD\px{\py \OP \pz} 
&= \px \opp \partial\py \opp \partial^2\pz \opp \aa2
\ss1 \opp \partial\px\inv
= \LD\px\py \opp \partial\px \opp \ss1\inv
\opp \partial^2\pz \opp \aa2 \ss1 \opp
\partial\px\inv\\
&= \LD\px\py \opp \partial\px  \opp \partial^2\pz 
\opp \ss1\inv\aa2 \ss1 \opp
\partial\px\inv
= \LD\px\py \opp \partial(\LD\px\pz) \opp
\partial^2\px \opp \ss2\inv \ss1\inv\aa2 \ss1 \opp
\partial\px\inv\\
&= \LD\px\py \opp \partial(\LD\px\pz) \opp
\partial^2\px \opp \aa1 \opp
\partial\px\inv
= \LD\px\py \opp \partial(\LD\px\pz) \opp \aa1
= \LD\px\py \OP \LD\px\pz,
\end{align*}
which completes the proof.
\end{proof}

The self-distributive structure so constructed will
be instrumental in the sequel.

\subsection{Diagram colouring}

We now come back to proving
that the relations of Lemma~\ref{L:Relations} make a
presentation of the group~$\Bs$. The point is to
establish that the canonical morphism of~$\BBs$
to~$\Bs$ is injective. We shall do it by
showing that, for any word~$\ww$, the class
of~$\ww$ in~$\BBs$ can be recovered from the isotopy
class of any diagram~$\DD_\tt(\ww)$, which depends
only on the class of~$\ww$ in~$\Bs$. To this end, we
appeal to diagram colourings.

The principle, which can be traced back at least to
Alexander, is to fix a nonempty set~$\SetOfCol$ (the
colours),  to attribute colours from~$S$ to the
initial positions in a braid diagram~$\DD$, and to
push the colours along the strands. If the colours
never change, the output colours are a permutation
of the input colours, and we do not gain much
information about the diagram. Now, assume that the
set of colours~$S$ is equipped with two binary
operations, say $\xx, \yy \mapsto \LD\xx\yy$ and
$\xx, \yy \mapsto \LDbar\xx\yy$---the notation is chosen to
suggest that $\LD\xx\yy$ and $\LDbar\xx\yy$ are images
of~$\yy$ under~$\xx$. We require that, when an
$\xx$-coloured strand crosses over a $\yy$-coloured
strand, then the colour of the latter becomes~$\LD\xx\yy$ or $\LDbar\xx\yy$ according to the orientation of
the crossing:
\begin{center}
\begin{picture}(45,18)(0,0)
\put(0,3){\includegraphics{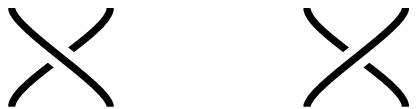}}
\put(0,15){$\xx$}
\put(11,15){$\yy$}
\put(30,15){$\yy$}
\put(41,15){$\xx$}
\put(-2,0){$\LD\xx\yy$}
\put(11,0){$\xx$}
\put(30,0){$\xx$}
\put(38,0){$\LDbar\xx\yy$}
\end{picture}
\end{center}
In this way, for each sequence of input colours and
each braid diagram~$\DD$, one obtains a sequence of
output colours, and some information
about~$\DD$ can be obtained by comparing the
input and output colours. One of the many facets of
the deep connection between braids and
self-distributivity is the following observation,
whose graphical verification is easy, and which
appears in different forms in~\cite{Brk, Joy,
Mat, Dgd, Dgr}:

\begin{lemm}\label{L:InvarianceColouring}
Assume that $\SetOfCol$ is a rack. Then $\SetOfCol$-colourings
are invariant under Reidemeister moves II and~III
in the sense that, for every diagram~$\DD$ and every
sequence of input colours, the corresponding output
colours depend only on the isotopy class of~$\DD$.
\end{lemm}

In order to control colourings in our current
framework, it is convenient to introduce coloured
trees. If $\DD$ is an ordinary $n$~strand braid
diagram, defining an $\SetOfCol$-colouring of~$\DD$
means attributing colours from~$\SetOfCol$ to the
$n$~input positions $1, \pp, n$, \ie, choosing a
sequence in~$\SetOfCol^n$. Propagating the colours along
the strands of~$\DD$ gives an output sequence that
lives in~$\SetOfCol^n$ again. Parenthesized braid diagrams
are similar, but the positions belong to~$\NNNs$
rather than to~$\NNN$, and they form a tree rather
than a sequence. Hence the objects to consider
are trees of $\SetOfCol$-coloured positions, \ie, 
{\it $\SetOfCol$-coloured trees}, defined to be trees
(of positions) in which colours from~$\SetOfCol$ are
attributed to the leaves. We shall use bold letters
like~$\ct$ for coloured trees.

\begin{defi} \label{I:Skeleton}
For $\xx$ in~$\SetOfCol$, we denote by~$\et_\xx$ the tree
with one single $\xx$-coloured node. For~$\ct$ an
$\SetOfCol$-coloured tree, we define the {\it
skeleton}~$\ct^\dag$ of~$\ct$ to be the uncoloured
tree~$\tt$ obtained by forgetting the colours
in~$\ct$; in this case, we say that $\ct$ is a
colouring {\it of~$\tt$}.
\end{defi}

Every $\SetOfCol$-coloured tree admits a unique
decomposition as a product of~$\et_\xx$ with~$\xx$
in~$\SetOfCol$. In particular, the sequence of
positions~$1, \pp, n$ with the colours~$\xx_1, \pp,
\xx_n$, as used for an ordinary
$\SetOfCol$-coloured $n$~strand braid diagram, corresponds
to the $\SetOfCol$-coloured right vine $\et_{\xx_1} \op
(\et_{\xx_2} \op \pp \op (\et_{\xx_n} \et) \pp
)$---as the last leaf encodes no position, we give
it no colour; if needed, we may assume that some
distinguished colour~$\xx_0$ is fixed and identify
an uncoloured tree with a tree uniformly
coloured~$\xx_0$.

\begin{figure}[htb]
\begin{picture}(105,35)(0,0)
\put(6,3){\includegraphics[scale=1]{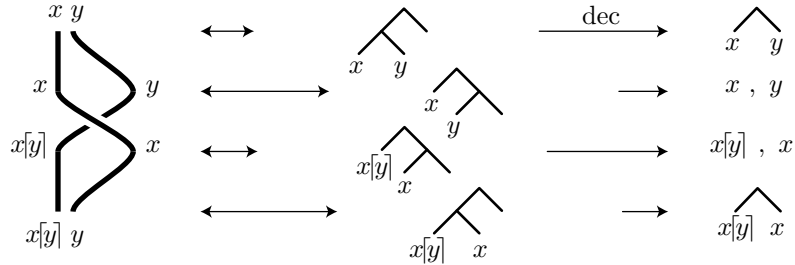}}
\put(5, 32){$x$}
\put(8, 32){$y$}
\put(3, 22){$x$}
\put(18, 22){$y$}
\put(0, 14){$\LD{x\!}{\!y\!}$}
\put(18, 14){$x$}
\put(2, 2){$\LD{x\!}{\!y\!}$}
\put(8, 2){$y$}
\put(45, 24.5){$x$}
\put(51, 24.5){$y$}
\put(55, 19.5){$x$}
\put(58, 16.5){$y$}
\put(46, 11){$\LD{x\!}{\!y\!}$}
\put(51.5, 8.5){$x$}
\put(53, 0){$\LD{x\!}{\!y\!}$}
\put(61.5, 0){$x$}
\put(76,31){$\Dec$}
\put(95, 27.5){$x$}
\put(101, 27.5){$y$}
\put(95, 22){$x~,~y$}
\put(93, 14){$\LD{x\!}{\!y\!}~,~x$}
\put(94, 3){$\LD{x\!}{\!y\!}$}
\put(101, 3){$x$}
\end{picture}
\caption{\smaller Correspondence between sets of
$\SetOfCol$-coloured positions and $\SetOfCol$-coloured trees:
here we start from $(1)$ and $(1, 1)$ coloured $x$ and $y$,
\ie, from the coloured tree
$(\et_x \op \et_y) \op \et$; then we go to $(1)$
and $(2)$ coloured $x$ and $y$, \ie, to $\et_x
(\et_y \et)$, \etc; on the right, we show the decomposition
of the trees, \ie, the subtrees under the right branch, the
last leaf excepted}
\label{F:ActionColouredTrees}
\end{figure}

Propagating $\SetOfCol$-colours along the strands of
a parenthesized braid diagram~$\DD$ amounts to defining a
partial action of~$\DD$ on $\SetOfCol$-coloured trees,
since, assuming that $\tt$ is the initial set of
positions in~$\DD$ and $\ttt$ is the final one, we
can associate with every $\SetOfCol$-colouring of~$\tt$ an
$\SetOfCol$-colouring of~$\ttt$
(Figure~\ref{F:ActionColouredTrees}):

\begin{defi} \label{I:Action}
For $\DD$ a parenthesized braid diagram with initial set of
positions~$\Pos\tt$ and $\ct$ an $\SetOfCol$-colouring
of~$\tt$, we denote by~$\ct \act \DD$ the
$\SetOfCol$-coloured tree obtained by propagating the
colours of~$\ct$ through~$\DD$. When $\DD$ has the
form~$\DD_\tt(\ww)$ for some word~$\ww$, we write $\ct \act \ww$ for~$\ct \act
\DD_\tt(\ww)$.
\end{defi}

It is easy to explicitly describe the action of~$\ss
i$ and~$\aa i$ on coloured trees.

\begin{lemm}\label{L:ActionPositionsTree}
Assume that $\ct$ is a coloured tree with $\Dec(\ct) =
(\ct_1, \pp, \ct_n)$. Then the coloured trees $\ct
\act \ss i$ and $\ct \act \aa i$ are defined for $i <
n$, and we have then
\begin{gather}
\label{E:ActionSOnColouredTree}
\Dec(\ct \act \ss i) = (\ct_1, \pp, \ct_{i-1},
\LDD{\ct_i}{\ct_{i+1}} ,
\ct_i, \ct_{i+2}, \pp, \ct_n),\\
\label{E:ActionAOnColouredTree}
\Dec(\ct \act \aa i) = (\ct_1, \pp, \ct_{i-1}, \ct_i
\op \ct_{i+1}, \ct_{i+2}, \pp, \ct_n),
\end{gather}
where $\LDD{\ct_i}{\ct_{i+1}}$ denotes the tree
obtained from~$\ct_{i+1}$ by replacing every
colour~$\xx$ with the corresponding colour
$\LD{\xx_1}{\LD{\pp}{\LD{\xx_p}{\xx}}\pp}$, where $\xx_1,
\pp, \xx_p$ form the left-to-right enumeration of the 
colours in~$\ct_i$.
\end{lemm}

\begin{proof}
First, we observe that the rules
of~\eqref{E:ActionSOnColouredTree}
and~\eqref{E:ActionAOnColouredTree} extend those
of~\eqref{E:ActionSOnTree}
and~\eqref{E:ActionAOnTree}: this is
natural, as, when we forget the colours, we must
find the previously defined action on families of
positions, \ie, on trees. So it only remains to look
at colours. For~\eqref{E:ActionAOnColouredTree}, the
result is clear as colours are not changed. As
for~\eqref{E:ActionSOnColouredTree}, the result of
applying~$\ss i$ is that each strand corresponding
to~$\ct_{i+1}$ goes under all strands corresponding
to~$\ct_i$, and it meets the latter from right to
left: the first one corresponds to the rightmost
position in~$\ct_i$, and the last one corresponds to
the leftmost position in~$\ct_i$. Applying the rule
for changes of colours at crossings, we deduce that
the strand with initial colour~$\xx$ eventually gets
the colour 
$\LD{\xx_1}{\LD{\pp}{\LD{\xx_p}{\xx}}\pp}$.
\end{proof}

\subsection{Using left cancellative LD-systems}

Lemma~\ref{L:InvarianceColouring} states that, if
$\SetOfCol$ is a rack, then, for each $\SetOfCol$-coloured
tree~$\ct$, the tree $\ct \act \DD$ depends on the
isotopy class of~$\DD$ only. It follows that, if two
words~$\ww, \ww'$ are
$\RRs$-equivalent and $\ct \act \ww$ and $\ct
\act \ww'$ are defined, the latter are equal.

In the sequel, we shall consider a more general
situation, namely when the set of
colours is a left cancellative LD-system, but
not necessarily a rack. In this case, all pairs of
colours need not be eligible for negative crossings:
we can still define~$\LDbar x y$ to be the unique
element~$z$ satisfying $\LD x z = y$ when it exists,
but the operation $\LDbar{}~$ need not be
everywhere defined. The following lemma gathers the
results we need:

\begin{lemm} \label{L:ColouringExist}
Let $\SetOfCol$ be a left cancellative
LD-system. Assume that $\ww_1, \pp, \ww_r$ are
words and $\tt$ is a tree such that $\tt \act
\ww_k$ exists for each~$k$. Then there exists at
least one colouring~$\ct$ of~$\tt$ such that
$\ct \act \ww_k$ exists for every~$k$.
\end{lemm}

\begin{proof}
If $\SetOfCol$ is a rack, {\it any} $\SetOfCol$-colouring
is convenient, as the colours can always be
propagated. When $\SetOfCol$ is only supposed to be a
left cancellative LD-system, we must be more
careful. First, we observe that, if the
word~$\ww$ is left $\RRs$-reversible to~$\www$,
and $\ct \act \www$ exists for some
$\SetOfCol$-coloured tree~$\ct$, then $\ct \act \ww$
exists as well, as can be checked by considering
the various cases---the point is that left
reversing creates no~$\sss i\ss i$. Hence, as
every word is left reversible to a
negative--positive word, it suffices to prove the
result when each~$\ww_k$ is such a word.
Moreover, positive words create no problem, so it
is even sufficient to consider the case when
each~$\ww_k$ is a negative word. Putting $\vv_k
= \ww_k\inv$, our problem is to prove
that, if $\vv_1, \pp, \vv_r$ are positive words,
then there exist $\SetOfCol$-coloured trees~$\ct_1,
\pp, \ct_r$ such that $\ct_k \act
\vv_k$ exists and is equal to some tree~$\ctt$
independent of~$k$. Now, by
Proposition~\ref{P:LcmB}, the elements
of~$\BBsp$ represented by~$\vv_1, \pp,
\vv_r$ admit a left common multiple, hence there
exist positive words~$\uu_1, \pp, \uu_r$ such
that the words $\uu_k \vv_k$ all are positively
$\RR$-equivalent (\ie, without introducing any
negative letter) to some positive word~$\ww$. Let
$\tt$ be a tree large enough to guarantee
that $\tt \act \ww$ exists, and let $\ct$ be any
$\SetOfCol$-colouring of~$\tt$. Put $\ct_k = \ct \act
\uu_k$. Then, by construction, $\ct_k \act
\vv_k$ exists and is equal to~$\ct \act \ww$ for
every~$k$.
\end{proof}

\begin{lemm} \label{L:Colouring}
Let $\SetOfCol$ be a left cancellative LD-system.
Assume that the parenthesized braid diagrams $\DD_\tt(\ww)$
and~$\DD_\tt(\ww')$ are isotopic. Then there exists at least one
$\SetOfCol$-colouring~$\ct$ of~$\tt$ such that
$\ct\act\ww$ and $\ct\act\www$ exist and are
equal. 
\end{lemm}

\begin{proof}
If $\SetOfCol$ is a rack, we can take for~$\ct$ any
$\SetOfCol$-colouring of~$\tt$. Then the colours can
be propagated without problem, \ie, $\ct\act\ww$
and~$\ct\act\www$ exist. The hypothesis that the
diagrams are isotopic implies in particular that
the final positions are the same, hence $\tt\act
\ww = \tt\act \www$ holds. On the other hand,
Lemma~\ref{L:InvarianceColouring} guarantees
that the sequences of output colours are the same
in both diagrams, \ie, the leaves
of~$\ct\act\ww$ and~$\ct\act\www$ have the same
colours. Hence
$\ct\act\ww$ and~$\ct\act\www$ are equal.

When $\SetOfCol$ is only supposed to be a
left cancellative LD-system, an arbitrary
$\SetOfCol$-colouring need not be convenient. Now,
the hypothesis that $\DD_\tt(\ww)$
and~$\DD_\tt(\ww')$ are isotopic implies that
there exists a finite sequence $\ww_1 =
\ww$, $\ww_1$, \pp, $w_r = \www$ such that, for
each~$k$, the diagram $\DD_\tt(\ww_{k+1})$ is
obtained from~$\DD_\tt(\ww_k)$ by one
Reidemeister move. By
Lemma~\ref{L:ColouringExist}, there exists an
$\SetOfCol$-colouring~$\ct$ of~$\tt$ such that $\ct
\act \ww_k$ is defined for each~$k$. Now, the
same argument as for
Lemma~\ref{L:InvarianceColouring} shows that
the final colours in two adjacent diagrams are
the same, hence in $\ct \act \ww$ and~$\ct \act
\www$, and we conclude as above.
\end{proof}

\subsection{Using $\BBs$-colourings}

As $\BBs$ equipped with its bracket is a left
cancellative LD-system, we can use it to colour
parenthesized braids. Here we use such colourings to
answer the pending question of whether the
relations~$\RRs$ present~$\Bs$. The key tool is a
certain function that associates with every
$\BBs$-coloured tree a specific element of~$\Bs$
constructed using the operation~$\OP$. 

\begin{defi} \label{I:Evaluation}
$(i)$ For~$\ct$ a $\BBs$-coloured tree, we
denote by~$\ev(\ct)$ the $\OP$-{\it evaluation} of~$\ct$,
\ie, the image of~$\ct$ under the mapping inductively
defined by
\begin{equation}
\ev(\et_x) = x
\qquad\text{and}\qquad
\ev(\ct \op \ctt) = \ev(\ct) \OP \ev(\ctt).
\end{equation}
The definition is extended to uncoloured trees by
identifying~$\et$ with~$\et_1$.

$(ii)$ For~$\ct$ a $\BBs$-coloured tree with
$\Dec(\ct) = (\ct_1, \pp, \ct_n)$, we put
\begin{equation}
\Evv(\ct) = \ev(\ct_1) \opp \partial\ev(\ct_2) \opp
\pp \opp \partial^{n-1}\ev(\ct_n).
\end{equation}
\end{defi}

For instance, for~$\tt$ the right vine of
size~$n+1$, we have $\ev(\tt) = \aa n \aa{n-1}
\pp \aa1$, while $\ev((\et\et)\et)$ is~$\aa1^2$.
We shall determine the action of the
generators~$\aa i$ and~$\ss i$ on the evaluation 
mapping~$\Evv$. First we begin with an auxiliary
result about ALD-systems.

\begin{lemm} \label{L:Auxi}
Assume that $\SetOfCol$ is an ALD-system. Then,
for all $\SetOfCol$-coloured trees~$\ct, \ctt$, we have
\begin{equation} \label{E:EvalEnhanced}
\ev(\LDD\ct\ctt) = \LD{\ev(\ct)}{\ev(\ctt)}.
\end{equation}
\end{lemm}

\begin{proof}
We use induction on the cumuled sizes of~$\ct$
and~$\ctt$. If both $\ct$ and~$\ctt$ have size~$1$,
the result follows from the definition
of~$\LDD\ct\ctt$ directly. Otherwise, the
definition gives
$$\LDD{(\cti\op\ctii)}{\ctt} =
\LDD{\cti}{\LDD{\ctii}{\ctt}}
\qquad \text{and} \qquad
\LDD{\ct}{\ctti \op \cttii} =
(\LDD{\ct}{\ctti}) \op (\LDD{\ct}{\cttii}).$$
Applying the evaluation morphism, we
deduce for $\ct = \cti \op \ctii$
\begin{align*}
\ev(\LDD\ct\ctt)
&= \ev(\LDD\cti{\LDD\ctii\ctt})
= \LD{\ev(\cti)}{\ev(\LDD\ctii\ctt)}\\
&= \LD{\ev(\cti)}{\LD{\ev(\ctii)}{\ev(\ctt)}}
= \LD{(\ev(\cti)\OP\ev(\ctii))}{\ev(\ctt)}
= \LD{\ev(\ct)}{\ev(\ctt)}
\end{align*}
using the induction hypothesis and the
first relation in~\eqref{E:ELD}.
Similarly, for $\ctt = \ctti\op \cttii$, we
find
\begin{align*}
\ev(\LDD\ct\ctt)
&= \ev((\LDD\ct\ctti) \op (\LDD\ct\cttii))
= \ev(\LDD\ct\ctti) \OP \ev(\LDD\ct\cttii)\\
&= \LD{\ev(\ct)}{\ev(\ctti)} \OP
\LD{\ev(\ct)}{\ev(\cttii)}
= \LD{\ev(\ct)}{\ev(\ctti) \OP \ev(\cttii)}
= \LD{\ev(\ct)}{\ev(\ctt)}
\end{align*}
using the induction hypothesis and the second
relation in~\eqref{E:ELD}.
\end{proof}

Then the following technical result is crucial, as
it shows that the mapping~$\Evv$ transforms the
action of diagrams on trees into a multiplication in
the group~$\BBs$.

\begin{lemm}\label{L:MainEval}
For~$\ct$ a $\BBs$-coloured tree~$\ct$ and~$\ww$ 	
a word such that $\ct \act \ww$ exists, we
have
\begin{equation}\label{E:MainEval}
\Evv(\ct \act \ww) = \Evv(\ct) \opp \cl\ww, 
\end{equation}
where $\cl\ww$ denotes the element of~$\BBs$
represented by~$\ww$.
\end{lemm}

\begin{proof}
For an induction, it is sufficient to establish
\eqref{E:MainEval} when $\ww$ consists of one single
letter~$\ss i$ or~$\aa i$. Let us assume
$\ev(\Dec(\ct))= (x_1, \pp, x_n)$, where $\ev((\ct_1,
\pp, \ct_n))$ stands for $(\ev(\ct_1), \pp,
\ev(\ct_n))$. First, we find
\begin{gather}
\label{E:Es}
\ev(\Dec(\ct \act \ss i)) = (x_1, \pp, x_{i-1},
\LD{x_i}{x_{i+1}}, x_i, x_{i+2}, \pp, x_n),\\
\label{E:Ea}
\ev(\Dec(\ct \act \aa i)) = (x_1, \pp, x_{i-1}, x_i
\OP x_{i+1}, x_{i+2}, \pp, x_n).
\end{gather}
Indeed, \eqref{E:Es} follows
from~\eqref{E:ActionSOnColouredTree}
using~\eqref{E:EvalEnhanced}, and \eqref{E:Ea}
follows from~\eqref{E:ActionAOnColouredTree}. Then
we find
\begin{align*}
\Evv(\ct \act \ss i)
&= x_1 \opp \pp \opp \partial^{i-2} x_{i-1} \opp
\partial^{i-1}(\LD{x_i}{x_{i+1}}) \opp \partial^i x_i \opp
\partial^{i+1} x_{i+2} \opp \pp \opp \partial^{n-1} x_n\\
&= x_1 \opp \pp \opp \partial^{i-2} x_{i-1} \opp
\partial^{i-1} x_i \opp \partial^i x_{i+1} \opp \ss i \opp
\partial^i x_i\inv
\opp \partial^i x_i \opp \partial^{i+1}x_{i+2} \opp \pp
\opp \partial^{n-1} x_n\\
&= x_1 \opp \pp \opp \partial^{i-2} x_{i-1} \opp
\partial^{i-1} x_i \opp \partial^i x_{i+1} \opp \ss i \opp
\partial^{i+1} x_{i+2} \opp \pp 
\opp \partial^{n-1} x_n\\
&= x_1 \opp \pp \opp \partial^{n-1} x_n \opp \ss i
= \Evv(\ct) \opp \ss i,
\end{align*}
as $\ss i \opp \partial^k x = \partial^k x \opp \ss i$ holds
for $k \ge i + 1$. For~$\aa i$, we find similarly
\begin{align*}
\Evv(\ct \act \aa i)
&= x_1 \opp \pp \opp \partial^{i-2} x_{i-1} \opp
\partial^{i-1}({x_i}\OP{x_{i+1}}) \opp
\partial^i x_{i+2} \opp \pp \opp \partial^{n-2} x_n\\
&= x_1 \opp \pp \opp \partial^{i-2} x_{i-1} \opp
\partial^{i-1} x_i \opp \partial^i x_{i+1} \opp \aa i \opp
\partial^i x_{i+2} \opp \pp \opp \partial^{n-2} x_n\\
&= x_1 \opp \pp \opp \partial^{n-1} x_n \opp \aa i
= \Evv(\ct) \opp \aa i,
\end{align*}
as $\aa i \opp \partial^k x = \partial^{k+1} x \opp \aa i$
holds for $k \ge i$.
\end{proof}

We are now able to conclude:

\begin{prop} \label{P:Presentation}
The groups~$\Bs$ and~$\BBs$ are isomorphic, \ie,
$(\XXsa, \RRs)$ is a presentation for the
group~$\Bs$ of parenthesized braids.
\end{prop}

\begin{proof}
Assume that $\ww$ and~$\ww'$ are words and there
is a tree~$t$ such that the diagrams~$\DD_t(\ww)$
and~$\DD_t(\ww')$ are isotopic. We have to prove
that~$\ww$ and~$\ww'$ are $\RRs$-equivalent,
\ie, they represent the same element of~$\BBs$.
Lemma~\ref{L:Colouring} guarantees that there
exists at least one $\BBs$-colouring~$\ct$
of~$\tt$ such that $\ct \act \ww$ and $\ct \act
\ww'$ are defined and equal. 
Now---this is the point---\eqref{E:MainEval}
implies that both $\ww$ and~$\www$
represent~$\Evv(\ct)\inv \opp
\Evv(\ct \act
\ww)$.
\end{proof}

All algebraic results about~$\BBs$ established in
Section~\ref{S:AlgebraicProp} are therefore
valid for~$\Bs$. In the sequel, we shall no longer
distinguish between~$\Bs$ and~$\BBs$, and use~$\Bsp$
for~$\BBsp$. In particular, we consider that $\Bi$
and~$F$ are included in~$\Bs$; the elements of~$F$
are called {\it Thompson elements}.

\subsection{Special decompositions}

Besides its group operation, the set~$\Bs$ is now
equipped with two binary operations, namely~$\LD{}~$
and~$\OP$. For each parenthesized braid~$\px$, the 
parenthesized braids that can be constructed
from~$\bx$ using these operations form a
sub-ALD-system of~$\Bs$. In particular, we can
start from the trivial braid~$1$, and introduce
what will be called special parenthesized braids. 

\begin{defi} \label{I:Special}
A braid (\resp a Thompson element, \resp a parenthesized braid)
is called {\it special} if it belongs to the closure
of~$\{1\}$ under~$\LD{}~$ (\resp under~$\OP$,
\resp under both~$\LD{}~$ and~$\OP$).
\end{defi}

For instance, $1$, $\ss1$, $\aa1$, and $\aa1 \ss2
\ss1 \aaa2$ are special parenthesized braids , as we can write
$$\ss1= \LD11, \quad 
\aa1 = 1\OP1, \quad
\aa1 \ss2 \ss1 \aaa2 =
\LD{\aa1}{\ss1}= \LD{(1\OP1)}{\LD11}.$$
We will see that every parenthesized braid admits
decompositions in terms of special parenthesized braids. The
following geometric characterization of special
parenthesized braids is crucial for uniqueness arguments. It
shows that special parenthesized braids are the
ones that produce themselves starting from a
right vine with trivial colours. To improve
readability, we skip some parentheses in trees
according to the convention that $\xx\yy\zz$
stands for~$\xx(\yy\zz)$; thus, for instance, a
right vine is denoted~$\et\et \pp \et$.

\begin{lemm} \label{L:CharactSpecial}
A parenthesized braid~$\pz$ is special if and only if it
admits an expression~$\pw$ such that each
sufficiently large
$\Bs$-coloured vine $(\et_1\et_1 \pp \et_1) \act \pw$
exists and has the form $\ct \et_1 \pp \et_1$. In
this case, all colours in~$\ct$ are special braids,
and we have $\pz =
\ev(\ct)$.
\end{lemm}

\begin{proof}
We first prove that the condition is necessary. As it
is true for~$\pz = 1$ with $\pw = \e$, it suffices
to prove that, if the condition is true for~$\pzi$
and $\pzii$, then it is for~$\LD\pzi\pzii$ and
$\pzi \OP \pzii$. So we assume that $\pw_i$ is an
expression of~$\pz_i$, that $(\et_1 \et_1 \pp
\et_1) \act \pw_i = \ct_i \et_1 \pp \et_1$ holds,
and, in addition, we have $\ev(\ct_i) = \pz_i$
and all colours in~$\ct_i$ are special braids.
Then
$\pwi \opp \partial\pwii \opp \ss1 \opp \partial\pwi\inv$
represents $\LD\pzi\pzii$, and, using the induction
hypothesis, we find
\begin{align*}
(\et_1 \et_1 \pp \et_1) 
&\act (\pwi \opp \partial \pwii \opp \ss1 \opp
\partial\pwi\inv)  = (\cti \et_1 \pp \et_1) \act (\partial
\pwii \opp
\ss1 \opp \partial\pwi\inv)\\ 
&= (\cti \ctii \et_1 \pp \et_1) \act (\ss1 \opp
\partial\pwi\inv)
=((\LDD\cti\ctii) \cti \et_1 \pp \et_1) \act
\partial\pwi\inv
=(\LDD\cti\ctii) \et_1 \pp \et_1. 
\end{align*}
Similarly, $\pwi \opp \partial\pwii \opp \aa1$
represents $\pzi \OP \pzii$, and we find
$$(\et_1 \et_1 \pp \et_1) \act (\pwi \opp \partial \pwii
\opp \aa1)
= (\cti \et_1 \pp \et_1) \act (\partial \pwii \opp
\aa1) 
= (\cti \ctii \et_1 \pp \et_1) \act \aa1
=(\cti \op \ctii) \et_1 \pp \et_1. 
$$

Conversely, by~\eqref{E:MainEval}, any equality
$(\et_1 \et_1 \pp \et_1) \act \pw = \ct \et_1 \pp
\et_1$ implies
$$\cl\pw
= \Evv(\et_1 \et_1 \pp \et_1) \opp \cl\pw
=  \Evv((\et_1 \et_1 \pp \et_1) \act \pw)
= \Evv(\ct \et_1 \pp \et_1) 
=\ev(\ct).$$
By definition, if the colours in~$\ct$ are special
braids (or, more generally, special parenthesized braids), the
evaluation~$\ev(\ct)$ is a special parenthesized braid. So, it
only remains to show that, whenever
$(\et_1 \et_1 \pp \et_1) \act \pw$ exists, then all
colours in the latter tree are special braids. Now we
can assume without loss of generality that $\pw$ is
tidy. Indeed, pushing the letters~$\aaa i$ to the
left and the letters~$\aa i$ to the right does
not change the negative crossings in the
associated braid diagram, and no obstruction may
appear. Now the hypothesis that
$(\et_1\et_1\et_1 \pp) \act \pw$ is defined implies
that there is no initial~$\aaa i$ in~$\pw$,
\ie, that $\pw$ consists of a braid word~$\vv$
followed by~$\aa i$'s. By \cite{Dgd},
Propositions~VI.5.8 and~5.12, if $\vv$ is a
$\ws$-word and $(\et_1\et_1 \pp\et_1)
\act \vv$ is defined, then the latter has the
form $\et_{\bz_1} \et_{\bz_2} \pp \et_{\bz_n}$ where
$\bz_1, \pp, \bz_n$ are special braids. The
subsequent~$\aa i$'s do not change the colours.
\end{proof}

We give now a complete description of special
Thompson elements. Note that, by definition of the
operation~$\OP$, such elements must be positive.

\begin{prop} \label{P:SpecialThompson}
(i) A Thompson element not equal to~$1$ is
special if and only if it has an expression
$\aa{i_1} \pp \aa{i_k}$ satisfying $i_{k+1} \ge
i_k - 1$ for each~$k$ and $i_r = 1$. This
expression is unique.

(ii) The mapping~$\ev$ establishes a one-to-one
correspondence between finite binary trees of
size~$n+1$ and special Thompson elements of
length~$n$. So, in particular, there are
$\frac1{n+1}{2n \choose n}$ special Thompson
elements of length~$n$.
\end{prop}

\begin{proof}
The existence of a decomposition as in~$(i)$ is true
for~$1$, and for $\fxi \OP
\fxii$ whenever it is for~$\fxi$ and~$\fxii$.
Hence it is true for every special Thompson
element. Conversely, if $\fx$ admits
an expression	$\ww$ as above, there is a unique
way of expressing~$\fx$ as $\fxi \OP \fxii$,
namely defining~$\fxi$ to be the element
represented by the largest prefix~$\wwi$ of~$\ww$
that finishes with~$\aa1$ if it exists, and~$1$
otherwise. Then $\fxi$ and~$\fxii$ have the same
syntaxic property as~$\fx$, and the parsing
continues.

Then, by definition, the mapping~$\ev$ establishes
a surjective mapping from trees to special
Thompson elements. To prove injectivity, we
observe that, for every tree~$\tt$, we have
\begin{equation}
(\et \et \et \pp ) \act \ev(\tt) = (\tt) \et \et \pp
\end{equation}
provided we start with a large enough vine, as
shows an easy induction on the size of~$\tt$.
Thus
$\ev(\tt)$ determines~$\tt$. This proves~$(ii)$,
and the uniqueness of the decomposition of~$(i)$
follows.
\end{proof}

\begin{lemm} \label{L:SpecialEval}
For each $\Bs$-coloured tree~$\ct$, we have 
\begin{equation} \label{E:SpecialPBraid}
\ev(\ct) = \pz_1 \opp \partial\pz_2 \opp \pp \opp
\partial^{n-1}\pz_n \opp \ev(\ct^\dag),
\end{equation}
where $(\pz_1, \pp, \pz_n)$ is the left-to-right
enumeration of the colours in~$\ct$.
\end{lemm}

\begin{proof}
First, for every special Thompson
element~$\fx$ of length~$n$ and every parenthesized braid~$\pz$,
we have
\begin{equation} \label{E:SpecialCommut}
\fx \opp \partial\pz = \partial^{1 + n}\pz \opp\fx.
\end{equation}
Indeed, the equality inductively follows from the
relation $\aa1 \opp \partial\pz = \partial^2\pz \opp \aa1$,
as the decomposition
of Proposition~\ref{P:SpecialThompson} guarantees
that, when pushing the letters~$\aa i$
of~$\fx$ to the right, one always meets
letters~$\aaaa k$ or~$\ssss k$ with~$k \ge i+1$.

Now we prove~\eqref{E:SpecialPBraid} using induction
on~$\ct$. The result is clear when
$\ct$ has size~$1$. For~$\ct = \cti \op \ctii$,
assuming that the colours in~$\ct_i$ are $\pz_{1, i},
\pp, \pz_{n_i, i}$ and using the induction
hypothesis, we find
$$\ev(\ct) = \pz_{1,1} \opp  \pp \opp
\partial^{n_1-1}\pz_{n_1,1} \opp \ev(\ct_1^\dag) \opp
\partial\pz_{1,2} \opp \pp \opp
\partial^{n_2}\pz_{n_2,2} \opp \partial\ev(\ct_2^\dag) \opp
\aa1.$$
By construction, $\ev(\ct_1^\dag)$ is a special
Thompson element of length~$n_i-1$.
Applying~\eqref{E:SpecialCommut} repeatedly, we
push $\ev(\ct_1^\dag)$ to the right, and obtain 
$$\ev(\ct) = \pz_{1,1} \opp  \pp \opp
\partial^{n_1-1}\pz_{n_1,1} \opp 
\partial^{n_1}\pz_{1,2} \opp  \pp \opp
\partial^{n_1 + n_2-1}\pz_{n_2,2} \opp
\ev(\ct_1^\dag) \opp \partial\ev(\ct_2^\dag)
\opp
\aa1,
$$
and \eqref{E:SpecialPBraid} follows using 
$\ev(\ct_1^\dag) \opp \partial\ev(\ct_2^\dag)
\opp
\aa1 = \ev(\ct^\dag)$.
\end{proof}

We can now express special parenthesized braids in terms of
special braids and Thompson elements.

\begin{prop} \label{P:DecompSpecial}
Every special parenthesized braid~$\pz$ admits a unique
decomposition
\begin{equation} \label{E:DecompSpecial}
\pz = \bx_1 \opp \partial\bx_2 \opp \pp \opp
\partial^{n-1}\bx_n
\opp \fz,
\end{equation}
where $\bx_1, \pp, \bx_n$ are special braids, and
$\fz$ is a special Thompson element of length~$n-1$.
\end{prop}

\begin{proof}
Let~$\pz$ be a special parenthesized braid. By
Lemma~\ref{L:CharactSpecial}, there exists a
$\Bs$-coloured tree~$\ct$, where all colours are special
braids, satisfying $\pz = \ev(\ct)$. Then
Lemma~\ref{L:SpecialEval} gives a decomposition of
the expected form. Next,
Proposition~\ref{P:Subgroups} first implies the
uniqueness of~$\fz$, as $\bx \opp
\fz = \bx' \opp \fz'$ implies $\bx\inv \bx' = \fz'
\fz\inv \in \Bi \cap F$. Then, when $\bx_1, \pp,
\bx_n$ are special braids, the product
$\bx_1 \opp \partial\bx_2 \opp \pp \opp \partial^{n-1}\bx_n$
determines each factor~$\bx_i$ as, by
Lemma~\ref{L:CharactSpecial} again, we have
$(\et_1 \pp \et_1) \act (\bx_1 \opp \pp \opp
\partial^{n-1}\bx_n) = \et_{\bx_1} \op \pp
\op \et_{\bx_n}$---note that we only use the easy
direction of Lemma~\ref{L:CharactSpecial}, and not
the more delicate converse that resorts to the
fine study of self-distributivity.
\end{proof}

Finally, we obtain canonical decompositions for
arbitrary positive parenthesized braids in terms of special
parenthesized braids, hence in terms of special braids and
special Thompson elements.

\begin{prop} \label{P:SpecialDecomp}
Every positive parenthesized braid~$\px$ admits two
unique decompositions:
\begin{gather} \label{E:SpecialDecomp}
\px = \pz_1 \opp \partial\pz_2 \opp \pp \opp
\partial^{p-1}\pz_p,\\
\label{E:SpecialDecompBis}\px = \bx_1 \opp
\partial\bx_2 \opp \pp \opp
\partial^{n-1}\bx_n
\opp \fz_1 \opp \partial\fz_2 \opp \pp \opp
\partial^{n-1}\fz_n,
\end{gather}
where $\pz_1, \pp, \pz_p$ are special parenthesized braids,
$\bx_1, \pp, \bx_n$ are special braids, and
$\fz_1, \pp, \fz_n$ are special Thompson
elements. 
\end{prop}

\begin{proof}
Let $\px$ be a positive parenthesized braid. By hypothesis,
$\px$ admits an expression~$\pw$ with no~$\sss i$
or~$\aaa i$. As $\pw$ contains no~$\sss i$, every
$\Bs$-colouring of a tree~$\tt$ such that $\tt \act
\pw$ is defined can be propagated along the strands
of the diagram~$\DD_\tt(\pw)$. Thus $\ct
\act \pw$ is defined for each $\Bs$-colouring~$\ct$
of~$\tt$, and \eqref{E:MainEval} then implies
$\px = \cl\pw = \Evv(\ct)\inv \opp \Evv(\ct \act
\pw)$. 

As $\pw$ contains no letter~$\aaa i$, we may
choose $\tt$ to be a right vine $\et \pp \et$, and
$\ct$ to be the corresponding colouring $\et_1 \pp
\et_1$. Then, by definition, we have $\Evv(\ct) =
1$, hence $\bx = \cl\pw = \Evv(\ct \act \pw)$.
Moreover, by construction, each colour in~$\ct
\act \pw$ belongs to the closure of~$\{1\}$
under the bracket operation, hence it is a
special braid. Then the
$\OP$-evaluation of the trees occurring in the
decomposition of~$\ct \act \pw$ are iterated
$\OP$-products of special braids, hence they are
special parenthesized braids. So, by definition, $\Evv(\ct \act
\pw)$ is a shifted product of special parenthesized braids, and
we obtain for~$\px$ a decomposition as
in~\eqref{E:SpecialDecomp}.

Now, if $\bx_1, \pp, \bx_n$ are special parenthesized braids,
Lemma~\ref{L:CharactSpecial} implies that, for
each~$k$, there exists an expression~$\pw_k$
of~$\pz_k$ satisfying
$(\et_1 \pp \et_1) \act \pw_k = (\ct_k) \et_1 \pp
\et_1$, where $\ct_k$ is a $\Bs$-coloured tree
satisfying $\ev(\ct_k) = \pz_k$. Provided the
initial  right vine is large enough, this implies
$$(\et_1 \et_1 \pp \et_1) \act (\pw_1 \opp \partial\pw_2
\opp \pp \opp \partial^{n-1}\pw_n) = (\ct_1) \pp (\ct_n)
\et_1 \pp \et_1.$$
This shows that the shifted product $\pz_1 \opp
\pp
\opp \partial^{n-1}\pz_n$ determines each tree~$\ct_k$,
hence each factor~$\pz_k$, thus proving the
uniqueness of the
decomposition~\eqref{E:SpecialDecomp}---we did not
prove here the (true) result that replacing~$\pw$
with an equivalent word~$\pww$ necessarily leads to
the same tree~$\ct$: this result is not
needed here, as we only use $\ev(\ct)$, which
is~$\px$ in any case.

Applying Proposition~\ref{P:DecompSpecial} to each
factor in~\eqref{E:SpecialDecomp} and
using~\eqref{E:SpecialCommut} to push the Thompson
factors to the right easily gives a decomposition
as in~\eqref{E:SpecialDecompBis}. For the
uniqueness of the latter, the same argument as
for Proposition~\ref{P:DecompSpecial} shows that
the braid part and the Thompson part are
determined, and that each special braid~$\bx_k$ is
determined by the shifted product $\bx_1 \opp
\pp
\opp \partial^{n-1}\bx_n$, so it only remains to verify
that the uniqueness of the special Thompson factors.
The latter follows from the equality
$$(\et_1\et_1 \pp \et_1) \act (\fz_1 \opp \partial\fz_2
\opp \pp \opp \partial^{n-1}\fz_n) = (\tt_1) \pp (\tt_n)
\et_1 \pp \et_1$$
for $\fz_k = \ev(\tt_k)$, again a consequence of
Lemma~\ref{L:CharactSpecial}.
\end{proof}

In the case of Thompson elements we have obtained
the following result, which provides a unique
normal form in~$\Fp$:

\begin{coro} \label{C:ThompsonDecomp}
Every positive Thompson element~$\fx$ admits a unique
decomposition
\begin{equation}
\fx = \fz_1 \opp \partial\fz_2 \opp \pp \opp
\partial^{p-1}\fz_p
\end{equation}
where $\fz_1, \pp, \fz_p$ are special Thompson
elements.
\end{coro}

By Proposition~\ref{P:GroupOfFractions}, every
parenthesized braid is a left fraction~$\px\inv \py$
with $\px, \py$ in~$\Bsp$, so another
consequence of Proposition~\ref{P:SpecialDecomp}
is:

\begin{coro}
Every parenthesized braid~$\px$ admits decompositions
\begin{gather} \label{E:SpecialDecompositionBis}
\px = \partial^{q-1}\pz'_q\Inv \opp \pp
\opp \partial\pz'_2\Inv \opp \pz'_1\Inv
\opp \pz_1 \opp \partial\pz_2 \opp \pp \opp
\partial^{p-1}\pz_p,\\
\px =  \partial^{n-1}\fz'_n\Inv \opp \pp
\opp \fz'_1\Inv
\opp \partial^{n-1}\bxx_n\Inv \opp \pp
\opp 
\bxx_1\Inv \opp \bx_1 \opp \pp \opp
\partial^{n-1}\bx_n \opp \fz_1 \opp \pp
\opp
\partial^{n-1}\fz_n,
\end{gather}
where $\pz_1, \pp, \pz'_q$
are special parenthesized braids, $\bx_1, \pp, \bxx_n$ are
special braids, and $\fz_1, \pp, \fz'_n$ are special
Thompson elements.
\end{coro}

\section{A linear ordering on~$\Bs$}
\label{S:Order}

Artin's braid group~$\Bi$ admits a distinguished
linear ordering that is compatible with
multiplication on one side and admits a number of
equivalent constructions \cite{Dgr}. On the other
hand, it is easy to construct on Thompson's
group~$F$ a linear ordering that is compatible with
multiplication on both sides. Merging these orderings
leads to ordering parenthesized braids.

\subsection{An ordering on~$\Fp$}

One can easily order~$F$ by attaching a piecewise
linear homeomorphism of~$[0,1]$ (or of the real
line) to each element and comparing the derivatives.
An equivalent construction involves trees.
We recall that, for~$\tt$ a tree, $\Dya(\tt)$
denotes the set of endpoints in the dyadic
decomposition of~$[0,1]$ attached to~$\tt$.

\begin{defi} \label{I:TreeOrder}
For $\tt, \ttt$ trees, we say that $\tt \prec
\ttt$ is true if $\Dya(\tt)$ follows~$\Dya(\ttt)$
in the lexicographical ordering.
\end{defi}

For instance, the sequences attached to
$\et(\et\et)$ and
$(\et\et)\et$ are $(0, \frac12, \frac34, 1)$ and
$(0, \frac14, \frac12, 1)$. The first
entries both are~$0$; the second entries are
$\frac12$ and $\frac14$, respectively: the former
is larger, so we declare $\et(\et\et) \prec
(\et\et)\et$.

\begin{lemm}
The relation~$\prec$ is a linear ordering
on trees. An alternative definition is: $\et
\prec
\tti\ttii$ is always true, and $\tti \op \ttii \prec
\ttti \op
\tttii$ is true if and only if $\tti
\prec \ttti$ is true, or $\tti = \ttti$ and $\ttii
\prec
\tttii$ are.
\end{lemm}

By Proposition~\ref{P:SpecialThompson}, the
evaluation mapping~$\ev$ establishes a one-to-one
correspondence between finite binary trees and
special Thompson elements. Moreover,
Corollary~\ref{C:ThompsonDecomp} shows that every
positive Thompson element admits a unique decomposition in
terms of special Thompson elements, hence in terms
of a sequence of trees. We can therefore carry the
tree ordering to~$\Fp$.

\begin{defi} \label{I:ThompsonOrder}
For $\fx, \fxx$ special Thompson elements,
we say that $\fx \lF^{sp} \fxx$ holds if and only
if we have $\ev\inv(\fx) \prec \ev\inv(\fxx)$.
For $\fx, \fxx$ in~$\Fp$, we say that $\fx
\lF \fxx$ holds if the (unique) special sequence
$(\fx_1, \pp, \fx_p)$ satisfying $\fx = \fx_1 \opp
\partial\fx_2 \opp \pp \opp \partial^{p-1}\fx_p$ is
lexicographically $\lF^{sp}$-smaller than the
special sequence $(\fxx_1, \pp, \fxx_q)$ satisfying
$\fxx = \fxx_1 \opp \partial\fxx_2 \opp \pp \opp
\partial^{q-1}\fxx_q$.
\end{defi}

For instance, we have $\aa2 \lF \aa1$, as the
special decomposition of~$\aa2$ is $1 \opp
\partial\aa1$, while $\aa1$ is special. Now
$\et(\et\et) \prec (\et\et)\et$ implies $1 = \ev(\et(\et\et))
\lF^{sp} \aa1 = \ev((\et\et)\et)$, and,
therefore, the sequence $(1, \aa1)$ is
lexicographically smaller than the
sequence~$(\aa1)$.

There is a canonical way of attaching
to each element~$\fx$ of Thompson's group~$F$ a
piecewise linear homeo\-morphism~$\hom\fx$ of
the unit inverval
\cite{CFP}---because of our conventions, we have
$\hom{\fx\fxx} =
\hom\fxx \circ \hom\fx$. The derivatives in~$\hom\fx$
make a finite sequence of dyadic numbers, \eg,
$(\frac12, 1, 2)$ in the case of~$\aa1$.

\begin{prop} \label{I:Homeo}
The relation~$\lF$ is a linear ordering on~$\Fp$. It
is compatible with multiplication on both sides. For
$\fx, \fxx$ in~$\Fp$, the relation $\fx \lF
\fxx$ holds if and only if the first derivative
not equal to~$1$ in~$\hom{\fx\inv \fxx}$ is smaller
than~$1$. 
\end{prop}

\begin{proof}
It is clear that $\lF$ is a linear ordering. 
The correspondence between~$\lF$ and the
homeomorphisms of~$[0,1]$ is as follows. If
$\ww$ is a positive $\wa$-word representing an
element~$\fx$, then $(\et\et\pp)
\act \ww$ is defined provided the initial vine
is large enough. Let $(\et\et\pp) \act
\ww = (\tt_1) \pp (\tt_p) \et \pp$. Then the
special decomposition of~$\fx$ is the shifted product
$\ev(\tt_1) \opp \partial\ev(\tt_2) \opp \pp$
Define~$\Dya(\fx)$ to be the union of the
sets~$\Dya(\tt_i)$ contracted from~$[0,1]$ to $[1 -
\frac1{2^{i-1}}, 1 - \frac1{2^i}]$ when $i$~varies.
Then $\fx \lF \fxx$ is equivalent to $\Dya(\fx)$
being larger than~$\Dya(\fxx)$ in the 
lexicographical order. Now the
homeomorphism~$\hom{\fx\inv \fxx}$ maps~$\Dya(\fx)$
to~$\Dya(\fxx)$, so the first divergence
between~$\Dya(\fx)$ and~$\Dya(\fxx)$ results in $\Dya(\fx)$
being declared larger if and only if the first
derivative~$\not=1$ in~$\hom{\fx\inv
\fxx}$ is less than~$1$.

Owing to the latter characterization, it is
clear that $\lF$ is compatible with multiplication on
the left. It is also compatible with multiplication
on the right, as the graph of~$\hom{\fx\fxx\fx\inv}$
is obtained from the graph of~$\hom\fxx$ by using~$\hom\fx$
to rescale the source and target intervals, which does not
change the fact that the graph diverges from the diagonal
downwards or upwards.
\end{proof}

For instance, the special decompositions of~$1,
\aa1$, and~$\aa2$ are $(\et, \et, \pp)$, $(\et\et,
\et, \et, \pp)$, and $(\et, \et\et, \et, \pp)$,
respectively. So we obtain $\Dya(1) = (0,
\frac12, \frac34, \frac78, \pp)$, 
$\Dya(\aa1) = (0, \frac14, \frac12, \frac34,
\pp)$,
and $\Dya(\aa2) = (0, \frac12,
\frac58, \frac34, \pp)$, hence $1 \lF \aa2 \lF
\aa1$.

\subsection{The ordering on~$\Bsp$}

As every element in~$\Bsp$ admits a
unique decomposition in terms of elements of~$\Bip$
and~$\Fp$, we deduce a linear order on~$\Bsp$
from any linear orders on~$\Bip$ and~$\Fp$. We
recall that $\Bi$ is equipped with a
distinguished linear ordering:

\begin{prop} \cite{Dgd, Dgr}
For $\bx, \bxx$ in~$\Bi$, say that $\bx \lB \bxx$
holds if and only if $\bx\inv \bxx$ admits an
expression in which the generator~$\ss i$ with
minimal index~$i$ occurs positively only, \ie, $\ss
i$ occurs but $\sss i$ does not. Then the
relation~$\lB$ is a linear ordering on~$\Bi$, and it
is compatible with multiplication on the left.
\end{prop}

\begin{defi} \label{I:PosOrder}
For $\px, \pxx$ in~$\Bsp$, we say that
$\px \lp \pxx$ holds if we have either $\bx \lB
\bxx$, or $\bx = \bxx$ and $\fx \lF \fxx$, where
$\px = \bx \fx$ and $\pxx = \bxx \fxx$ are the $\Bip
\times \Fp$-decompositions of~$\px$ and~$\pxx$.
\end{defi}

For instance, we have
$$\pp \lp \aa2 \lp \aa1 \lp \pp \lp \ss2 \lp
\ss1.$$ 
Indeed, we saw above that $\aa i \lF
\aa j$ holds for $i > j$ (in the case $i=1$,
$j=2)$. Then, we have $1 \lB \ss j$, hence $\aa i
\lp \ss j$ for all~$i, j$---and, more generally,
$\fx \lp \bx$ for all~$\fx$ in~$\Fp$ and $\bx$
in~$\Bip \setminus\{1\}$. Finally, $\ss i
\lF \ss j$ holds for $i > j$, as we have $\ss i \lB
\ss j$ since $\sss i \ss j$ is a braid word in which
the generator with smallest index, here~$\ss j$,
occurs positively and not negatively.

\begin{lemm}
The relation~$\lp$ is a linear order on~$\Bsp$,
compatible with left multiplication.
\end{lemm}

\begin{proof}
As both $\lB$ and $\lF$ are linear orders and the
$\Bip\times\Fp$-decomposition is unique, $\lp$ is
a linear order. To prove compatibility
with multiplication on the left, assume $\bx \fx
\lp \bxx \fxx$. Assume first $\bx \lB \bxx$. As
the braid ordering is compatible with
left multiplication, we have $\ss k \bx
\lB \ss k \bxx$ for every~$k$, hence $\ss k \opp
\bx \fx \lp \ss k
\opp \bxx \fxx$. On the other hand,
\eqref{E:Zappa} gives
\begin{equation}\label{E:CompatOrder}
\aa k \opp \bx \fx = \dbl_k(\bx)
\opp \aa{\bx\inv[k]} \fx
\quad\text{and}\quad
\aa k \opp \bxx \fxx = \dbl_k(\bxx)
\opp \aa{\bxx\Inv[k]} \fxx.
\end{equation}
To compare the braids $\dbl_k(\bx)$
and~$\dbl_k(\bxx)$, we consider $\dbl_k(\bx)\inv 
\dbl_k(\bxx)$. By construction, the latter is
$\dbl_{\bx\inv[k]}(\bx\inv \bxx)$. The
hypothesis $\bx \lB \bxx$ means that we can
represent~$\bx\inv \bxx$ by a braid diagram in
which the leftmost crossings all are positively
oriented. When we double a strand, the
latter property is preserved. So $\dbl_k(\bx) \lB
\dbl_k(\bxx)$ holds, and we deduce $\aa k \opp \bx
\fx \lp \aa k
\opp \bxx \fxx$. Hence, in this case,
$\px \opp \bx\fx \lp \px \opp \bxx\fxx$ holds
for every parenthesized braid~$\px$.

Assume now $\bx = \bxx$ and $\fx \lF \fxx$.
Then $\ss k \opp \bx\fx \lp \ss k \opp \bxx\fxx$
holds trivially for every~$k$. As for multiplication
by~$\aa k$, we use \eqref{E:CompatOrder} again:
$\bx = \bxx$ implies $\dbl_k(\bx) =
\dbl_k(\bxx)$, and $\fx \lF \fxx$ implies
$\aa{\bx\inv[k]} \fx \lF \aa{\bxx\Inv[k]}
\fxx$, because $\bx\inv[k] = \bxx\Inv[k]$ holds
and $\lF$ is compatible with multiplication on the
left. So, again, 
$\px \opp \bx\fx \lp \px \opp \bxx\fxx$ holds
for every parenthesized braid~$\px$.
\end{proof}

\subsection{The ordering on~$\Bs$}

As every parenthesized braid is a quotient of two
positive parenthesized braids, we can now easily deduce an
ordering on~$\Bs$ from the previous ordering
on~$\Bsp$.

\begin{defi} \label{I:Cone}
We denote by~$\Cone$ the set of all elements in~$\Bs$
that can be written as $\px\inv \pxx$ with $\px,
\pxx$ in~$\Bsp$ and $\px \lp \pxx$.
\end{defi}

\begin{lemm} \label{L:Cone}
The set~$\Cone$ is a positive cone, \ie, we
have $\Cone \opp \Cone \ince \Cone$ and $\Cone \cap \Cone\inv =
\emptyset$.
\end{lemm}

\begin{proof}
Consider two elements of~$\Cone$, say $\px_1\inv
\pxx_1$ and $\px_2\inv \pxx_2$ with $\px_i, \pxx_i$
in~$\Bsp$ and $\px_i \lp \pxx_i$ for $i = 1,2$. The
elements~$\pxx_1$ and~$\px_2$ admit a common left
multiple in~$\Bsp$, say $\py \pxx_1 = \pyy
\px_2$. Then we have $(\px_1\inv \pxx_1) \opp
(\px_2\inv \pxx_2) = (\py \px_1)\inv \opp (\pyy
\pxx_2)$. Using the compatibility of~$\lp$ with
left multiplication, we find
$\py \px_1 \lp \py \pxx_1 = \pyy \px_2 \lp
\pyy \pxx_2$, hence $(\px_1\inv \pxx_1) \opp
(\px_2\inv \pxx_2) \in \Cone$, and $\Cone \opp \Cone
\ince \Cone$.

Assume $\px \in \Cone \cap \Cone\inv$. Then we have
$1 = \px \opp \px\inv \in \Cone \opp \Cone$, hence $1
\in \Cone$ by the above result. So there must exist
$\bx, \bxx$ in~$\Bip$, and $\fx, \fxx$
in~$\Fp$ with $\bx \fx = \bxx \fxx$ and $\bx
\lB \bxx$, or $\bx = \bxx$ and $\fx \lF
\fxx$, contradicting the uniqueness of the
$\Bip \times \Fp$ decomposition in~$\Bsp$ in both
cases.
\end{proof}

\begin{defi} \label{I:Order}
For $\px, \pxx$ in~$\Bs$, we say that $\px < \pxx$
holds if $\px\inv \pxx$ belongs to~$\Cone$.
\end{defi}

For instance, we have $\ss2 < \aaa1 \ss1 \aa1 <
\ss1$. Indeed, we find $(\ss2)\inv (\aaa1 \ss1 \aa1)
= \aaa1 \sss3 \ss1 \aa1$, and $\ss3 \lB \ss1$ implies
$\ss3 \aa1 \lp \ss1 \aa1$. Similarly, we have 
$(\aaa1 \ss1 \aa1)\inv (\ss1) = 
\aaa1 \sss1 \aa1 \ss1 =\aaa1 \sss1 \ss2 \ss1 \aa2$,
and $\ss1 \lB \ss2\ss1$ implies $\ss1\aa1 \lp
\ss2\ss1 \aa2$.

\begin{prop} \label{P:Order}
The relation~$<$ is a linear ordering on~$\Bs$ that
is compatible with multiplication on the left, and
with the shift endomorphism~$\partial$. This linear
ordering extends the orders~$\lp$ on~$\Bsp$,
$\lB$ on~$\Bi$ and~$\lF$ on~$F$.
\end{prop}

\begin{proof}
Lemma~\ref{L:Cone} guarantees that $<$
is a partial order on~$\Bs$. This order is linear,
because $\lp$ is a linear order on~$\Bsp$, so, for
all~$\px, \pxx$ in~$\Bs$, either~$\px\inv \pxx$
or~$(\px\inv \pxx)\inv$, \ie, $\pxx\Inv \px$, belongs
to~$\Cone$. The order is compatible with
multiplication on the left by definition.
Then $\partial$ preserves the orders~$\lB$ and~$\lF$,
hence the order~$\lp$ on~$\Bsp$. This implies
$\partial\Cone \ince \Cone$, hence
$\px < \pxx$ implies, and, therefore, is equivalent
to, $\partial\px < \partial \pxx$.

Assume $\px, \pxx$ in~$\Bsp$ with $\px \lp \pxx$.
Then, by definition, $\px\inv \pxx$ belongs to~$\Cone$,
and, therefore, we have $\px < \pxx$ in~$\Bs$. As
$\lp$ is a linear ordering, the implication is an
equivalence.

Assume now $\bx, \bxx$ in~$\Bi$ with $\bx \lB \bxx$.
Then there exists a positive braid~$\bx_0$ such that
$\bx_0\bx$ and $\bx_0\bxx$ belong to~$\Bip$, and $\bx \lB
\bxx$ implies $\bx_0\bx \lB \bx_0\bxx$, hence $\bx_0\bx \lp
\bx_0\bxx$. Then $\bx\inv \bxx =
(\bx_0\bx)\inv (\bx_0\bxx)$ implies $\bx\inv \bx
\in \Cone$, hence $\bx < \bxx$. Once again, as
$\lB$ is a linear ordering, the implication is an
equivalence. Finally, for $\fx, \fxx$ in~$F$ with
$\fx \lF \fxx$, the same argument shows
that $\fx < \fxx$ holds in~$\Bs$. Hence $<$
restricted to~$F$ coincides with~$\lF$.
\end{proof}

\begin{coro}
The group~$\Bs$ is left-orderable. 
The group algebra~$\CCC[\Bs]$ has no zero divisor.
\end{coro}

\subsection{Syntaxic characterization}

We now describe the order on~$\Bs$ in terms of words.

\begin{defi} \label{I:Positive}
A $\wsa$-word is called {\it $\ss i$-positive} if
it contains~$\ss i$, but no~$\sss i$ or~$\ssss j$
with~$j < i$.
\end{defi}

\begin{prop} \label{P:Partition}
For~$\px$ a parenthesized braid not in~$F$, the
following are equivalent:

(i) We have $\px > 1$, \ie, $\px \in \Cone$;

(ii) There exists~$i$ such that $\px$ admits a tidy
$\ss i$-positive expression.
\end{prop}

\begin{proof}
Let $\px$ be an arbitrary parenthesized braid. By
Proposition~\ref{P:GroupOfFractions}, we can write
$\px = \fx\inv \bx \fxx$ with $\fx, \fxx \in F$ and
$\bx \in \Bi$. Then $\px \not\in F$ is equivalent to
$\bx \not= 1$. In that case, $\px \in \Cone$ is
equivalent to $\bx >_B 1$. By the results
of~\cite{Dgd}, the latter is equivalent to~$\bx$
admitting at least one $\ss i$-positive
expression. 
\end{proof}

The example of the word $\aa1 \ss2 \aa1\inv
\ss3\inv$, which is $\ss2$-positive but
represents~$1$ in~$\Bs$, shows that considering
tidy words is important. However, the case
of~$\ss1$ is particular, as we have:

\begin{prop}
If a parenthesized braid~$\px$ admits a $\ss1$-positive
expression, then $\px > 1$ holds.
\end{prop}

\begin{proof}
Let $\pw$ be a $\ss1$-positive word. We can
transform~$\pw$ into an equivalent tidy word by
pushing the letters~$\aa i$ to the right, and the
letters~$\aaa i$ to the left. The point is that,
in the process, the letters~$\ss1$ cannot vanish,
and no letter~$\sss1$ can appear. Indeed,
according to~\eqref{E:Zappa}, the rules for the
transformation are
\begin{equation*}
\aa k \ss i \mapsto \dbl_k(\ss i) \aa{\sss i[k]}
\qquad \text{and} \qquad
\ss i \aaa k \mapsto \aaa{\ss i[k]} \dbl_{\ss
i[k]}(\ss i).
\end{equation*}
By definition of the operation of doubling a
strand, the generator~$\ss i$ may be replaced
with~$\ss{i+1}$ in the case $k< i$, but this cannot
happen in the case~$i = 1$. Thus we always obtain
$\ss1$-positive words, and we finish with a tidy
$\ss1$-positive word.
\end{proof}

A direct consequence is:

\begin{prop} \label{P:CompatBracket}
For all~$\px, \py$ in~$\Bs$, one has $\px <
\LD\px\py$.
\end{prop}

\begin{proof}
By definition, we have $(\px)\inv \opp (\LD\px\py) =
\partial\px \opp \ss1 \opp \partial\py\inv$, an expression
with one~$\ss1$ and no~$\sss1$.
\end{proof}

\begin{coro}
Let $\px$ be an arbitrary element of~$\Bs$. Then the
closure of~$\{\px\}$ under the bracket operation is
a free LD-system.
\end{coro}

\begin{proof}
According to the so-called Laver's criterion
(\cite{Dgd}, Proposition~V.6.4), an
LD-system~$\SetOfCol$ with one generator is free if and
only if no equality of the form $\px = \LD\px
{\py_1} \LD{\pp}{\py_r}$ is possible in~$\SetOfCol$.
Now Proposition~\ref{P:CompatBracket} gives
$$\px < \LD\px{\py_1} <
\LD{\LD{\px}{\py_1}}{\py_2} <
\pp < \LD{\px}{\py_1}\LD{\pp}{\py_r}$$
for all~$\px, \py_1, \pp, \py_r$, hence $\px
\not= \LD\px {\py_1} \LD{\pp}{\py_r}$.
\end{proof}

\begin{ques}
Is the LD-system generated by~$1, \aa1, \aa2, \pp,
\aa{r-1}$ a free LD-system of rank~$r$? 
\end{ques}

\begin{rema}
There is no similar characterization of the
order~$\lF$ on~$F$ in terms of particular
decompositions. However, sufficient conditions exist.
Let us say that an $\wa$-word~$\fw$ is {\it $\aa
i$-positive} if $\fw$ contains~$\aa
i$, but no~$\aaa i$ or~$\aaa j$ with
$j < i$. Then an $\aa i$-positive word always
represents an element larger than~$1$, but,
conversely, $\aaa1 \aa2
\aa1$ is an example of an element larger than~$1$
that admits no
$\aa i$-positive expression.
\end{rema}

\subsection{The subword property}

The braid ordering is not compatible with
multiplication on the right, and, more generally,
there exists no linear ordering on~$\Bi$ that is
compatible with multiplication on both sides. So
the same holds for~$\Bs$, and $\Bs$ is not
bi-orderable. 

However, we shall now prove a partial
compatibility result involving conjugacy. In
general, a conjugate of an element~$\px$
satisfying~$\px > 1$ need not be larger than~$1$:
consider for instance $\ss1\sss2$ and its
conjugate~$\ss2\sss1$. We prove that this cannot
happen for~$\px$ in~$\Bip$.

We begin with a technical
result about the $\aa k$-conjugates of~$\ss i$
or, more generally, of any braid
$\dbl_i^p(\dbl_{i+1}^p(\ss i))$, or 
$\dbl_i^p\dbl_{i+1}^p(\ss i)$ for short, obtained
from~$\ss i$ by multiplying each strand by~$p+1$.

\begin{lemm} \label{L:Conjugate}
For all positive~$i, k$, and $p \ge 0$, there
exists~$i', k'$ and $e$ in~$\{0,1\}$ satisfying
\begin{equation} \label{E:DDouble}
\aa k \opp \dbl_i^p\dbl_{i+1}^p(\ss i)\opp \aaa k
= 
\aa{k'}^{-e}\opp
\dbl_{i'}^{p+e}\dbl_{i'+1}^{p+e}(\ss{i'})
\opp \aa{k'}^e.
\end{equation}
\end{lemm}

\begin{proof}
In the braid diagram $\dbl_i^p\dbl_{i+1}^p(\ss
i)$, the strands $i$ to~$i+p$ cross over the
strands $i+p+1$ to~$i + 2p + 1$. Hence
\eqref{E:DDouble} is clear for $k < i$ and
$k > i + 2p + 1$ with $e = 0$ and $i' = i$ or $i'
= i-1$. For $i
\le k \le i+p$, mult"plying by~$\aa k$ amounts to
doubling one more strand in the first block
of~$p$, so we have
$\aa k \opp \dbl_i^p\dbl_{i+1}^p(\ss i) = 
\dbl_i^{p+1}\dbl_{i+1}^p(\ss i) \opp \aa{k+p+1}$.
Then $\aa{k+p+1} \aaa k$ is $\aaa k \aa{k+p+2}$.
For the same geometric reason, we have
$\dbl_i^{p+1}\dbl_{i+1}^p(\ss i) \opp \aaa k = 
\aaa{k+p+2} \opp \dbl_i^{p+1}\dbl_{i+1}^{p+1}(\ss
i)$, which is~\eqref{E:DDouble} with $e = 1$,
$k' = k+p+1$ and $i' = i$. The computation is
similar for $i+p+1 \le k \le i + 2p + 1$, leading
now to $e = 1$, $k' = k$ and $i' = i$. 
\end{proof}

\begin{prop} \label{P:PropertyS}
For each parenthesized braid~$\px$ in~$\Bs$ and each~$i$, we have
$\px \ss i \px\inv > 1$.
\end{prop}

\begin{proof}
Write $\px = \fx\inv \bx \fxx$ with $\fx, \fxx \in
\Fp$ and $\bx$ in~$\Bi$. Then we have 
$$\px \ss i
\px\inv = \fx\inv \bx \fxx \ss i \fxx\Inv \bx\inv
\fx.$$ By Lemma~\ref{L:Conjugate}, we have $\fxx \ss
i \fxx\Inv = \fy\inv
\dbl_{i'}^p\dbl_{i'+1}^p(\ss{i'})
\fy$ for some~$\fy$ in~$\Fp$ and
some~$i', p$. Then we have $\bx \fy\inv =
\fyy\Inv \bxx$ for some $\fyy$ in~$\Fp$
and $\bxx$ in~$\Bi$, hence
$$\px \ss i \px\inv = \fx\inv \fyy\Inv \bxx
\dbl_{i'}^p\dbl_{i'+1}^p(\ss{i'}) \bxx\Inv \fyy
\fx.$$
By construction, the braid
$\dbl_{i'}^p\dbl_{i'+1}^p(\ss{i'})$ belongs
to~$\Bip$. By~\cite{Dgr}, Proposition~1.2.15, every
conjugate of a braid in~$\Bip$ is larger than~$1$.
Hence $\bxx \dbl_{i'}^p\dbl_{i'+1}^p(\ss{i'})
\bxx\Inv$ is a $\ss j$-positive braid for
some~$j$, and $\px \ss i \px\inv$ belongs
to~$\Cone$.
\end{proof}

\begin{coro}
For each parenthesized braid~$\px$, every parenthesized braid represented
by a word obtained from an expression of~$\px$ by
inserting letters~$\ss i$ is larger than~$\px$.
\end{coro}

\begin{proof}
It suffices to consider the addition of one~$\ss
i$, \ie, to compare elements of the form
$\px \py$ and $\px \ss i \py$. Now, we have $(\px
\py)\inv (\px \ss i \py) = \py\inv \ss i
\py$. By Proposition~\ref{P:PropertyS}, the latter
belongs to~$\Cone$.
\end{proof}

The previous property does
not extend to the letters~$\aa i$: for instance, we
have $\ss1 \sss1 = 1$ and $\ss1 \aa1
\sss1 = \ss1 \sss2 \sss1 \aa2 = \sss2 \sss1 \ss2
\aa2$, an expression that is $\ss1$-negative, hence
represents an element of~$\Cone\inv$. So, in this case,
inserting~$\aa1$ diminishes the element.

\subsection{Order and colourings}

The order on parenthesized braids can also be characterized in
terms of colourings by special braids. 

\begin{defi} \label{I:Col}
For~$\ct$ a $\Bs$-coloured tree, we denote
by~$\Col(\ct)$ the left-to-right enumeration of the
colours in~$\ct$. We denote by~$\Bi^{sp}$ the set of all
special braids.
\end{defi}

\begin{prop}
For all words~$\pw, \pww$, the
following are equivalent:

(i) We have $\cl\pw < \cl\pww$;

(ii) There exists a $\Bisp$-coloured tree~$\ct$
satisfying
\begin{equation} \label{E:Lexico}
\Col(\ct \act \pw) <^{Lex} \Col(\ct \act \pww) 
\quad\text{or}\quad
\Col(\ct \act \pw) = \Col(\ct \act \pww) 
~\text{and}~
(\ct \act \pw)^\dag \prec (\ct \act \pww)^\dag.
\end{equation}

(iii) For every $\Bi$-coloured tree~$\ct$ such
that $\ct \act \pw$ and $\ct \act \pww$
exist, \eqref{E:Lexico} holds.
\end{prop}

\begin{proof}
Assume $(ii)$. Put
$$(\bx_1, \pp, \bx_n) = \Col(\ct \act \pw), \quad
(\bxx_1, \pp, \bxx_n) = \Col(\ct \act \pww),
\quad (\bxxx_1, \pp, \bxxx_n) = \Col(\ct).$$
Then Lemma~\ref{L:SpecialEval} gives
\begin{gather*}
\ev(\ct \act \pw) = \bx_1 \opp \pp \opp
\partial^{n-1}\bx_n \opp \ev((\ct \act
\pw)^\dag),\\
\ev(\ct \act \pww) = \bxx_1 \opp \pp \opp
\partial^{n-1}\bxx_n \opp \ev((\ct \act
\pww)^\dag),\\
\ev(\ct) = \bxxx_1 \opp \pp \opp
\partial^{n-1}\bxxx_n \opp
\ev(\ct^\dag).
\end{gather*}
Next, \eqref{L:MainEval} gives
$\cl\pw = \ev(\ct)\inv \opp \ev(\ct \act \pw)$,
hence
$$\cl\pw\Inv \opp \cl\pww = 
\ev((\ct \act \pww)^\dag)\inv 
\opp \partial^{n-1}\bx_n\inv
\opp \pp \opp \partial\bx_2\inv
\opp \bx_1\inv \opp \bxx_1 \opp
\partial\bxx_2 \opp
\pp \opp \partial^{n-1}\bxx_n \opp 
\ev((\ct \act \pww)^\dag).$$
If $\Col(\ct \act \pw) <^{Lex} \Col(\ct \act
\pww)$ holds, there exists~$k$ such that
$\bx_i = \bxx_i$ holds for~$i < k$,
and $\bx_k < \bxx_k$ holds. As $\bx_k$
and $\bxx_k$ are special braids, this implies that
$\bx_k\inv \opp \bxx_k$ is $\ss1$-positive,
hence $\cl\pw\Inv \opp \cl\pww$ is $\ss
i$-positive, and $(i)$ is true.
On the other hand, if $\Col(\ct \act \pw)$
and $\Col(\ct \act \pww)$ coincide, there remains
$\cl\pw\Inv \opp \cl\pww = \ev((\ct \act
\pw)^\dag)\inv\opp \ev((\ct \act \pww)^\dag)$,
and, by definition, $(\ct \act \pw)^\dag \prec (\ct
\act \pww)^\dag$ implies $\ev((\ct \act \pw)^\dag)
\lF \ev((\ct \act \pww)^\dag$, hence $\cl\pw <
\cl{\pww}$. So $(ii)$ implies~$(i)$.

Assume now~$(iii)$. For~$\tt$ a large enough
tree, $\tt \act \ww$ and~$\tt \act \www$ are
defined, hence, by Lemma~\ref{L:ColouringExist},
there exists at least one $\Bs$-coloured
tree~$\ct$ such that $\ct \act \pw$ and $\ct \act
\pww$ exist. Hence~$(ii)$ holds.

Finally, assume that $(iii)$ fails. By the
argument above, there exists~$\ct$ such that $\ct
\act \pw$ and $\ct \act \pww$ exist and
\eqref{E:Lexico} fails. Because
$<^{Lex}$ and $\prec$ are linear orders, this
implies that either $\ww$ and~$\www$ are
equivalent, or \eqref{E:Lexico} with $\pw$
and~$\pww$ exchanged is true. We saw above that
this implies $\cl\pw > \cl\pww$. So, in any
case, $(i)$ fails.
\end{proof}

\section{Homeomorphisms of a punctured sphere}
\label{S:MCG}

Artin's braid group~$B_n$ can be realized as the
mapping class group of a disk with $n$~punctures
\cite{Bir}, and the induced action on the
fundamental group gives Artin's representation
of~$B_n$ in the automorphisms of a rank~$n$ free
group. In this section, we prove similar results for
the group~$\Bs$. We observe that
$\Bs$ can be mapped to the mapping class group of a
sphere with a Cantor set of punctures, and deduce
that $\Bs$ embeds in the groups of automorphisms of
a free group of countable rank using the ordering
of~Section~\ref{S:Order}.

\subsection{The mapping class group of a sphere
with a Cantor set of punctures}

We aim at mapping~$\Bs$ into the homeomorphisms of
a punctured space. As $\Bs$ includes~$\Bi$, disks
with infinitely many punctures are to be expected.
Moreover the tree-like structure of~$\Bs$ should
make it natural to meet the Cantor set. A suitable
choice is to collapse the boundary of the disk, \ie,
to start with a $2$-sphere, and to remove a Cantor
set of punctures. Note that the complement of a Cantor set
consists of a countable collection of open intervals
naturally indexed by dyadic numbers. 

\begin{figure}[htb]
\begin{picture}(113,42)(0, 0)
\put(0,0){\includegraphics{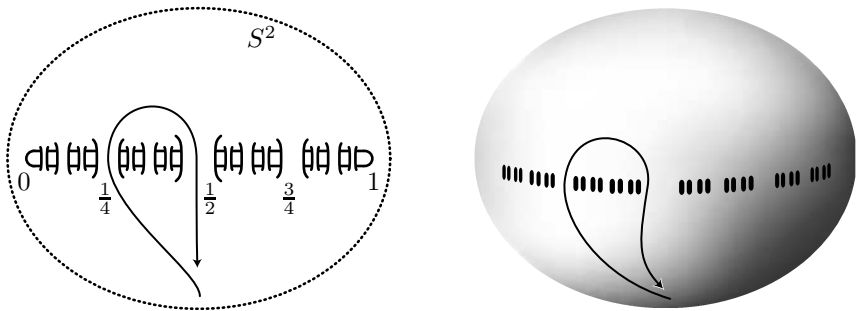}}
\put(32,35){$S^2$}
\put(1.5,16){$0$}
\put(12,14){$\frac14$}
\put(26,14){$\frac12$}
\put(36.5,14){$\frac34$}
\put(48,16){$1$}
\end{picture}
\caption{\smaller The space~$\SK$: a sphere with a
Cantor set removed from the equator, or, equivalently,
two hemispheres connected by a countable family of
bridges indexed by dyadic numbers; the loop represents
the element~$x_{1,1}\inv x_1$ of the fundamental group:
it starts from the South pole, crosses the bridge
at~$\frac14$ to the North hemisphere, and returns to the
South pole by the bridge at~$\frac12$.}
\label{F:SCantor}
\end{figure}

\begin{defi} (Figure~\ref{F:SCantor}) \label{I:Sphere}
We fix a real number~$\rho$ in~$(0, 1)$---for instance
$\rho = 1/3$---and we denote by~$\Can$ the Cantor
subset of~$[0, 1]$ obtained by iteratively removing
the median intervals of size~$\rho^k$. We define $\SK$ to be the topological space obtained
from the disk of diameter $[-\rho,1+\rho]$
in~$\RRR^2$ by removing the points of~$\Can$ and
collapsing the outer circle.
\end{defi}

We denote by~$\MCG(\SK)$ the mapping class group
of~$\SK$, \ie, the group of all homeomorphisms of~$\SK$
up to isotopy. As in the case of a finite set of punctures, a
continuous motion in the disk that maps~$\Can$ to
itself determines an element of~$\MCG(\SK)$.
Imitating the standard constructions, we can
define elements of~$\MCG(\SK)$ corresponding to
Dehn's half-twists on the one hand, and to
Thompson's piecewise linear homeomorphisms on
the other hand. 

\begin{figure}[htb]
\begin{picture}(45,26)(0,0)
\put(0,3){\includegraphics{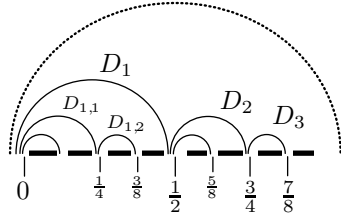}}
\put(7,12){$\scriptstyle D_{1,\!1}$}
\put(13,9.5){$\scriptstyle D_{1,\!2}$}
\put(1,0){$0$}
\put(21,-0.5){$\frac 1 2$}
\put(31,-0.5){$\frac 3 4$}
\put(36,-0.5){$\frac 7 8$}
\put(11,1){$\scriptstyle\frac 1 4$}
\put(16,1){$\scriptstyle\frac 3 8$}
\put(26,1){$\scriptstyle\frac 5 8$}
\put(12,17){$D_1$}
\put(28,12){$D_2$}
\put(35,10){$D_3$}
\end{picture}
\caption{\smaller The disks~$D_s$: essentially,
$D_s$ is the disk based on~$s$ and its immediate
successor in the lexicographical ordering: for
instance, $D_1$ is essentially the disk with
diameter~$[0, \frac 1 2]$, and
$D_{1,1}$ is essentially the disk with diameter~$[0, \frac
1 4]$; the adjustments guarantee that the disks $D_{s, i}$
are disjoint and nested in~$D_s$}
\label{F:DisksDs}
\end{figure}

\begin{figure}[htb]
\begin{picture}(100,26)(0,0)
\put(3,0){\includegraphics{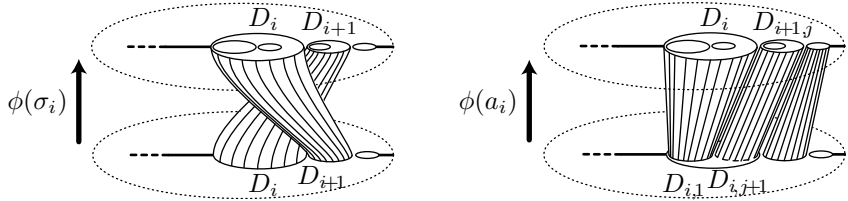}}
\put(-5,12){$\phi(\ss i)$}
\put(55,12){$\phi(\aa i)$}
\put(27,23){$D_i$}
\put(34,22.5){$D_{i+1}$}
\put(27,1){$D_i$}
\put(33.5,2){$D_{i\!+\!1}$}
\put(87,23){$D_i$}
\put(94,22.5){$D_{i\!+\!1,\!j}$}
\put(82,1){$D_{i,\!1}$}
\put(88,1.5){$D_{i\!,j\!+\!1}$}
\end{picture}
\caption{\smaller Homeomorphisms of~$\SK$ associated
with~$\ss i$ and~$\aa i$: a Dehn half-twist, and a
dilatation--contraction}
\label{F:Homeomorphisms}
\end{figure}

\begin{defi} \label{I:Disk}
$(i)$ (Figure~\ref{F:DisksDs}) Let~$s$ be a finite
sequence of positive integers, say $s = (i_1, \pp,
i_r)$. Put $\rho_s:= 2^{-i_1 - \cdots -i_r} \rho$. Then
$D_s$ is defined to be the (image in~$\SK$ of the) disk with
diameter
$$[ 0.1^{i_1-1}01^{i_2-1}0\dots 01^{i_r-1} - \rho_s\ ,\  
0.1^{i_1-1}01^{i_2-1}0\dots 01^{i_{r-1}-1}01^{i_r} -
\rho_s/2]$$
(referring to the dyadic expansion of rationals; $\rho$ is
the constant used in the realization of the Cantor
set~$\Can$, \eg, $1/3$).

$(ii)$ (Figure~\ref{F:Homeomorphisms})
For~$i \ge 1$, we define $\phi(\ss i)$ to be the class
in~$\MCG(\SK)$ of a clockwise half-turn (with rescaling)
that exchanges $D_i$ and~$D_{i+1}$ and is the identity on
all other~$D_j$'s. We define $\phi(\aa i)$ to  be the class
in~$\MCG(\SK)$ of a motion that fixes~$D_j$ for~$j < i$,
dilates~$D_{i,1}$ to~$D_i$, translates~$D_{i, j+1}$
to~$D_{i+1, j}$ for every~$j$, and contracts~$D_j$
to~$D_{j+1}$ for~$j > i$.
\end{defi}

An immediate verification shows that
all relations in~$\RRs$ induce isotopies, so we
have:

\begin{lemm} \label{L:MCGr}
The mapping~$\phi$ induces a morphism of~$\Bs$
into~$\MCG(\SK)$.
\end{lemm}

\subsection{Action on the fundamental group}

The homeomorphisms of~$\SK$ induce automorphisms of
its fundamental group, and those coming from the
elements of~$\Bs$ can be described explicitly.
We first identify~$\pi_1(\SK)$. 

\begin{figure}[htb]
\begin{picture}(144,46)(0, 0)
\put(0,0){\includegraphics{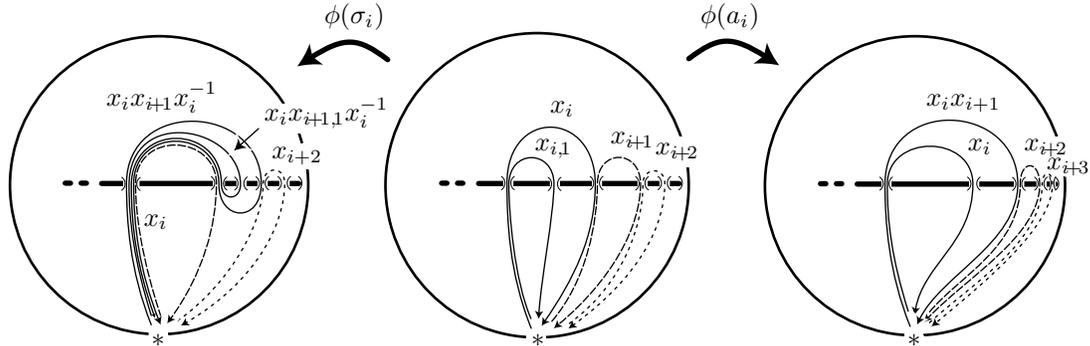}}
\put(19,0){$*$}
\put(69.5,0){$*$}
\put(119.5,0){$*$}
\put(42,43){$\phi(\ss i)$}
\put(92,43){$\phi(\aa i)$}
\put(72,31){$x_i$}
\put(70,26){$x_{i,\!1}$}
\put(80,27){$x_{i\!+\!1}$}
\put(86,25.5){$x_{i\!+\!2}$}
\put(13,32){$x_i x_{i\!+\!1} x_i\inv$}
\put(34,30){$x_i x_{i\!+\!1,\!1} x_i\inv$}
\put(35,25){$x_{i+2}$}
\put(18,16){$x_i$}
\put(122,32){$x_i x_{i+1}$}
\put(127.5,26){$x_i$}
\put(135,26){$x_{i\!+\!2}$}
\put(138,23.5){$x_{i\!+\!3}$}
\end{picture}
\caption{\smaller Generators of~$\pi_1(\SK)$,
and action of $\phi(\ss i)$ and $\phi(\aa i)$ on
these generators}
\label{F:ActionPiOne}
\end{figure}

\begin{defi} \label{I:Loop}
(Figures~\ref{F:SCantor} and~\ref{F:ActionPiOne})
For~$s$ a finite nonempty sequence of positive integers, we
define~$x_s$ to be the class in~$\pi_1(\SK)$ of a
loop that starts from the South pole of~$\SK$,
reaches the South pole of~$D_s$, turns
around~$D_s$ clockwise, and returns to the South
pole of~$\SK$. We define $x_s$ to be~$1$ for $s$ the empty
sequence. 
\end{defi}

\begin{lemm} \label{L:Base}
The fundamental group of~$\SK$ is the free
group~$\FGs$ based on the~$x_s$'s.
\end{lemm}

\begin{proof}
As $\SK$ is open in~$S^2$, a loop, which is compact,
may cross the equator only finitely many times.
So, in order to prove that $\pi_1(\SK)$ is generated
by the~$x_s$'s, it is sufficient to show that, for
every sequence~$s$, the loop~$\g_s$ that starts from
the South pole, crosses the equator at the left
of~$0$ and returns to the South hemisphere by the
bridge immediately at the right of~$D_s$ can be expressed as
a product of~$x_s$'s. Indeed, as $S^2$ has no
boundary, the loop crossing near~$0$ and returning
near~$1$ is trivial, and, if we can obtain~$\g_s$,
then, by using loops of the form $\g_s\inv \g_{s'}$, we
obtain every loop crossing the equator twice, and,
from there, every loop crossing the equator finitely
many times. Now, one easily checks that, for $s =
(i_1, \pp, i_r)$, one can take for~$\g_s$ any loop
representing
$$(x_1 x_2 \pp x_{i_1-1})
(x_{i_1,1} x_{i_1,2} \pp x_{i_1,i_2-1})
\pp (x_{{i_1},\pp,i_{r-1},1} 
x_{{i_1},\pp,i_{r-1},2}\pp
x_{{i_1},\pp,i_{r-1},i_r-1}).$$

It remains to show that the~$x_s$'s form a
free family. Assume that we have a relation
in~$\pi_1(\SK)$, say $w(x_{s_1}, \pp, x_{s_n}) = 1$
with~$w$ a freely reduced word. If the disks
$D_{s_1}$, \pp, $D_{s_n}$ are pairwise disjoint,
collapsing each of them to a point induces a
surjective homomorphism of the subgroup
of~$\pi_1(\SK)$ generated by~$x_{s_1}$,
\pp, $x_{s_n}$ onto the fundamental group of a disk
with $n$~punctures. The latter is a free group of
rank~$n$, so $w$ must be trivial.

Assume now that some disk~$D_{s_i}$ includes another
disk~$D_{s_j}$. This means that $s_i$ is a prefix
of~$s_j$. For each such~$i$, we define $y_i = 
x_{s_i, 1} x_{s_i, 2} \pp x_{s_i, p_i}$, where
$p_i$ is the minimal~$p$ such that $(s_i, p)$ is a
prefix of no other index~$s_j$. Note that the process
creates no new inclusion. Let~$\varphi$ be the
result of collapsing all~$x_{s_i, p}$'s with~$p >
p_i$. By construction, we have $\varphi(x_{s_i})
= y_i$, and, therefore, $w(x_{s_1},
\pp, x_{s_n}) = 1$ implies $w(y_1, \pp, y_n) = 1$.
Now, for each~$i$, the variable~$x_{s_i, p_i}$
occurs in~$y_i$ only, and the disks~$D_{s_i,
p_i}$ are disjoint. Then the same argument as
above shows that $w$ must be trivial.
\end{proof}

The homeomorphisms of~$\SK$ induce automorphisms of
its fundamental group~$\FGs$, and we obtain a
morphism of~$\MCG(\SK)$ into~$\Aut(\pi_1(\SK))$,
\ie, into~$\Aut(\FGs)$. 

\begin{prop} \label{P:AutoFree} \label{I:Psi}
Let~$\psi$ denote the composition of the above morphism
of~$\MCG(\SK)$ to~$\Aut(\FGs)$ with the morphism~$\phi$
of~$\Bs$ to~$\MCG(\SK)$. Then $\psi$ maps~$\Bs$
into~$\Aut(\FGs)$, and we have
\begin{gather} 
\label{E:Artins}
\psi(\ss i): 
\quad x_{j,s} \mapsto x_{j,s}
\text{\quad for $j \not= i, i+1$},
\qquad x_{i,s} \mapsto x_i x_{i+1, s} x_i\inv,
\qquad  x_{i+1, s} \mapsto x_{i,s},\\
\label{E:Acax}
\psi(\aa i): 
\begin{cases}
\quad x_{j,s} \mapsto x_{j,s}
\text{\ for $j < i$},
\quad x_{j,s} \mapsto x_{j+1,s}
\text{\ for $j > i$},
\\
\quad x_i \mapsto x_i x_{i+1},
\qquad x_{i, 1, s} \mapsto x_{i,s},
\qquad x_{i, j+1, s} \mapsto x_{i+1,j, s}
\text{\ for $j \ge 2$}.
\end{cases}
\end{gather}
\end{prop}

\begin{proof}
That $\psi$ is a morphism follows from the
construction---or from a direct verification, once
the explicit formulas for~$\psi(\ss i)$
and~$\psi(\aa i)$ are known. The latter can be read
in Figure~\ref{F:ActionPiOne}.
\end{proof}

\subsection{Determining the automorphism}

Once the automorphisms attached with~$\ss i$
and~$\aa i$ are known, we can determine the
automorphism of~$\FGs$ associated with any~$\px$
in~$\Bs$ by composing the automorphisms
associated with the successive letters of any
word representing~$\px$. Here we give an
alternative description involving
$\FGs$-coloured trees, \ie, finite binary trees in
which the leaves wear colours from~$\FGs$.

\begin{figure} [htb]
\begin{picture}(146,26)(0, 0)
\put(5,0){\includegraphics{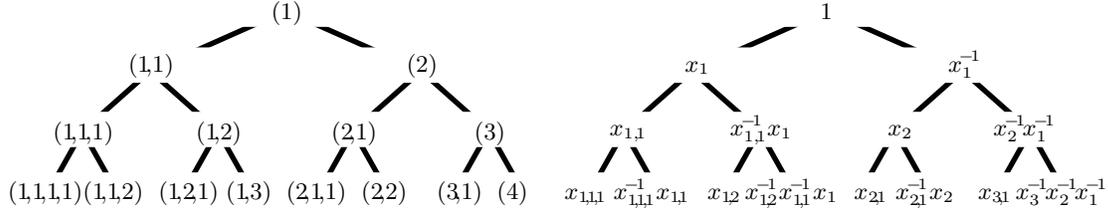}}
\put(0,0){\smaller$(1\!,\!1\!,\!1\!,\!1)$}
\put(10,0){\smaller$(1\!,\!1\!,\!2)$}
\put(20,0){\smaller$(1\!,\!2\!,\!1)$}
\put(29,0){\smaller$(1\!,\!3)$}
\put(37,0){\smaller$(2\!,\!1\!,\!1)$}
\put(47,0){\smaller$(2\!,\!2)$}
\put(57,0){\smaller$(3\!,\!1)$}
\put(65,0){\smaller$(4)$}
\put(6,8){\smaller$(1\!,\!1\!,\!1)$}
\put(25,8){\smaller$(1\!,\!2)$}
\put(43,8){\smaller$(2\!,\!1)$}
\put(62,8){\smaller$(3)$}
\put(16,17){\smaller$(1\!,\!1)$}
\put(53,17){\smaller$(2)$}
\put(35,24){\smaller$(1)$}
\put(108,24){\smaller$1$}
\put(90,17){\smaller$x_1$}
\put(125,17){\smaller$x_1\invv$}
\put(80,8.5){\smaller$x_{1,\!1}$}
\put(96,8.5){\smaller$x_{1,\!1}\invv x_1$}
\put(117,8.5){\smaller$x_2$}
\put(131,8.5){\smaller$x_2\invv\! x_1\invv$}
\put(74,0){\smaller$x_{1\!,\!1\!,\!1}$}
\put(80.5,0){\smaller$x_{1\!,\!1\!,\!1}\invv x_{1\!,\!1}$}
\put(93,0){\smaller$x_{1\!,\!2}$}
\put(98,0){\smaller$x_{1\!,\!2}\invv x_{1\!,\!1}\invv x_1$}
\put(112.5,0){\smaller$x_{2\!,\!1}$}
\put(118,0){\smaller$x_{2\!,\!1}\invv x_2$}
\put(129,0){\smaller$x_{3\!,\!1}$}
\put(134,0){\smaller$x_3\invv\!x_2\invv\!x_1\invv$}
\end{picture}
\caption{\smaller Addresses for the nodes in
trees, and the associated natural
$\FGs$-colours; for each~$s$, the variable~$x_s$
is the natural colour of the node with
address~$(s, 1)$; we recall that $x_s$ is $1$ for $s$ the
empty sequence, whence the colours on the right branch}
\label{F:Addresses}
\end{figure}

\begin{defi} \label{I:Address}
We use finite sequences of positive integers as addresses
for the nodes in binary trees, as described
in Figure~\ref{F:Addresses}. Moreover, we define
for each node its {\it natural
$\FGs$-colour} to be $x_{s, k-1}\inv  x_{k-2}\inv \pp x_{s,
1}\inv x_s$ for the node with address~$(s, k)$.
\end{defi}

In the sequel, it will be convenient to consider
trees in which not only the leaves, but also the
inner nodes are given $\FGs$-colours. 

\begin{defi} \label{I:Coherent}
An $\FGs$-coloured tree will be called {\it
coherent} if the colour at each inner node is the
product of the colours of the left and right sons
of the node (in this order).
\end{defi}

By construction, when we give to each node in a tree~$\tt$
its natural $\FGs$-colour, we obtain a coherent
$\FGs$-coloured tree that will be called the {\it
natural $\FGs$-colouring} of~$\tt$. 

We now introduce a partial action of words on
$\FGs$-coloured trees extending the action on
uncoloured trees. As in the case of
$\Bs$-coloured trees, the point is to specify how
colours behave.

\begin{defi} \label{I:ActionBis}
For~$\ct$ a coherent $\FGs$-coloured tree with
$\Dec(\ct) = (\ct_1, \pp, \ct_n)$ and $n > i$, the
trees~$\ct \act \ss i$ and $\ct \act \aa i$ are
determined by:
\begin{gather}
\Dec(\ct \act \ss i) =
(\ct_1, \pp, \ct_{i-1}, \ctt, \ct_i, \ct_{i+2},
\pp, \ct_n), \\
\Dec(\ct \act \aa i) =
(\ct_1, \pp, \ct_{i-1}, \ct_i \op \ct_{i+1},
\ct_{i+2}, \pp, \ct_n),
\end{gather}
where $\ctt$ is the tree obtained from~$\ct_{i+1}$
by replacing each colour~$y$ with $xyx\inv$, where
$x$ is the colour of the root in~$\ct_i$. Then, for
$\pw$ a word, $\ct \act \pw$ is
defined so that $\ct \act \pw\inv = \ctt$ is
equivalent to $\ctt \act \pw = \ct$, and $\ct
\act (\pwi \pwii) = (\ct \act \pwi) \act
\pwii$ holds.
\end{defi}

It is easy to check that the previous action preserves
coherence. Then we have the following effective method for
determining the automorphism of~$\FGs$
associated with a word~$\pw$. 

\begin{prop} \label{P:ComputingAuto} \label{I:Hat}
For~$\pw$ a parenthesized braid word, put $\Hat\pw
= \psi(\cl\pw)$\ \footnote{where we recall
$\cl\pw$ denotes the element of~$\Bs$ represented
by~$\pw$}. Then $\auto(\pw)$ can be determined as
follows:

(i) Choose a tree~$\tt$ that is large enough to
ensure that $\tt \act \ww$ exists;

(ii) Compute $\ct \act \ww$, where $\ct$ is
the natural $\FGs$-colouring of~$\tt$;

(iii) Then $\auto(\pw)$ maps the natural colour
of every node in~$\ct \act \pw$ to its actual colour
in~$\ct \act \pw$.
\end{prop}

\begin{figure}[htb]
\begin{picture}(145,25)(0, 0)
\put(0,0){\includegraphics{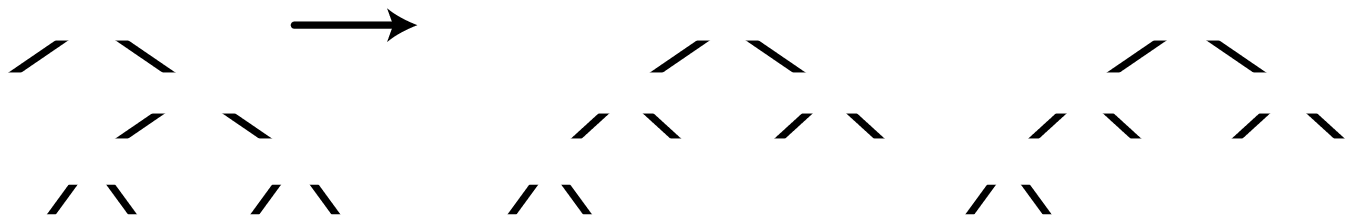}}
\put(4,22){$\ct:$}
\put(10.5,22){\smaller$1$}
\put(53,22){$\ct \act \aa2\ss1:$}
\put(75,22){\smaller$1$}
\put(110,22){$\ctt:$}
\put(122,22){\smaller$1$}
\put(1,15){\smaller$x_1$}
\put(19,15){\smaller$x_1\invv$}
\put(60,15){\smaller$x_1 x_2 x_3 x_1\invv$}
\put(78,15){\smaller$x_1 x_3\invv x_2\invv
x_1\invv$}
\put(111,15){\smaller$x_1$}
\put(130,15){\smaller$x_1\invv$}
\put(10,8){\smaller$x_2$}
\put(28,8){\smaller$x_2\invv\! x_1\invv$}
\put(53.5,8){\smaller$x_1 \! x_2 x_1\invv$}
\put(66,8){\smaller$x_1 \! x_3 x_1\invv$}
\put(79,8){\smaller$x_1$}
\put(85,8){\smaller$x_3\invv\!x_2\invv\!x_1\invv$}
\put(102,8){\smaller$x_{1,1}$}
\put(114,8){\smaller$x_{1,1}\invv  x_1$}
\put(125,8){\smaller$x_2$}
\put(135,8){\smaller$x_2\invv\!x_1\invv$}
\put(4,0){\smaller$x_{2,1}$}
\put(12,0){\smaller$x_{2,1}\invv x_2$}
\put(24,0){\smaller$x_3$}
\put(30,0){\smaller$x_3\invv\!x_2\invv\!x_1\invv$}
\put(46,0){\smaller$x_1 \! x_{2,\!1} \!
x_1\invv$}
\put(58,0){\smaller$x_1\!x_{2,\!1}\invv\!x_2x_1\invv$}
\put(95,0){\smaller$x_{1,1,1}$}
\put(104,0){\smaller$x_{1,1,1}\invv x_{1,1}$}
\put(34,20){$\aa2\ss1$}
\put(87,22){compare with}
\end{picture}
\caption{\smaller Computing the
automorphism of~$\FGs$ associated with~$\aa2\ss1$:
we let $\aa2 \ss1$ act on a tree~$\ct$ with
natural $\FGs$-colours, and compare the colours
in~$\ct \act \aa2\ss1$ with the natural ones:
the node with natural colour~$x$ has colour
$\psi(\aa2\ss1)(x)$ in~$\ct \act \aa2\ss1$. For
instance, $x_1$ is mapped to~$x_1
x_2 x_3 x_1\inv$ and that $x_{1,1}\inv x_1$ is
mapped to~$x_1 x_3 x_1\inv$.}
\label{F:ComputingAuto}
\end{figure}

\begin{proof} 
(See Figure~\ref{F:ComputingAuto} for an
example). For $\ct$ an $\FGs$-coloured tree and
$\theta$ a mapping of~$\FGs$ into itself, we
denote by~$\ct^\theta$ the tree obtained
from~$\ct$ by replacing each colour~$x$
with~$\theta(x)$. What we want to prove is the
equality $\ct \act \pw = {\ctt}^{\auto(\pw)}$
where $\ctt$ is the natural $\FGs$-colouring
of~$(\ct \act \pw)^\dag$.

A direct inspection shows that the result is true
when $\pw$ is a single letter~$\ssss i$
or~$\aaaa i$. So the point is to show that the
result is true for~$\pw = \pwi \pwii$ when
it is for~$\pwi$ and~$\pwii$. Assume that $\tt
\act \pw$ exists. Denote by~$\cti$ the natural
$\FGs$-colouring of~$\tt\act\pwi$. By induction
hypothesis, we have
$\ct \act \pwi = \cti^{\auto(\pwi)}$, hence
$\ct \act \pw = \cti^{\auto(\pwi)} \act \pwii$.
By induction hypothesis again, we have 
$\cti \act \pwii = \ctt{}^{\auto(\pwii)}$, which
means that each node with colour~$x$ in~$\ctt$,
has colour~ $\auto(\pwii)(x)$ in~$\cti \act
\pwii$. By construction, this colour is an
expression $E(x_{s_1}, \pp, x_{s_p})$ involving
some variables~$x_{s_1}, \pp, x_{s_p}$ with
products and inverses. When we substitute~$\cti$
with~$\cti^{\auto(\pwi)}$ and let~$\pwii$ act,
the result is the corresponding expression
$E(\auto(\pwi)(x_{s_1}), \pp,
\auto(\pwi)(x_{s_p}))$, which is also
$\auto(\pwi)(E(x_{s_1}, \pp, x_{s_p}))$ as
$\auto(\pwi)$ is a group automorphism. This means
that $\cti^{\auto(\pwi)} \act \pwii$, which is
$\ct \act \pw$, is $\ctt{}^{\auto(\pwi)
\circ \auto(\pwii)}$, \ie, $\ctt{}^{\auto(\pw)}$, as
expected.
\end{proof}

\begin{rema}
The (partial) actions of~$\Bs$ on $\FGs$-
and $\Bs$-coloured trees extends to all
$\SetOfCol$-coloured trees where
$\SetOfCol$ is a left cancellative ALD-system.
\end{rema}

\subsection{The injectivity result}

Artin's representation of~$\Bi$ is an
embedding~\cite{Bir}. We extend the result
to~$\Bs$, so obtaining a realization of~$\Bs$ as a
group of automorphisms of a free group. 

\begin{prop} \label{P:Embe}
The representation~$\psi$ of~$\Bs$ in~$\Aut(\FGs)$ is an
embedding.
\end{prop}

\begin{coro}
The morphism~$\phi$ of~$\Bs$ into~$\MCG(\SK)$ is
injective.
\end{coro}

The method for proving Proposition~\ref{P:Embe}
relies on the possibility of considering words~$\pw$
of a specific form, in connection with the linear
ordering of~$\Bs$ constructed in
Section~\ref{S:Order}. In the case of braids, the
method was first used by D.~Larue in~\cite{Lra}, and
it gives a powerful method for proving the possible
injectivity of a representation~\cite{Shp, CrP}. 

\begin{defi} \label{I:Red}
For~$\uu$ a word in the letters~$x_s\pmo$,
we denote $\red(\uu)$ for the freely reduced word
obtained from~$\uu$ by removing all pairs~$x x\inv$
and $x\inv x$. 
\end{defi}

Thus $\FGs$ identifies with the set of all freely reduced
words. We recall that $\Hat\pw$ denotes the
automorphism~$\psi(\cl\pw)$ of~$\FGs$ associated
with~$\pw$. 

We begin with two auxiliary results. The first one
is similar to Proposition~5.1.6 of~\cite{Dgr} for
braids. The only change is that variables~$x_s$
with $s$ of length more than~$1$ may occur, but
this does not change the argument.

\begin{lemm} \label{L:Embe1}
The image of a word ending
with~$x_i\inv$ under~$\hs_i$ or~$\hs_j\pmo$ with
$j > i$ ends with~$x_i\inv$.
\end{lemm}

\begin{proof}
Assume that $\uu$ ends with~$x_i\inv$, say $\uu =
\uu' x_i\inv$. Then we have
\begin{equation} \label{E:Larue}
\hs_i(\uu) = \red(\hs_i(\uu') x_i x_{i+1}\inv
x_i\inv).
\end{equation}
In order to prove that the word above ends
with~$x_i\inv$, it is sufficient to check that
the final~$x_i\inv$ cannot be cancelled during the
reduction by some~$x_i$ coming
from~$\hs_i(\uu')$. By~\eqref{E:Artins}, an~$x_i$
in~$\hs_i(\uu')$ must come from some~$x_i$,
$x_i\inv$, or~$x_{i+1}$ in~$\uu'$. We consider the 
three cases, displaying the supposed
involved letter in~$\uu'$. For $\uu' = \uu''
x_i \uu'''$,
\eqref{E:Larue} becomes
\begin{equation*}
\hs_i(\uu) =
\red(\hs_i(\uu'') x_i x_{i+1} x_i\inv
\hs_1(\uu''') x_i x_{i+1}\inv x_i\inv).
\end{equation*}
The assumption that the first~$x_i$ cancels
the final~$x_i\inv$ implies $\hs_i(\uu''') =
\e$, hence $\uu''' = \e$, contradicting the
hypothesis that $\uu'' x_i
\uu''' x_i\inv$ is  reduced.  For $\uu' = \uu''
x_i\inv \uu'''$,
\eqref{E:Larue} is
\begin{equation*}
\hs_i(\uu) =
\red(\hs_i(\uu'') x_i x_{i+1}\inv x_i\inv
\hs_1(\uu''') x_i x_{i+1}\inv x_i\inv).
\end{equation*}
The assumption that the first~$x_i$ cancels
the final~$x_i\inv$ implies now that
$x_{i+1}\inv x_i\inv \hs_i(\uu''') x_i
x_{i+1}\inv$ reduces to~$\e$,
hence $\hs_i(\uu''') = x_i x_{i+1}^2
x_i\inv$, and, therefore, $\uu''' = x_i^2$,
again contradicting the hypothesis that $\uu''
x_i\inv \uu'''$ is reduced.
Finally, for $\uu' = \uu'' x_{i+1} \uu'''$, 
\eqref{E:Larue} says
\begin{equation*}
\hs_i(\uu) =
\red(\hs_i(\uu'') x_i
\hs_1(\uu''') x_i x_{i+1}\inv x_i\inv).
\end{equation*}
The assumption that the first~$x_i$ cancels
the final~$x_i\inv$ implies that
$\hs_i(\uu''') x_i x_{i+1}\inv$ reduces to~$\e$,
hence
$\hs_i(\uu''') = x_{i+1} x_i\inv$, and,
then, $\uu''' = x_{i+1}\inv
x_i$, contradicting the hypothesis
that $\uu'' x_{i+1} \uu'''$ is reduced.
We similarly consider the action
of~$\hs_j^e$ with $j > i$ and $e = \pm 1$. We
find
\begin{equation} \label{E:Larue2}
\hs_j(\uu) = \red(\hs_j^e(\uu') x_i\inv),
\end{equation}
and aim at proving that the
final~$x_i\inv$ cannot vanish in reduction. Now
it could do it only with some~$x_i$
in~$\hs_j^e(\uu')$, itself coming from some~$x_i$
in~$\uu'$. For a contradiction, we display the
latter as $\uu' = \uu'' x_i \uu'''$. Then
\eqref{E:Larue2} becomes
$\hs_j(\uu) = \red(\hs_j^e(\uu'') x_i
\hs_j^e(\uu''') x_i\inv)$. As above, we must have
$\hs_j^e(\uu''') = \e$, hence $\uu''' = \e$,
contradicting the hypothesis that $\uu'' x_i \uu'''
x_i\inv$ is reduced.
\end{proof}

The second preliminary result is specific to our
current situation.

\begin{defi} \label{I:SpecialBis}
A word in the letters~$x_s\pmo$ is said to
be {\it special} if it is freely reduced and it
admits a suffix of the form $x_s\inv
x_{s,j_1,s_1} \pp x_{s, j_r, s_r}$ with $r \ge
0$, where $s, s_1, \pp, s_r$ are sequences, and
$j_1, \pp, j_r$ are positive integers.
\end{defi}

Thus $x_1\inv$ and $x_1 x_2\inv x_{2,1}$ are
special words.

\begin{lemm} \label{L:Embe2}
For each~$i$, the image of a special word
under~$\ha_i\inv$ is a special word.
\end{lemm}

\begin{proof}
Let $\uu = \uu' x_{j,s}\inv x_{j,s,j_1,s_1} \pp
x_{j,s, j_r, s_r}$ be a special word. We consider
the image of~$\uu$ under~$\ha_i\inv$,
according to the mutual positions of~$i$
and~$j$. Assume first $j < i$. Then we have
$\ha_i\inv(x_{j,s}) = x_{j,s}$, and, similarly,
$\ha_i\inv(x_{j,s, j_k, s_k}) = x_{j,s, j_k, s_k}$
for each~$k$, hence
\begin{equation} \label{E:Emb1}
\ha_i\inv(\uu) =  \red(\ha_i\inv(\uu')
x_{j,s}\inv x_{j,s, j_1, s_1} \pp x_{j,s,
j_r, s_r}).
\end{equation}
In order to conclude that this word is special, it
suffices to prove that the displayed letter
$x_{j,s}\inv$ cannot vanish during reduction. Now
assume it does. The letter~$x_{j,s}\inv$ is
cancelled by some letter~$x_{j,s}$ coming
from~$\ha_i\inv(\uu')$. The explicit formulas
for~$\ha_i\inv$ are
\begin{equation}
\label{E:Acaxi}
\ha_i\inv: 
\begin{cases}
\quad x_{j,s} \mapsto x_{j,s}
\text{\ for $j < i$},
\quad x_{j+1,s} \mapsto x_{j,s}
\text{\ for $j > i$},\\
\quad x_{i,s} \mapsto x_{i,1,s},
\qquad x_{i+1} \mapsto x_{i,1}\inv x_i,
\qquad x_{i+1, j, s} \mapsto x_{i,j+1, s}
\text{\ for $j \ge 2$}.
\end{cases}
\end{equation}
So a letter~$x_{j,s}$ in~$\ha_i\inv(\uu')$ must come
from a letter~$x_{j,s}$ of~$\uu'$. Let us display the
considered letter and write $\uu' = \uu'' x_{j,s}
\uu'''$. Then \eqref{E:Emb1} becomes
$$\ha_i\inv(\uu') = \red(\ha_i\inv(\uu'') x_{j,s}
\ha_i\inv(\uu''')).$$
The assumption that the
final~$x_{j,s}\inv$ in $\ha_i\inv(\uu') x_{j,s}\inv$ is
cancelled by the displayed~$x_{j,s}$ implies
$\ha_i\inv(\uu''') = \e$, hence $\uu''' = \e$ as
$\ha_i$ is an automorphism. This means that $\uu'$
finishes with~$x_{j,s}$, contradicting the
hypothesis that $\uu' x_{j,s}\inv$ is reduced.

The argument is similar for~$x_j$ with $j >
i+1$, and, more generally, it works for
all~$x_{t}$'s except~$x_i$ and~$x_{i+1}$. Indeed,
in these cases,
$\ha_i\inv$ maps~$x_{j,s}$ to a (possibly different)
letter~$x_{j',s'}$ so that a letter~$x_{j',s'}$
in~$\ha_i\inv(\vv)$ must come from a~$x_{j,s}$
in~$\vv$. Then, the previous argument shows that the
letter~$x_{j,s}\inv$ witnessing for
specialness becomes a
letter~$x_{j',s'}\inv$ that cannot be
cancelled. On the other hand,
\eqref{E:Acaxi} shows that, in all considered
cases, the final letters~$x_{j,s, j_k, s_k}$ become
letters~$x_{j',s', j_k, s_k}$, so the
word~$\ha_i\inv(\uu)$ is special.

There remain the cases of~$x_i$ and~$x_{i+1}$. To
simplify reading, we assume $i = 1$. Let us first
consider~$x_1$,  \ie, $\uu =
\uu' x_1\inv x_{1,j_1, s_1} \pp x_{1,j_r, s_r}$,
which gives 
\begin{equation} \label{E:Embe}
\ha_1\inv(\uu) =
\red(\ha_1\inv(\uu') x_{1,1}\inv x_{1,1,j_1,
s_1} \pp x_{1,1,j_r,s_r}).
\end{equation}
If the displayed $x_{1,1}\inv$ does not vanish
during reduction, the above word is
special. We shall see now that $x_{1,1}\inv$ may
vanish, but one nevertheless obtains a special
word. Indeed, \eqref{E:Acaxi} shows that
an~$x_{1,1}$ in $\ha_1\inv(\uu')$ comes either
from an~$x_1$ or from an~$x_2\inv$ in~$\uu'$. By
the same argument as above, $x_1$ is excluded. So
assume $\uu' = \uu'' x_2\inv
\uu'''$. Then
\eqref{E:Embe} becomes
$$\ha_1\inv(\uu) =
\red(\ha_1\inv(\uu'') x_1\inv
x_{1,1} \ha_1\inv(\uu''' ) x_{1,1}\inv
x_{1,1,j_1, s_1} \pp x_{1,1,j_r,s_r},$$
and the assumption is
$\ha_1\inv(\uu''' ) = \e$. As above, we deduce
$\uu''' = \e$, hence $\uu' = \uu''
x_2\inv$---which is not forbidden. In this case, we
find
\begin{equation} \label{E:Emb4}
\ha_1\inv(\uu) = \red(\ha_1\inv(\uu'')
x_1\inv x_{1,1,s_1} \pp x_{1,1,s_r}).
\end{equation}
To show that this word is special, it is
sufficient to prove that the~$x_1\inv$
cannot disappear. Now the only way
$x_1\inv$ could vanish is with some~$x_1$ in
$\ha_1\inv(\uu'')$, necessarily coming from
some~$x_2$ in~$\uu''$. Write $\uu'' = \uu'''
x_2 \uu''''$. As above, we obtain
$\ha_1\inv(\uu'''') = \e$, hence $\uu'''' = \e$,
implying that $\uu''$ finishes with~$x_2$, and
contradicting the hypothesis that $\uu'' x_2\inv$ is
reduced. So the study for~$x_1$ is complete.

Finally, let us consider the case of~$x_2$. The
problem here is that $\ha_1\inv$ maps $x_2$
to~$x_{1,1}\inv x_1$, which is not a single letter.
So assume $\uu = \uu' x_2\inv x_{2,j_1,s_1} \pp
x_{2,j_r,s_r}$. We obtain 
\begin{equation} \label{E:Emb5}
\ha_1\inv(\uu) =
\red(\ha_1\inv(\uu') x_1\inv x_{1,1}
x_{1,j_1+1,s_1} \pp x_{1,j_r+1,s_r}).
\end{equation}
In order to show that this word is special, it
suffices to prove that the letter~$x_1\inv$ cannot
vanish. Now  a letter~$x_1$ in~$\ha_1\inv(\uu')$
must come from a letter~$x_2$ in~$\uu'$, and we
argue as above. 
\end{proof}

We can now prove the injectivity of the
homomorphism~$\psi$ of~$\Bs$ into~$\Aut(\FGs)$.

\begin{proof} [Proof of
Proposition~\ref{P:Embe}]
Our aim is to show that, if $\pw$ is a word
that represents a non-trivial element of~$\Bs$, then
the automorphism~$\Hat\pw$ (\ie, $\psi(\cl\pw)$)
is not the identity mapping, \ie, there exists at
least one letter~$x_s$ such that $\Hat\pw(x_s)$ is
not~$x_s$. By Proposition~\ref{P:Partition}, at the
expense of replacing~$\pw$ by an equivalent word and
possibly exchanging~$\pw$ and~$\pw\inv$, we may
assume that $\pw$ is either $\ss i$-positive or
is a non-trivial $\wa$-word.

{\bf Case 1:} $\pw$ is $\ss i$-positive. By
definition, we can write $\pw = \pwi\inv \pwii
\pwiii$, where $\pwi$ and $\pwiii$ are positive
$\wa$-words, and $\pwii$ is a $\ss i$-positive
$\ws$-word. First, because $\pwiii$ contains positive
letters~$\aa k$ only, there exists a vine~$\tt$
such that $\tt \act \pwiii$ is defined and we may
assume in addition that the right height of~$\tt$ is
at least~$i+1$. Let~$\ct$ be the natural
$\FGs$-colouring of~$\tt$. By construction,
$x_i$ is a colour in~$\ct$, hence in~$\ct \act
\pwiii$, and Proposition~\ref{P:ComputingAuto}
implies that there must exist~$x$ in~$\FGs$ such
that $\Hat\pwiii$ maps~$x$ to~$x_i$. All colours in a
natural $\FGs$-colouring are not single variables,
but this is always the case for nodes with
addresses ending with~$1$. So, in any case, the
left son of the node where~$x_i$ occurs has
colour~$x_{i, 1}$ in~$\ct
\act \pwiii$, and colour~$x_s$ for some~$s$ in the
natural colouring of~$\ct \act \pwiii$. In other
words, there exists~$s$ satisfying $\Hat\pwiii(x_s) =
x_{i,1}$.

We now consider $\Hat\pwii(\Hat\pwiii(x_s))$,
\ie, $\Hat\pwii(x_{i,1})$. Write $\pwii = \pww_0
\ss i \pww_1 \ss i \pp \ss i \pww_r$, where
$\pww_k$ contains no~$\ss j\pmo$ with $j \le
i$. Then $\Hat{\pww_r}$ fixes~$x_{i,1}$, while
$\ss i $ maps it to~$x_i x_{i+1,1} x_i\inv$, a
reduced word ending with~$x_i\inv$. Applying
Lemma~\ref{L:Embe1} repeatedly, we deduce that the
final~$x_i\inv$ cannot disappear, and, so,
$\Hat\pwii(\Hat\pwiii(x_s))$ is a reduced word
ending with~$x_i\inv$.

Consider now the action of~$\Hat\pwi\inv$ on the
latter word. Every reduced word ending with~$x_i\inv$
is a special word, hence, by Lemma~\ref{L:Embe2}, its
image under~$\Hat\pwi\inv$ is a special word.
Hence $\Hat\pw(x_s)$ is a special word. As $x_s$ is
not a special word,
$\Hat\pw$ cannot be the identity mapping.

{\bf Case 2:} $\pw$ is a non-trivial
$\wa$-word. Let $\tt, \ttt$ be trees satisfying
$\ttt = \tt \act
\pw$. The hypothesis that $\pw$ is non-trivial
implies $\ttt \not= \tt$. Then there must exist an
address~$s$ such that $(s, 1)$ is an address
in~$\ttt$ and not in~$\tt$. Then $x_s$ occurs in the
natural $\FGs$-colouring~$\ctt$ of~$\ttt$, and not in
the natural $\FGs$-colouring~$\ct$ of~$\tt$.
Proposition~\ref{P:ComputingAuto} implies
that $\Hat\pw(x_s)$ is a combination of colours
occurring in~$\ctt$, so it cannot be~$x_s$, and
$\Hat\pw$ is not the identity mapping.
\end{proof}

An application of Proposition~\ref{P:Embe} is an
alternative proof of the fact that the
relations~$\RRs$ make a presentation of the
group~$\Bs$. Indeed, ignoring the injectivity of
$\pi: \BBs \to \Bs$, we can construct a
morphism~$\widetilde\psi$ of~$\BBs$
to~$\Aut(\FGs)$ using the explicit formulas of
Proposition~\ref{P:AutoFree}. Then
Proposition~\ref{P:ComputingAuto} shows that,
for each word~$\pw$, the
automorphism~$\widetilde\psi(\cl\pw)$ can be
recovered from the action of~$\pw$ on
$\FGs$-coloured trees. Now the latter can in turn be
deduced from the diagram~$\DD(\pw)$ using
$\FGs$-colourings, hence from the isotopy class
of~$\DD(\pw)$ as isotopy preserves colours. So
$\widetilde\psi(\cl\pw)$ depends on the image
of~$\pw$ in~$\Bs$ only, \ie, $\widetilde\psi$
factors through~$\Bs$:
$$\begin{picture}(67,17)(0, 0)
\put(20,3){\includegraphics{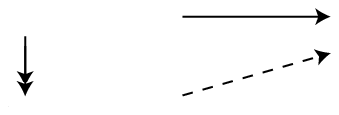}}
\put(10,12){$\BBs = \Gr(\XXsa; \RRs)$}
\put(58,12){$\Aut(\FGs)$}
\put(18,7){$\pi$}
\put(43,8){$\psi$}
\put(43,14){$\widetilde\psi$}
\put(-5,0){$\Bs=\{\text{parenthesized
diagrams}\}/\text{isotopy}$}
\end{picture}$$
What Proposition~\ref{P:Embe} shows is that
$\widetilde\psi$ is injective, which implies
that both~$\pi$ and~$\psi$ are injective.

\section{Miscellani}

We conclude with a few additional remarks
about~$\Bs$.

\subsection{Pure parenthesized braids}

Each braid induces a permutation of positive
integers, which leads to a surjective
homomorphism of~$\Bi$ onto the group~$\Symi$ of
eventually trivial permutations. The
group~$\Symi$ is the quotient of~$\Bi$ under the
relations~$\ss i^2 = 1$, and the kernel is the
pure braid group~$\PBi$. The situation is similar
with~$\Bs$. The quotient of~$\Bs$ obtained by
adding the relations~$\ss i^2 = 1$ is the
subgroup~$\Syms$ of Thompson's group~$V$ made of
the elements that, in the action of~$V$ on the
Cantor set~$\Can$, preserve the right endpoint;
see~\cite{Dhb}, and~\cite{Bri1, Bri2} where this group
is called~$\widehat{V}$. Then the kernel of the
projection $\Bs \to \Syms$ is a non-trivial
normal subgroup~$\PBs$ of~$\Bs$, whose elements
can be called {\it pure} parenthesized braids.

\begin{prop}
We have $\PBs = (\Fp)\inv \opp \PBi \opp \Fp$.
\end{prop}

One inclusion is trivial, and the other follows from
the equality $\Bs = (\Fp)\inv \opp \Bi \opp \Fp$.

\subsection{Alternative presentations}

Alternative presentations of~$\Bs$ have been
considered. On the one hand, exactly as Thompson's
group~$F$ is generated by the two elements here
denoted~$\aa1$ and~$\aa2$, the group~$\Bs$ is
generated by~$\ss1, \ss2, \aa1, \aa2$, and it is
a finitely presented group \cite{Bri2}. 

On the other hand, large presentations may also
of interest. The presentation~$(\XXsa,
\RRs)$ gives different roles to the left
and right sides. This in particular implies that
$\Bs$ is a group of left fractions of~$\Bsp$
only, and that right common multiples need not
exist in~$\Bsp$. As shown in~\cite{Dhb}, $\Bs$, as
well as Thompson's groups~$F$ and~$V$, can be given a
balanced presentation. The
principle is to consider new generators similar
to~$\ss i$ and~$\aa i$ but acting at any possible
address in a tree, and not only at addresses on
the rightmost branch. In the current framework,
it is natural to denote by~$\ss s$ and~$\aa s$
such generators, with
$s$ a finite sequence of positive integers. For
instance, $\ss{1,1}$ corresponds to applying~$\ss1$
at the address~$(1,1)$ (in the sense of
Figure~\ref{F:Addresses}) instead of
at~$(1)$, which amounts to defining
$\ss{1,1} = \aaa1 \aaa2 \ss1 \aa2 \aa1$. We obtain
in this way an extended double family of
generators~$\ss s$, $\aa s$, and, using the
techniques of~\cite{Dhb}, one can show:

\begin{prop}
In terms of the generators~$\ss s$ and~$\aa s$, a
presentation of~$\Bs$ is:
\begin{gather}
x_{s, i, s'} \op y_{s, j, s''} 
= y_{s, j, s''} \op x_{s, i, s'}
\qquad \text{for $j \not= i$},\\
\ss{s, i} \op x_{s, j, s'} 
= x_{s, j, s'} \op \ss{s, i}
\qquad
\aa{s, i} \op x_{s, j, s'} 
= x_{s, j-1, s'} \op \aa{s, i}
\qquad \text{for $j \ge i+2$},\\
x_{s, i, j, s'} \op \ss{s, i}
= \ss{s, i} \op x_{s, i+1, j, s'}, 
\qquad
x_{s, i+1, j, s'} \op \ss{s, i}
= \ss{s, i} \op x_{s, i, j, s'},\\
x_{s, i, 1, s'} \op \aa{s, i}
= \aa{s, i} \op x_{s, i, 1, 1, s'}, 
\qquad
x_{s, i+1, j, s'} \op \aa{s, i}
= \aa{s, i} \op x_{s, i, j+1, s'},\\
\ss{s, i} \op \ss{s, i+1} \op \ss{s, i}
= \ss{s, i+1} \op \ss{s, i} \op \ss{s, i+1},
\ 
\ss{s, i+1} \op \ss{s, i} \op \aa{s, i+1}
= \aa{s, i} \op \ss{s, i},
\ \ss{s, i} \op \ss{s, i+1} \aa{s, i},
= \aa{s, i+1} \op \ss{s, i},\\
\label{E:Pentagon}
\ss{s, i} \op \aa{s, i+1} \op \aa{s, i}
= \aa{s, i+1} \op \aa{s, i} \op \aa{s, i, 1},
\qquad
\aa{s, i} \op \aa{s, i}
= \aa{s, i+1} \op \aa{s, i} \op \aa{s, i, 1},
\end{gather}
with $i, j$ positive integers, $s, s', s''$ sequences
of positive integers, and  $x, y$ denoting any
of~$\ss¥$ or~$\aa¥$.
\end{prop}

Despite its apparent complexity, the above
presentation is simple: in addition to the
relations of~$\RRs$, it only contains more or
less trivial commutation relations, plus the
last relation in~\eqref{E:Pentagon}, which is
MacLane's pentagon relation~\cite{Mac}.
The advantage of this presentation is that it
restores the symmetry between left and right---this
becomes more evident when sequences of~$0$'s
and~$1$'s are used as addresses \cite{Dhb}. In
particular, the presentation leads to a new
monoid, larger than~$\Bsp$, in which both left
and right lcm's exist, and $\Bs$ is both a group
of left and right fractions of this monoid.

\subsection{Artin's representation of the group~$BV$}

In~\cite{Bri1, Bri2}, M.\,Brin introduces two
groups denoted~$BV$ and~$\widehat{BV}$, for which he
establishes presentations. The presentation
of~$\widehat{BV}$ shows that this group is
isomorphic to~$\Bs$. The group~$BV$, which is an
extension of Thompson's group~$V$,
includes~$\widehat{BV}$, hence~$\Bs$, as a subgroup,
but, at the same time, it identifies with the
subgroup~$\Bs^{(1)}$ of~$\Bs$ consisting of the
parenthesized braids in which only the strands starting at a
positions~$(1, s)$---\ie, $1$ or infinitely
close---may be braided. For instance, $\aaa1 \ss1
\aa1$ is a typical element of~$\Bs^{(1)}$. By using
the Artin representation of~$\Bs^{(1)}$, we obtain a
representation of the group~$BV$ into~$\Aut(\FGs)$.
>From the point of view of an action on trees,
$BV$ can be obtained from~$\Bs$ by adding new
generators~$c_i$, $i \ge 1$, whose effect is to
switch the subtrees~$\tt_i$ and $\tt_{i+1} \pp
\tt_n\et$ of the right decomposition. 

\begin{prop}
Defining
$$\psi(c_i) : x_{j, s} \mapsto x_{j, s} \text{~for $j
< i$},\quad  x_i \mapsto x_i\inv, \quad
x_{i, j, s} \mapsto x_i x_{i+j, s} x_i\inv, \quad
x_{i+j, s} \mapsto x_{i, j, s}$$
extends the embedding~$\psi$ of
Proposition~\ref{P:AutoFree} to the group~$BV$. 
\end{prop}

\subsection{Further questions}

Owing to the many results about~$\Bi$ and~$F$, in
particular in terms of (co)-homology and geometry
of the Cayley graph, investigating~$\Bs$ in these
directions seems a promising project. 

\section{Appendix: The cube condition for
the presentation $(\XXsa, \RRs)$}

The algebraic results of
Section~\ref{S:AlgebraicProp} rely on the fact that
the presentation $(\XXsa, \RRs)$ satisfies the
so-called left and right cube conditions. Verifying
these combinatorial properties requires that we consider all
possible triples of letters. There are infinitely many of
them, but only finitely many different patterns may
appear, and the needed verifications are finite
in number. Here we give some details.

\subsection*{The left cube condition}

The left cube condition for a triple of
letters~$(\gx, \gy, \gz)$ claims that, whenever the
word~$\gx \gy\inv \gy \gz\inv$ is left reversible
to some word~$\vv\inv \uu$ with $\uu, \vv$ containing
no negative letter, then $\vv \gx  \gz\inv \uu\inv$
is left reversible to the empty word~$\e$.

In the presentation $(\XXsa, \RRs)$,
there exists exactly one relation
$\uu\gx = \vv \gy$ for each pair of letters~$\gx,
\gy$, hence there exists at most one way to
reverse a word~$\pw$ to a word of the
form~$\vv\inv \uu$ with
$\uu, \vv$ positive. We shall denote by $\uu \cL
\vv$ the unique positive~$\uu'$ such that $\uu
\vv\inv$ is left reversible to~$\vv'\Inv \uu'$ for
some positive~$\vv'$, if such words exist.
If $\pw$ is left reversible to~$\pw'$, then
$\pw\inv$ is left reversible to~$\pw'\Inv$, and
therefore, if $\uu \vv\inv$ is left reversible
to~$\vv'\Inv \uu'$, the latter is $(\vv \cL
\uu)\inv (\uu \cL \vv)$. So, for instance, we
have $\ss1 \cL \ss2 = \ss2
\ss1$ and $\ss2 \cL \ss1 = \ss1 \ss2$, and
\eqref{E:LeftRevSA} rewrites as
\begin{equation} \label{E:AuxRev0}
\ss i \cL \aa j = \dbl_j(\ss i), \qquad 
\aa j \cL \ss i = \aa{\ss i[j]}.
\end{equation}
In the case of two $\aa i$'s, the
formula for~$/$ always takes the form $\aa i \cL \aa
j = \aa{i'}$. The index~$i'$ will be denoted~$i \cL
j$. For instance, one has $1 \cL 2 = 1$ and $2 \cL
1 = 3$. It is then easy to verify the equalities
\begin{equation}
\label{E:AuxRev1}
\dbl_k(\ss i) \cL \dbl_k(\ss j) \equiv
\dbl_{\ss j[k]}(\ss i \cL \ss j), \qquad
\ss k[i] \cL \ss k[j] = \dbl_j(\ss k) [i \cL j],
\end{equation}
where $\equiv$ denotes $\RRs$-equivalence.
Let us
write
\begin{picture}(16,6)(0,5)
\put(3,2){\includegraphics{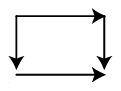}}
\put(0,5){$\vv'$}
\put(14,5){$\vv$}
\put(6.5,5){$\revL$}
\put(7,0){$\uu$}
\put(7,10){$\uu'$}
\end{picture}
\vrule width0pt depth5mm
when $\uu \vv\inv$ is left reversible
to~$\vv'\Inv \uu'$. The left cube
condition for~$(\gx, \gy, \gz)$ means that,
when we fill the diagram\ 
\vrule width0pt depth8mm
\begin{picture}(26,6)(0, 8)
\put(3,2){\includegraphics{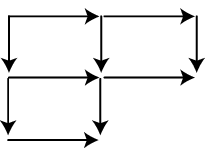}}
\put(-0.5,6){$\vv_1$}
\put(-0.5,12){$\vv_2$}
\put(14.5,5.5){$\gy$}
\put(17,6.5){$\gy$}
\put(24,12){$\gz$}
\put(6.5,5){$\revL$}
\put(6.5,11.5){$\revL$}
\put(16,11.5){$\revL$}
\put(7,0){$\gx$}
\put(7,16.5){$\uu_2$}
\put(16.5,16.5){$\uu_1$}
\end{picture}, then the word $\vv_1 \vv_2 \gx
\gz\inv \uu_1\inv \uu_2\inv$ must be left
reversible to~$\e$, \ie, filling the
corresponding diagram leads to
$\e$~edges on the left and the top side. 

\begin{figure}[htb]
\begin{picture}(110,45)(0,-4)
\put(-1,-1){\includegraphics{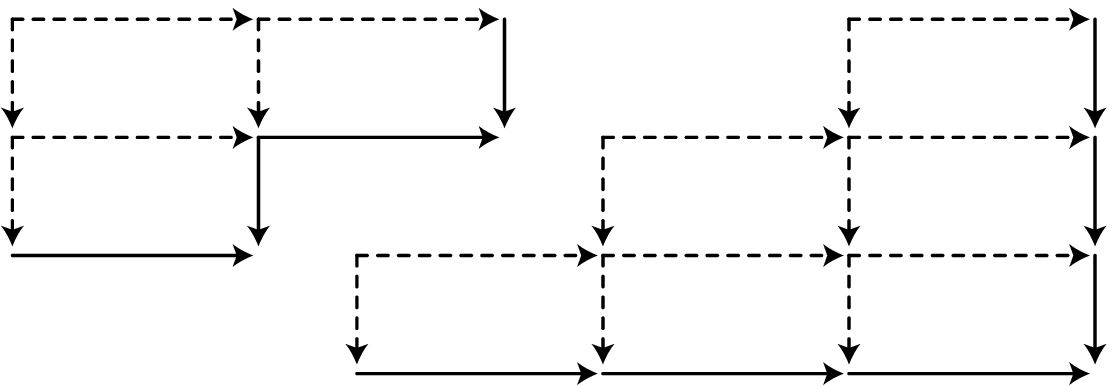}}
\put(-18,30){$\aa{\ss i \cL \ss j[\ss j[k]]}$}
\put(-10,18){$\ss j \cL \ss i$}
\put(12,9){$\ss i$}
\put(8,26){$\ss i \cL \ss j$}
\put(27,18){$\ss j$}
\put(36,21){$\ss j$}
\put(52,30){$\aa k$}
\put(3,38){$\scriptstyle\dbl_{\ss j[k]}(\ss i \cL \ss j)$}
\put(33,38){$\scriptstyle\dbl_k(\ss j)$}
\put(27,30){$\aa{\ss j[k]}$}
\put(38,-3){$\aa{\ss i \cL \ss j[\ss j[k]]}$}
\put(70,-3){$\ss j \cL \ss i$}
\put(97,-3){$\ss i$}
\put(111,30){$\scriptstyle\dbl_{\ss j[k]}(\ss i \cL \ss j)$}
\put(111,18){$\scriptstyle\dbl_k(\ss j)$}
\put(111,6){$\aa k$}
\put(94,14){$\scriptstyle\dbl_k(\ss i)$}
\put(87,26){$\scriptstyle\dbl_k(\ss j) \cL \dbl_k(\ss i)$}
\put(86,18){$\scriptstyle\dbl_k(\ss i) \cL \dbl_k(\ss j)$}
\put(63,14){$\scriptstyle\dbl_{\ss i[k]}(\ss j \cL \ss i)$}
\put(61,6){$\aa{\ss j \cL \ss i(\ss i[k])}$}
\put(86,6){$\aa{\ss i[k]}$}
\put(32,6){$\e$}
\put(57,18){$\e$}
\put(82,30){$\e$}
\put(48,14){$\e$}
\put(73,26){$\e$}
\put(98,38){$\e$}
\put(11,18){$\revL$}
\put(11,30){$\revL$}
\put(36,30){$\revL$}
\put(96,30){$\revL$}
\put(96,6){$\revL$}
\put(46,6){$\revL$}
\put(71,18){$\revL$}
\end{picture}
\caption{\smaller Left cube condition for the triples
$(\ss{},\ss{},\aa{})$: one first reverses $\ss i
\sss j \ss j \aa k\inv$ to $(\ss j
\cL \ss i)\inv (\aa{\ss i
\cL \ss j[\ss j[k]]})\inv (\dbl_{\ss j[k]}(\ss i \cL
\ss j)) (\dbl_k(\ss j))$, then restart from
$(\aa{\ss i \cL \ss j[\ss j[k]]}) \,
(\ss j \cL \ss i) \, (\ss i) \, (\ss k)\inv \,
(\dbl_k(\ss j))\inv \,
(\dbl_{\ss j[k]}(\ss i \cL \ss j))\inv$ and check
that the latter is left reversible to~$\e$; the
values follow from~\eqref{E:AuxRev1} and the fact
that the permutations associated with $(\ss i
\cL \ss j)\ss j$ and $(\ss j \cL \ss i)\ss i$
coincide, as both come from the left lcm of the
involved braid.}
\label{F:LCube1}
\end{figure}

We are ready to consider all possible triples of
letters. We sort them according to the numbers
of~$\ss{}$'s and~$\aa{}$'s. In the case of
three~$\ss¥$'s or of three $\aa¥$'s, the
condition is already known. So, we have only to
consider the four cases corresponding to
one~$\aa¥$ and two~$\ss¥$'s, or two~$\ss¥$'s and
one~$\aa¥$. The values follow from the formulas
of~\eqref{E:AuxRev0} and~\eqref{E:AuxRev1}.
Figure~\ref{F:LCube1} gives the details for the
$(\ss¥, \ss¥, \aa¥)$ case; the other three cases
are similar.

\subsection*{The right cube condition}

The verifications for the right cube condition are
similar, except that we use right reversing, \ie, we
push the negative letters to the right. Again, right
reversing leads to
at most one final word of the from~$\uu \vv\inv$ with
$\uu, \vv$ positive, but, in contrast to left
reversing, right reversing need not converge:
$\RRs$ contains no relation of the form $\aa i
\uu = \aa{i+1} \vv$ or $\ss i \uu = \aa i \vv$,
hence $\aa i\inv
\aa{i+1}$ and $\ss i\inv \aa i$ are not right
reversible.

\begin{figure}[h]
\begin{picture}(130,35)(0,-2)
\put(-1,-1){\includegraphics{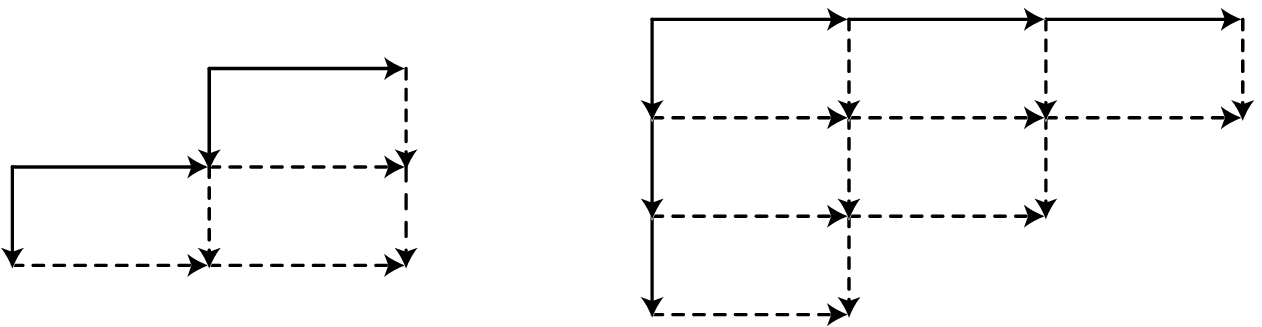}}
\put(-4,10){$\ss2$}
\put(8,17){$\ss1$}
\put(28,17){$\aa3$}
\put(16,20){$\ss1$}
\put(28,27){$\aa3$}
\put(6,2){$\ss1\ss2$}
\put(24,2){$\ss3\ss2\aa1$}
\put(21,10){$\ss2\ss1$}
\put(41,10){$\ss2\ss1$}
\put(41,20){$\ss1$}
\put(74,-3){$\e$}
\put(94,7){$\e$}
\put(114,17){$\e$}
\put(87,5){$\e$}
\put(107,15){$\e$}
\put(127,25){$\e$}
\put(61,25){$\ss2$}
\put(57,15){$\ss1\ss2$}
\put(54,5){$\ss3\ss2\aa1$}
\put(73,32){$\aa3$}
\put(93,32){$\ss1$}
\put(112,32){$\ss2\ss1$}
\put(86,25){$\ss2$}
\put(86,15){$\ss1\ss2$}
\put(106,25){$\ss2\ss1$}
\put(72,22){$\ss3\aa2$}
\put(70,12){$\ss3\ss2\aa1$}
\put(92,22){$\ss1\ss2$}
\put(73,25){$\revR$}
\put(93,25){$\revR$}
\put(115,25){$\revR$}
\put(73,15){$\revR$}
\put(95,15){$\revR$}
\put(73,5){$\revR$}
\put(8,10){$\revR$}
\put(30,10){$\revR$}
\put(28,20){$\revR$}
\end{picture}
\caption{\smaller Right cube
condition for the triple $(\ss2, \ss1, \aa3)$:
one first reverses $\ss 2\inv \ss1\ss1\inv \aa3$ to a
positive--negative word, here $\ss1\ss2\ss3\ss2\aa1
\ss1\inv \ss2\inv \ss1\inv$, and, then, one checks
that $\aa1\inv \ss2\inv \ss3\inv
\ss2\inv \ss1\inv \ss2\inv \aa3 \ss1 \ss2
\ss1$ is right reversible to~$\e$.}
\label{F:RightCube}
\end{figure}

It is possible to establish general formulas
similar to~\eqref{E:AuxRev0} and~\eqref{E:AuxRev1}.
Denote by $\uu \cR \vv$ and $\vv \cR \uu$ the unique
positive words such that $\uu\inv \vv$ is right
reversible to $(\uu \cR \vv) (\vv \cR \uu)\inv$, if
such words exist. Then, if $\uu, \vv$ are
$\ws$-words, $\uu \cR (\vv\aa j)$, when it
exists, is obtained from~$\uu \cR (\vv \ss j)$ by
replacing the final~$\ss k$ with the
corresponding~$\aa k$, and $\aa j \cR \uu$, when it
exists, is obtained from~$\uu$ by erasing the $j$-th
strand (in the braid diagram coded by~$\uu$).
However, such formulas are not very convenient as
they do not guarantee that the considered words 
exist, and it is actually easier to
systematically consider all possible cases, which
are not so many owing to symmetries and trivial
cases. Because of the above mentioned formula,
all words appearing have length~$6$ at most,
and the less trivial cases are when the indices
are neighbours. A typical example is given in
Figure~\ref{F:RightCube}; all other cases are
similar or more simple.

\section*{Index of terms and notation}

\makebox[\textwidth][s]{
\vtop{\parindent=0pt\hsize=7cm\smaller\smaller

\relax $\revr{}$ (word reversing): Def.~\ref{I:Reversing}

\relax $\LD\xx\yy$ (operation on~$\Bs$): Def.~\ref{I:Bracket}

\relax $\xx \OP \yy$ (operation on~$\Bs$):
Def.~\ref{I:Bracket}

\relax $\et_\xx$ (coloured tree): Def.~\ref{I:Skeleton}

\relax $\prec$ (tree ordering): Def.~\ref{I:TreeOrder}

\relax $\lF^{sp}$ (special order on~$F$):
Def.~\ref{I:ThompsonOrder}

\relax $\lp$ (positive ordering): Def.~\ref{I:PosOrder}

\relax $<$ (ordering): Def.~\ref{I:Order}

\relax $\aa*$: (family of all $\aa i$'s):
Def.~\ref{I:BBraids}

\relax address (node of a tree): Def.~\ref{I:Address}

\relax ALD-system (augmented LD-system):
Def.~\ref{I:LDSystem}

\relax $\wa$-word (parenthesized braid word):
Sec.~\ref{I:Diagrams}

\relax $\Bs$ (group of parenthesized braids):
Def.~\ref{I:Braids}

\relax $\BBs$ (group presented by~$\RRs$):
Def.~\ref{I:BBraids}

\relax $\BBsp$ (monoid presented by~$\RRs$):
Def.~\ref{I:Monoid}

\relax $\Bi^{sp}$ (special braids): Def.~\ref{I:Col}

\relax $c_n$ (right vine): Example~\ref{I:Vine}

\relax $\cc\tt\ttt$ (diagram completion):
Def.~\ref{I:Completion}

\relax $\Col(\ct)$ (colours in a tree): Def.~\ref{I:Col}

\relax $\Cone$ (positive cone): Def.~\ref{I:Cone}

\relax coherent ($\FGs$-coloured tree): Def.~\ref{I:Coherent}

\relax complete (presentation): Def.~\ref{I:Complete}

\relax completion (of a braid diagram):
Def.~\ref{I:Completion}

\relax $\partial$ (shift mapping): Def.~\ref{I:Shift}

\relax $\DD_\tt(w)$ (braid diagram):
Def.~\ref{I:Diagram}

\relax $\dbl_k(\ww)$ (strand doubling):
Def.~\ref{I:Double}

\relax $\Dec(t)$ (tree decomposition):
Def.~\ref{I:RD}

\relax $D_s$ (disk): Def.~\ref{I:Disk}

\relax $\Dya(\tt)$ (rationals associated with~$\tt$):
Def.~\ref{I:Tree}

\relax dyadic realization (sequence):
Def.~\ref{I:Position}

\relax equivalent (braid diagrams):
Def.~\ref{I:Equivalence}

\relax $\ev(\ct)$ (evaluation of a coloured tree):
Def.~\ref{I:Evaluation}

\relax $\Evv(\ct)$ (evaluation of a coloured tree):
Def.~\ref{I:Evaluation}

\vfill}\hfill\vtop{\parindent=0pt\hsize=7cm\smaller\smaller

\relax $\phi$ (morphism of $\Bs$ to $\MCG(\SK)$):
Def.~\ref{I:Disk}

\relax $\hom\fx$ (homeomorphism): Prop.~\ref{I:Homeo}

\relax homogeneous (presentation): Def.~\ref{I:Homogeneous}

\relax $\Can$ (Cantor set): Def.~\ref{I:Sphere}

\relax LD-system: Def.~\ref{I:LDSystem}

\relax $\NNNs$ (set of all positions):
Def.~\ref{I:Position}

\relax naturel (colouring): Def.~\ref{I:Address}

\relax parenthesized braid: Def.~\ref{I:Braids}

\relax $\Pos\tt$ (positions associated with~$\tt$):
Def.~\ref{I:Tree}

\relax position: Def.~\ref{I:Position}

\relax $\psi$ (morphism of~$\Bs$ to $\Aut(\FGs)$):
Prop.~\ref{I:Psi}

\relax $\rho$ (construction of a Cantor set):
Def.~\ref{I:Sphere}

\relax $\RRs$ (relations): Def.~\ref{I:BBraids}

\relax rack: Def.~\ref{I:LDSystem}

\relax $\red(u)$ (free reduced word): Def.~\ref{I:Red}

\relax reversing: Def.~\ref{I:Reversing}

\relax $\dd s$ (dyadic realization):
Def.~\ref{I:Position}

\relax special (parenthesized braid): Def.~\ref{I:Special}

\relax special (word): Def.~\ref{I:SpecialBis}

\relax $\ss*$ (family of all $\ss i$'s): Def.~\ref{I:BBraids}

\relax $\SK$ (sphere with a Cantor removed):
Def.~\ref{I:Sphere}

\relax $\ws$-word, $\wsa$-word (parenthesized braid word):
Sec.~\ref{I:Diagrams}

\relax $\ss i$-positive (word): Def.~\ref{I:Positive}

\relax skeleton (tree): Def.~\ref{I:Skeleton}

\relax $\tt^\dag$ (skeleton of~$\tt$): Def.~\ref{I:Skeleton}

\relax $\ct \act \ww$ (action on $\Bs$-coloured tree):
Def.~\ref{I:Action} 

\relax $\ct \act \ww$ (action on $\FGs$-coloured tree):
Def.~\ref{I:ActionBis} 

\relax tidy (word): Def.~\ref{I:Tidy}

\relax $\ww[k]$ (initial position):
Def.~\ref{I:Double}

\relax $\cl\ww$ (element represented by~$\ww$): 
Lemma~\ref{L:MainEval}

\relax $\widehat\ww$ (automorphism of~$\FGs$):
Prop.~\ref{I:Hat}

\relax $x_s$ (loop class): Def.~\ref{I:Loop}

\vfill}}

\end{document}